\documentclass{cmslatex}
\usepackage[paperwidth=7in, paperheight=10in, margin=.875in]{geometry}
\usepackage{amsfonts,amssymb}
\usepackage{amsmath}
\usepackage{cases}
\usepackage{graphicx}
\usepackage{cite}
\usepackage{enumerate}
\usepackage{bm}
\usepackage[backref,colorlinks,linkcolor=red,anchorcolor=green,citecolor=blue]{hyperref}
\usepackage[capitalize,nameinlink]{cleveref}
\allowdisplaybreaks[4]
\renewcommand{\theequation}{\arabic{section}.\arabic{equation}}

\newcommand{\inner}[2]{( #1, #2 ) }
\newcommand{\biginner}[2]{\big( #1, #2 \big) }
\newcommand{\Biginner}[2]{\Big( #1, #2 \Big) }
\newcommand{\pairing}[2]{\left\langle #1, #2 \right\rangle }
\newcommand{\pdif}[2]{ \mathrm{\partial}_{#2} #1}
\newcommand{\dif}[2]{ \frac{ \mathrm{d} }{ \mathrm{d #2} } #1 }
\newcommand{\intomega}[3]{\int_{#3} #1 \ \mathrm{d}#2}
\newcommand{\intinterval}[4]{\int_{#3}^{#4} #1 \ \mathrm{d}#2}
\newcommand{\norm}[1]{\| #1 \|}
\newcommand{\bignorm}[1]{\big\| #1 \big\|}
\newcommand{\Bignorm}[1]{\Big\| #1 \Big\|}
\newcommand{\wkp}[2]{W^{#1,#2}(\Omega)}
\newcommand{\wkpzero}[2]{W^{#1,#2}_0(\Omega)}
\newcommand{\wkpsigma}[2]{\bm{W}^{#1,#2}_{\sigma}(\Omega)}
\newcommand{\hs}[1]{H^{#1}(\Omega)}
\newcommand{\hszero}[1]{H^{#1}_0 (\Omega)}
\newcommand{\hssigma}[1]{\bm{H}^{#1}_{\sigma}(\Omega)}
\newcommand{\hst}[2]{H^{#1} \bracket{#2}}

\newcommand{\hstuloc}[2]{H^{#1}_{\mathrm{uloc}} \bracket{#2}}

\newcommand{\lpsigma}[1]{\bm{L}^{#1}_{\sigma}(\Omega)}
\newcommand{\lpt}[2]{L^{#1} \bracket{#2}}

\newcommand{\lptbig}[2]{L^{#1} \bigbracket{#2}}
\newcommand{\lptbigloc}[2]{L^{#1}_{\mathrm{loc}} \bigbracket{#2}}
\newcommand{\lptbiguloc}[2]{L^{#1}_{\mathrm{uloc}} \bigbracket{#2}}
\newcommand{\wkpt}[3]{W^{#1,#2} \bracket{#3}}

\newcommand{\wkptbig}[3]{W^{#1,#2} \bigbracket{#3}}
\newcommand{\wkptbigloc}[3]{W^{#1,#2}_{\mathrm{loc}} \bigbracket{#3}}
\newcommand{\wkptbiguloc}[3]{W^{#1,#2}_{\mathrm{uloc}} \bigbracket{#3}}
\newcommand{\rd}[1]{\mathbb{R}^{#1}}
\newcommand{\bracket}[1]{( #1 )}
\newcommand{\bigbracket}[1]{\big( #1 \big)}
\newcommand{\Bigbracket}[1]{\Big( #1 \Big)}
\newcommand{\halfof}[1]{\frac{#1}{2}}
\newcommand{\limit}[1]{\mathop{\mathrm{lim}}\limits_{#1}}

\newcommand{\abs}[1]{| #1 |}
\newcommand{\bigabs}[1]{\big| #1 \big|}
\newcommand{\Bigabs}[1]{\Big| #1 \Big|}
\newcommand{\lp}[1]{L^{#1}(\Omega)}
\newcommand{\set}[1]{\{ #1 \} }

\def\be{\begin{equation}}
\def\ee{\end{equation}}
\def\bea{\begin{eqnarray}}
\def\eea{\end{eqnarray}}
\def\bt{\begin{theorem}}
\def\et{\end{theorem}}
\def\bl{\begin{lemma}}
\def\el{\end{lemma}}
\def\bnum{\begin{numcases}{}}
\def\enum{\end{numcases}}
\def\br{\begin{remark}}
\def\er{\end{remark}}
\def\bpf{\begin{proof}}
\def\epf{\end{proof}}
\def\bpp{\begin{proposition}}
\def\epp{\end{proposition}}
\def\bc{\begin{corollary}}
\def\ec{\end{corollary}}
\def\bd{\begin{definition}}
\def\ed{\end{definition}}
\def\ba{\begin{assumption}}
\def\ea{\end{assumption}}
\def\p{\mathrm{\partial}}

\def\ltwo{L^2(\Omega)}
\def\hone{H^1(\Omega)}

\def\ls{\bm{L}^2_{\sigma}(\Omega)}

\def\sp{\mathbf{u}}
\def\half{\frac{1}{2}}
\def\div{\mathrm{div}\,}
\def\nv{\mathbf{n}}
\def\mhi{\mu_{\phi}}
\def\msi{\mu_{\psi}}
\def\ul{_{\mathrm{uloc}}}
\def\loc{_{\mathrm{loc}}}
\def\bmtheta{\bm{\theta}}
\def\bmeta{\bm{\eta}}
\def\eto{E_{\mathrm{tot}}}
\def\efr{E_{\mathrm{free}}}
\def\ehi{E_{\phi}}
\def\esi{E_{\psi}}

\def\ae{\mathrm{a.e.}}
\def\bmj{\mathbf{J}}
\def\gradphi{\nabla \phi}
\def\gradpsi{\nabla \psi}
\def\gradmhi{\nabla \mhi}
\def\gradmsi{\nabla \msi}
\def\mbhi{m_{\phi}}
\def\mbsi{m_{\psi}}

\def\pn{\p^-_{t,h}}

\def\r{\mathbb{R}}

\def\tphi{\widetilde{\phi}}
\def\tpsi{\widetilde{\psi}}

\def\psigma{\mathbb{P}_{\sigma}}
\def\jhi{\mathbf{J}_{\phi}}
\def\jsi{\mathbf{J}_{\psi}}
\def\whi{W_{\phi}}
\def\wsi{W_{\psi}}
\def\fhi{F_{\phi}}
\def\fsi{F_{\psi}}
\def\flnhi{F_{\mathrm{ln},\phi}}
\def\flnsi{F_{\mathrm{ln},\psi}}
\def\tfhi{\widetilde{F}_{\phi}}
\def\tfsi{\widetilde{F}_{\psi}}
\allowdisplaybreaks
\begin{document}
	
\title{Global Weak Solutions of a Thermodynamically Consistent Diffuse Interface Model for Nonhomogeneous Incompressible Two-phase Flows with a Soluble Surfactant
\thanks{Received date, and accepted date}}

\author{
Bohan Ouyang
\thanks{School of Mathematical Sciences,
Fudan University, Handan Road 220, 200433 Shanghai, P.R. China, Email: \texttt{22110180033@m.fudan.edu.cn}.}
\and
Maurizio Grasselli
\thanks{Dipartimento di Matematica,
Politecnico di Milano, Via E. Bonardi 9, 20133 Milano, Italy, Email: \texttt{maurizio.grasselli@polimi.it}.}
\and
Hao Wu
\thanks{Corresponding author. School of Mathematical Sciences
and Shanghai Key Laboratory for Contemporary Applied Mathematics, Fudan University, Handan Road 220, 200433 Shanghai, P.R. China, Email: \texttt{haowufd@fudan.edu.cn}.}
}
\pagestyle{myheadings} \markboth{INCOMPRESSIBLE TWO-PHASE FLOWS WITH A SOLUBLE SURFACTANT}{B.-H. OUYANG, M. GRASSELLI, H. WU}
\maketitle

\begin{abstract}
We study a thermodynamically consistent diffuse interface model that describes the motion of a two-phase flow of two viscous incompressible Newtonian fluids with unmatched densities and a soluble surfactant in a bounded domain of two or three dimensions. The resulting hydrodynamic system consists of a nonhomogeneous Navier--Stokes system for the (volume averaged) velocity $\sp$ and a coupled Cahn--Hilliard system for the phase-field variables $\phi$ and $\psi$ that represent the difference in volume fractions of the binary fluids and the surfactant concentration, respectively. For the initial boundary value problem with physically relevant singular potentials subject to a no-slip boundary condition for the fluid velocity and homogeneous Neumann boundary conditions for the phase-field variables and the chemical potentials, we first establish the existence of global weak solutions in the case of non-degenerate mobilities based on a suitable semi-implicit time discretization. Next, we prove the existence of global weak solutions for a class of general degenerate mobilities, with the aid of a new type of approximations for both the mobilities and the singular parts of the potential densities.
\end{abstract}



\section{Introduction}\label{intro}

Surfactants are ubiquitous in real-world applications appearing in nature and industry. There are various types of surfactants, such as wetting agents and emulsifiers, that are important in oil recovery, food production, and microfluidics \cite{RK2012}. A classic example is oil-water systems with detergents. Surfactants such as detergents preferentially adsorb at the free interfaces between oil and water. This selective adsorption is due to the amphiphilic structure of surfactant molecules: the hydrophilic heads interact with water, while the hydrophobic tails align with oil, anchoring the molecules at the interface. Surfactant adsorption reduces the interfacial tension between fluid components, thus facilitating their mixing, affecting wettability in coatings, and stabilizing emulsions/foams by inhibiting droplet/bubble coalescence \cite{Chang1995,MS2020,Petkova2012,Weinstein2024}. On the other hand, spatial variations in surfactant concentration can induce surface tension gradients that generate Marangoni forces. These forces, combined with capillary effects, influence interfacial stability and flow dynamics in fluid mixtures, playing a fundamental role in many hydrodynamic processes \cite{Anna2016}.

Modeling interfacial dynamics of fluid mixtures with surfactants is a challenging task due to the multiscale nature of the problem and the complex coupling to hydrodynamics \cite{MS2020}. In addition to the hydrodynamic system for the mixed fluids, an advection-diffusion equation should be introduced to describe the evolution of the surfactant concentration. For soluble surfactants, the adsorption and desorption kinetics at the interface have to be taken into account \cite{DA1996,CMS1997}. Numerous theoretical as well as numerical studies have been dedicated to modeling and analysis of two-phase flows with surfactants \cite{AGLW2017,AGW2019,GLS2014,JL2004,LTH2008,MT2008,XLLZ2006}. Among various approaches in the literature, the diffuse interface model turns out to be an efficient tool \cite{AGLW2017,GLS2014,EDAT2013,KK1997,La1992,LK2012,LZ2010,TSLV2011,vander2006}. In the diffuse interface framework, the two fluid phases are treated as macroscopically immiscible but allow partial mixing in a narrow interfacial layer. An order parameter (i.e., a phase-field function) is introduced to distinguish the bulk fluids and capture the transition between them without explicit tracking of the interface \cite{AMW1998}, that is, it takes distinct (almost) constant values in each of the fluid components and varies smoothly across the interfacial layer. Compared to the classical sharp interface model, in which the interfaces separating fluid components are modeled with evolving hypersurfaces, the diffuse interface model can naturally handle topological changes of free interfaces, for instance, breakup, reconnection, and tip-streaming driven by Marangoni forces \cite{AMW1998,LT1998}.

In this work, we analyze a thermodynamically consistent diffuse interface model for a two-phase flow of two viscous incompressible Newtonian fluids with a soluble surfactant in the general scenario of unmatched densities and variable mobilities. The model is derived through the first law of thermodynamics, associated thermodynamic relations, and the Onsager variational principle. For details, we refer to \cite{EDAT2013,ZKYWYS2019} and the references therein. The resulting hydrodynamic system can be regarded as a generalization of the non-homogeneous Navier--Stokes--Cahn--Hilliard system proposed in \cite{AGG2012}, incorporating additional interactions with a soluble surfactant and a suitable mass source. More precisely, given a bounded domain $\Omega$ with a smooth boundary $\p \Omega$ in $\rd{d}$, $d \in \{ 2,3\}$, we consider the following coupled system:
\begin{numcases}{}
    \pdif{(\rho (\phi) \sp)}{t} + \div\bracket{\sp \otimes (\rho (\phi) \sp + \mathbf{J}) } - \div (\nu (\phi) D\sp)+ \nabla\pi \nonumber\\
    \quad = \mu_{\phi} \nabla \phi + \mu_{\psi} \nabla \psi \label{eq:nsch1},
    \\
    \div \, \sp =0 \label{eq:nsch2},
    \\
    \pdif{\phi}{t} + \sp \cdot \nabla \phi +  \sigma_1 (\phi)\bigbracket{\overline{\phi} -c}  = \div (m_{\phi} (\phi) \nabla \mhi) \label{eq:nsch3},
    \\
    \mu_{\phi} = -\Delta \phi + \sigma_2 \mathcal{N} \bigbracket{\phi - \overline{\phi}} + F_{\phi}'(\phi) + \pdif{G}{\phi}(\phi,\psi) \label{eq:nsch4},
    \\
    \pdif{\psi}{t} + \sp \cdot \nabla \psi = \div (m_{\psi} (\psi) \nabla \msi) \label{eq:nsch5},
    \\
    \mu_{\psi} = -\beta \Delta \psi + F_{\psi}'(\psi) + \pdif{G}{\psi}(\phi,\psi) \label{eq:nsch6},
\end{numcases}
in $\Omega \times (0,\infty)$, subject to the boundary and initial conditions:
\begin{numcases}{}
    \sp = \mathbf{0}, \ \pdif{\phi}{\nv} = m_\phi(\phi)\pdif{\mu_{\phi}}{\nv} = \pdif{\psi}{\nv} = m_\psi(\psi)\pdif{\mu_{\psi}}{\nv} = 0, \label{eq:nschb}& on $\p \Omega \times (0,\infty)$,
    \\
    \sp|_{t=0} = \sp_0(x), \ \phi|_{t=0} = \phi_0(x), \ \psi|_{t=0} = \psi_0(x), \label{eq:nschi}& in $\Omega$.
\end{numcases}
Here, $\mathbf{n}=\mathbf{n}(x)$ denotes the unit outward normal vector on $\partial\Omega$ and $\partial_\mathbf{n}$ stands for the outward normal derivative on the boundary. The state variables of the system include the volume averaged (solenoidal) velocity $\sp: \Omega \times [0,\infty)\to \mathbb{R}^d$, the pressure of the fluid mixture $\pi: \Omega\times[0,\infty)\to \mathbb{R}$, and the phase-field variable $\phi: \Omega\times[0,\infty)\to [-1,1]$ that denotes the volume fraction difference between the two fluid components. As in \cite{La1992,KK1997,EDAT2013}, an additional phase-field variable $\psi: \Omega\times[0,\infty)\to [0,1]$ is introduced to denote the volume fraction of the surfactant. Unlike the classical multi-component system (see, e.g., \cite{Dong2014,EG1997}), the concentration of surfactant is usually very small, that is, $0<\psi\ll 1$. Thus, the (three) concentration components of the binary fluids and the surfactant should be understood as $(1-\psi)(1+\phi)/2, (1-\psi)(1-\phi)/2, \psi$, and then approximated with $(1+\phi)/2, (1-\phi)/2, \psi$, see \cite{EDAT2013}.

The equations \eqref{eq:nsch1}--\eqref{eq:nsch2} correspond to a nonhomogeneous Navier--Stokes system for the fluid velocity $\sp$, while \eqref{eq:nsch3}--\eqref{eq:nsch6} yield a coupled convective Cahn--Hilliard system for the phase-field variables $(\phi,\psi)$. In \eqref{eq:nsch1}, $D\sp =\frac{1}{2}(\nabla \sp+\nabla \sp^T)$ denotes the symmetrized gradient of $\sp$. Let $\rho_1$, $\rho_2$, and $\nu_1$, $\nu_2$ be homogeneous positive density and viscosity coefficients of the unmixed fluid components. Typical forms of the average density $\rho$ and the average viscosity $\nu$ are given by
\begin{align}
\rho(\phi)= \frac{\rho_1-\rho_2}{2}\phi+\frac{\rho_1+\rho_2}{2},
\quad \nu(\phi)= \frac{\nu_1-\nu_2}{2}\phi+\frac{\nu_1+\nu_2}{2}.
\notag
\end{align}
We assume for simplicity that the density and viscosity of the fluid mixture depend only on $\phi$ since $\psi$ is small \cite{ZKYWYS2019}. The scalar functions $\mhi:\Omega\times [0,\infty)\to \mathbb{R}$ and $\msi:\Omega\times [0,\infty)\to \mathbb{R}$ are the so-called chemical potentials for the two-phase flow and the surfactant, respectively. These chemical potentials can be obtained by taking variational derivatives of the free energy functional $\efr$ defined in \eqref{free} below. Next, the nonnegative scalar functions $m_\phi:[-1,1]\to [0,\infty)$ and $m_\psi:[0,1]\to [0,\infty)$ are mobility coefficients for the binary fluids and the surfactant, respectively. In general, they can depend on $\phi$ and $\psi$ and can even be degenerate in pure phases \cite{EDAT2013,Novick}. Finally, the relative flux
$$
\mathbf{J} = - \halfof{\rho_1 - \rho_2} m_{\phi}(\phi) \nabla \mu_{\phi}
$$
is related to diffusion of the fluid components. It only appears in the case of unmatched densities (i.e., $\rho_1\neq \rho_2$) and is crucial for consistency with thermodynamics \cite{AGG2012}.

The total free energy of the system \eqref{eq:nsch1}--\eqref{eq:nsch6} consists of three parts (see e.g., \cite{EDAT2013,KK1997,La1992}):
\begin{align}
 \label{free}
    &\efr \left(\phi,\psi\right) = \ehi(\phi) + \esi(\psi) + \int_{\Omega}G\left( \phi,\psi \right) \mathrm{d}x.
\end{align}
The first two terms on the right-hand side of \eqref{free} are the free energies associated with $\phi$ and $\psi$, respectively, while $G$ is the energy density function for interactions between the fluid mixture and the surfactant. Setting some constants (e.g., the interfacial thickness) equal to unity for simplicity, we take the following Ginzburg--Landau type free energies
 \begin{align*}
    & \ehi( \phi ) =  \intomega{\left(\half\vert\nabla \phi\vert^2 + \fhi(\phi) + \halfof{\sigma_2} \vert \nabla \mathcal{N}\bigbracket{\phi - \overline{\phi}}\vert^2\right)}{x}{\Omega},\\
    & \esi( \psi ) =  \intomega{\left(\halfof{\beta}\vert\nabla \psi\vert^2 + \fsi(\psi)\right)}{x}{\Omega},
\end{align*}
where $\sigma_2 \in \mathbb{R}$, $\beta>0$ are given constants. The nonlocal operator $\mathcal{N}$ stands for the inverse of the minus Laplace operator with homogeneous Neumann boundary condition and $\overline{\phi}$ denotes the spatial average of $\phi$ over $\Omega$, that is, $\overline{\phi}=|\Omega|^{-1}\int_\Omega \phi\,\mathrm{d}x$. When $\sigma_2\neq 0$, $\ehi$ yields an Ohta--Kawasaki type functional accounting for possible nonlocal interactions between the two fluid components (see, e.g.,\cite{CR2003,NO1995}). The parameter $\sigma_2$ measures the strength of nonlocal interaction, and its sign indicates suppression or enhancement of the chemical reaction for the phase separation and coarsening process \cite{CMW2011}. In $\ehi$ and $\esi$, the homogeneous free energy densities for the fluid mixture (i.e., $\fhi$) and the surfactant (i.e., $\fsi$) take the Flory--Huggins type potentials, namely,
\begin{align}
\label{FH3}
    & \fhi(s) = \halfof{\theta_1}  [(1+s) \ln(1+s) + (1-s) \ln(1-s)] + \halfof{\widetilde{\theta}_1} \bigbracket{1-s^2}, \quad s\in(-1,1),\\
\label{FH1}
    & \fsi(s) = \halfof{\theta_2}  [s \ln s + (1-s) \ln(1-s)] + \halfof{\widetilde{\theta}_2} s(1-s), \quad s\in(0,1),
\end{align}
with $\theta_1, \widetilde{\theta}_1,\theta_2,\widetilde{\theta_2}>0 $ being temperature-relevant constants.
The parameters $\widetilde{\theta}_1$, $\widetilde{\theta}_2$ measure the strength of the demixing effect. For sufficiently large $\widetilde{\theta}_1>\theta_1$, $\widetilde{\theta}_2>\theta_2$, both potentials $\fhi$ and $\fsi$ present a double-well structure, which leads to possible phase separation phenomena. We recall that in previous work \cite{EDAT2013,KK1997,La1992,vander2006,ZKYWYS2019}, $\fhi$ was taken as a fourth-order polynomial with double-well structure, that is, the classical approximation of the Flory--Huggins potential, while a logarithmic potential $\fsi$ was preferred for the surfactant, as it provides an isotherm relation. The choice of a regular potential in the fourth Cahn--Hilliard equation leads to easier treatments in both theoretic and numerical analysis; however, it cannot ensure that $\phi$ takes its values in the physical range $[-1,1]$, see \cite{M2019}. In this study, we consider physically relevant singular potentials for both $\phi$ and $\psi$, i.e., \eqref{FH3}, \eqref{FH1}. The singular nature of the potential $\fhi$ (resp. $\fsi$) guarantees that the phase-field variable $\phi$ (resp. $\psi$) always takes its values in the physical interval $[-1,1]$ (resp. $[0,1]$) as time evolves. In addition, we note that the competition between the gradient term (with preference to homogeneous mixing) and the double-well potential (with preference to unmixed phases) in the free energies can create complex pattern formation in equilibrium \cite{AMW1998,LT1998,Novick}.

Next, let us say some words about the bi-variate function $G$ in \eqref{free}, whose typical form is as follows:
\begin{equation}
    G(\phi,\psi) = \halfof{\gamma_1} \psi \phi^2 - \frac{\gamma_2}{4} \psi \bracket{1-\phi^2}^2,
    \label{ab-G}
\end{equation}
with $\gamma_1, \gamma_2>0$ being given constants. The two terms on the right-hand side of \eqref{ab-G} are complementary in the sense that the first term penalizes the presence of free surfactant in the respective bulk phases (that is, an enthalpic term measuring the cost of free surfactant and $\gamma_1$ is related to the bulk solubility), while the second term favors the surfactant to reside at the free interfaces between the two fluid components (see ``Model 3'' in \cite[Section 3.2.2]{EDAT2013}). The second term serves as a gradient-free approximation for the surface energy potential $-\gamma_3 \psi \abs{\nabla \phi}^2$ with $\gamma_3>0$ (see \cite{La1992,KK1997}), which accounts for the adsorption of the surfactant to the fluid interfaces. Here, the square gradient acts as a natural diffuse version of the sharp interface indicator function. It corresponds to the so-called model ``Model 1'' in \cite[Section 3.2.1]{EDAT2013}, in which an additional regularization term $\halfof{\beta}\vert\nabla \psi\vert^2$ with $\beta>0$ was introduced in the associated free energy $\esi$ to guarantee well-posedness. Another possible approximation of the gradient adsorption term could be $-\frac{\gamma_2}{4} \psi \bracket{1-\phi^2}$, see ``Model 2'' in \cite[Section 3.2.2]{EDAT2013}. Nevertheless, as suggested in \cite{EDAT2013}, the choice \eqref{ab-G} performs better in numerical experiments compared to the other two formulae. In the current framework, the adsorption and desorption kinetics of the soluble surfactant at the fluid interface are characterized by the interaction energy $\int_\Omega G(\phi,\psi) \,\mathrm{d}x$ that yields a single Cahn--Hilliard equation for the surfactant phase-field variable $\psi$. Alternative descriptions were given in \cite{AGLW2017,GLS2014} using three evolution equations for the bulk surfactant density in each fluid phase and the interface surfactant density, respectively.

In the system \eqref{eq:nsch1}--\eqref{eq:nsch6}, we also allow certain (reversible) chemical reactions between the fluid components \cite{Huo03,Huo04}. This corresponds to the inclusion of a mass source term of Oono's type, see \cite{GGM2017} for the single Cahn--Hilliard--Oono equation and \cite{DG2022,O2024} for the coupled Cahn--Hilliard system related to \eqref{eq:nsch3}--\eqref{eq:nsch6} without fluid interaction, see also \cite{BGM2014,MT2016} for the Navier--Stokes--Cahn--Hilliard--Oono system. We consider in this study a slightly more general case with possibly non-constant reaction coefficients, e.g.,  $\sigma_1(s)\ge0$ for all $s \in \r$. If $\sigma_1 =\sigma_2=\sigma\in [0,\infty)$ and $\mbhi = 1$, then it holds
$$
\sigma_1 \bigbracket{\overline{\phi} -c} - \mathrm{div}(m_\phi\nabla (\sigma_2\mathcal{N}(\phi-\overline{\phi})))=\sigma (\phi-c),
$$
so that the standard linear Oono term can be recovered ($c\in (-1,1)$ is a given constant). Inspired by \cite{HW2023}, here we have decomposed the standard Oono's interaction term into two parts: the term $\sigma_2\mathcal{N}(\phi-\overline{\phi})$ contributes not only to the free energy $\ehi$, but also to the capillary force at fluid-fluid interfaces (see \eqref{eq:nsch1}), while $\sigma_1 \bigbracket{\overline{\phi} -c}$ yields a possible mass transfer between two components of the fluid mixture and an additional contribution in the energy production (see \eqref{mass-phi}, \eqref{BEL-tot} below). In general, the analysis of Cahn--Hilliard type equations with a mass source term is a non-trivial task, since mass transfer can influence both the mass dynamics and the energy balance, see \cite{CFG2024,GLRS2022,JWZ2015,M2019} and the references therein.

The aim of this study is to establish the existence of global weak solutions to the initial boundary value problem \eqref{eq:nsch1}--\eqref{eq:nschi} in both cases of non-degenerate mobilities and degenerate mobilities (see Theorem \ref{wsnschndme} and Theorem \ref{wsnschdme}). Here, we choose to start with the theoretical analysis of the system \eqref{eq:nsch1}--\eqref{eq:nsch6} by considering the standard boundary conditions \eqref{eq:nschb} for the sake of simplicity. In recent work \cite{WLZ2024,ZKYWYS2019}, the generalized Navier boundary condition for $\sp$ and a dynamic boundary condition for $\phi$ were introduced to account for possible contact line dynamics at the solid boundary, which is important in practical applications. The associated theoretical analysis turns out to be more challenging (see \cite{GGM16,GGW2019,GK2023} and the references therein for related Navier--Stokes--Cahn--Hilliard systems on two-phase flows without surfactant interaction). This will be the goal of future studies.

Our work extends the previous results in \cite{ADG2013ndm, ADG2013dm} on the existence of global weak solutions to the Abels--Garcke--Gr\"{u}n (AGG) model \cite{AGG2012} for nonhomogeneous two-phase flows of incompressible viscous fluids without the influence of surfactants. We refer to \cite{AGG2024,G2021,G2022} for recent progress on the AGG model concerning the existence and uniqueness of local/global strong solutions, propagation of regularity for global weak solutions, and long-time behavior. See also \cite{AGP2024} for a multi-component extension of the AGG model, and \cite{GGW2019,GLW2024,GK2023} for analysis of extended AGG models with dynamic boundary conditions.

The proofs of Theorem \ref{wsnschndme} and Theorem \ref{wsnschdme} rely on two fundamental properties of the coupled system \eqref{eq:nsch1}--\eqref{eq:nsch6} subject to the boundary conditions in \eqref{eq:nschb}, that is, \textit{mass conservation} and \textit{energy balance}. More precisely, for sufficiently regular solutions to problem \eqref{eq:nsch1}--\eqref{eq:nschi}, integrating \eqref{eq:nsch3} and \eqref{eq:nsch5} over $\Omega$ respectively, we find
\begin{align}
&\frac{\mathrm{d}}{\mathrm{d}t}\int_\Omega \phi\,\mathrm{d}x + \bigbracket{\overline{\phi} -c}\int_\Omega \sigma_1(\phi)\,\mathrm{d}x=0,
\quad &\forall\, t>0,
\label{mass-phi}
\\
&\frac{\mathrm{d}}{\mathrm{d}t}\int_\Omega \psi\,\mathrm{d}x=0,\quad &\forall\, t>0.
\label{mass-psi}
\end{align}
On the other hand, the total energy of the problem \eqref{eq:nsch1}--\eqref{eq:nschi} is defined as
\begin{equation*}
    \eto (\sp,\phi,\psi) = \intomega{\halfof{\rho(\phi)} \abs{\sp}^2}{x}{\Omega} + \efr(\phi,\psi).
\end{equation*}
Then a straightforward calculation leads to the following energy identity:
\begin{align}
    & \frac{\mathrm{d}}{\mathrm{d}t} \eto (\sp,\phi,\psi) + \intomega{ \bracket{ \nu(\phi) \abs{D \sp}^2 + \mbhi(\phi) \abs{\nabla \mhi}^2 + \mbsi(\psi) \abs{\nabla \msi}^2 } }{x}{\Omega}
    \notag \\
    & \quad +  \bigbracket{\overline{\phi} -c} \intomega{\sigma_1(\phi) \Bigbracket{\mhi - \frac{\rho_1 - \rho_2}{4} \abs{\sp}^2 }}{x}{\Omega}  = 0,\qquad \forall\, t>0.
    \label{BEL-tot}
\end{align}

When the mobilities $\mbhi$ and $\mbsi$ are both non-degenerate according to hypothesis $(\mathbf{H3})$ (see Section \ref{sectionassumptions} below), we prove the existence of a global weak solution (see Theorem \ref{wsnschndme}) by adapting the strategy in \cite{ADG2013ndm}, namely, with the aid of a semi-implicit time discretization. However, additional efforts are required to handle interactions between the fluid mixture and the surfactant, as well as the influence due to the mass source. To construct a time discretization that preserves a discrete counterpart of the energy balance \eqref{BEL-tot}, we need to approximate the first-order derivatives of the interaction energy density $G$ by suitable difference quotients. Taking into account the mass source term $\sigma_1(\phi) \bigbracket{\overline{\phi} -c}$, although the uniform boundedness in time of global weak solutions holds unconditionally in the case of a nonnegative constant $\sigma_1$, the same property may not be valid in the general case with a non-constant $\sigma_1$. To overcome this difficulty, we impose some reasonable requirements regarding the integrability and decay properties of $\sigma_1(\phi) \bigbracket{\overline{\phi} -c}$ (see \eqref{additional_requirements_uniform_boundedness} and Remark \ref{add-sigphi}). Besides, to handle the extra term
\begin{equation}\label{additional_term_sp}
    \frac{\rho_1-\rho_2}{4} \bigbracket{\overline{\phi}-c} \intomega{ \sigma_1 \bracket{\phi} \abs{\sp}^2 }{x}{\Omega}
\end{equation}
on the right-hand side of the energy identity \eqref{BEL-tot}, we estimate the solution on a finite interval $[0,T]$ and the unbounded interval $[T,\infty)$ separately, with $T>0$ being sufficiently large but fixed. In this way, we are able to derive estimates for weak solutions on $[0,T]$ by Gronwall's lemma, while on $[T,\infty)$, we can control \eqref{additional_term_sp} by utilizing the decay property of $\sigma_1(\phi) \bigbracket{\overline{\phi} -c}$ (i.e., a consequence of \eqref{mass-phi}), together with the energy dissipation $ \intomega{ \nu\bracket{\phi} \abs{D \sp}^2 }{x}{\Omega}$ and Korn's inequality.

When the mobilities $\mbhi$ and $\mbsi$ are both degenerate at the pure phases $\phi=\pm 1$ and $\psi=0,1$, respectively (see hypothesis $(\mathbf{H3*})$ below), we prove the existence of a suitably defined global weak solution (see Theorem \ref{wsnschdme}) for a class of degenerate mobilities (cf. \cite{EG1996,B1999}) that are more general than those considered in \cite{ADG2013dm}, together with compatible singular potentials. More precisely, under our assumptions, we can handle mobility functions such as $\mbhi(s) = (1-s^2)^k$, $s \in [-1,1]$, as well as $\mbsi(s) = s^k(1-s)^k$, $s\in [0,1]$, for all $k \in [1,\infty)$, while in \cite{ADG2013dm} only mobility with linear degeneracy (and regular potentials), that is, $\mbhi(s)=1-s^2$, was considered. Taking advantage of the existence result in the non-degenerate case (i.e., Theorem \ref{wsnschndme}), we achieve our goal by extending the techniques developed in \cite{ADG2013dm}. The key point is to introduce suitable approximations of the potential functions $\fhi$, $\fsi$ to ensure that some of our structural assumptions can be inherited by the approximate problem. This allows us to apply Lebesgue's dominated convergence theorem to guarantee the convergence of some tricky terms in the weak formulation. More precisely, we regularize the degenerate mobilities by non-degenerate ones and approximate the original singular potentials by a family of regular potentials plus a small perturbation in the form of a specific given singular potential multiplied by a small approximating parameter. In this way, Theorem \ref{wsnschndme} can be applied to obtain the existence of global weak solutions to the approximate problem. To derive estimates for approximate solutions that are uniform with respect to the approximating parameter, we need an additional structural assumption on the coefficient $\sigma_1$. Thanks to our specific approximations for the mobilities as well as the potentials, we can obtain the required uniform estimates for approximate solutions by combining energy estimates and the so-called entropy estimates (see Lemma \ref{energyestimateapdm}).
Our work provides an alternative approach to construct global weak solutions for the Navier--Stokes--Cahn--Hilliard type systems with general degenerate mobilities, compared to the arguments in \cite{B1999} that require certain monotonicity property of the second order derivative for the singular potential. The same method with minor modification can be applied to the case where only one of the mobility functions is degenerate.

It is worth mentioning that ``Model 1'' proposed in \cite[Section 3.2.1]{EDAT2013} has been extensively analyzed in \cite{DGW2023} in the special case of matched densities, constant mobilities and a regular potential $\fhi$, with the aid of an additional higher-order regularizing term such as $\gamma_4 \vert\Delta \phi\vert^2$, $\gamma_4>0$, in the Helmholtz free energy $\ehi$. For the resulting hydrodynamic system subject to suitable boundary and initial conditions, the authors proved the existence of a global weak solution in both two and three dimensions. Moreover, with stronger regularity assumptions on the initial data, the existence of a unique global (resp. local) strong solution in two (resp. three) dimensions was achieved. In the two-dimensional case, the authors further established the uniqueness and instantaneous regularization of global weak solutions. In particular, it was shown that the surfactant concentration stays uniformly away from the pure states $0, 1$ after some positive time. In the current contribution, under the choice of an interaction energy density like \eqref{ab-G}, we are able to prove the existence of global weak solutions for the general case with unmatched densities, variable mobilities and two singular potentials, but without the above mentioned higher-order regularization. In fact, the preserved physical bounds of $\phi$ and $\psi$ allow us to handle the bi-variate function $G(\phi, \psi)$ in a rather general form. Further properties of solutions to problem \eqref{eq:nsch1}--\eqref{eq:nschi} such as strong well-posedness, propagation of regularity, and long-time behavior will be investigated in future studies.

\textit{Plan of this paper.} In Section \ref{mr}, we first introduce the functional settings, notation, and basic assumptions. After that, we define weak solutions and state our main results. Section \ref{wpndm} is devoted to the problem \eqref{eq:nsch1}--\eqref{eq:nschi} with non-degenerate mobilities. We construct an implicit time-discretization scheme and solve it with the help of the Leray--Schauder fixed point principle. Then we prove the existence of a global weak solution by passing to the limit. In Section \ref{wpdm}, we establish the existence of a global weak solution to the problem \eqref{eq:nsch1}--\eqref{eq:nschi} with degenerate mobilities.
In the Appendices, we show how to recover the pressure in both cases with non-degenerate mobilities or degenerate mobilities and give the proofs of some auxiliary results.

\section{Main Results}\label{mr}
\subsection{Preliminaries.}\label{pre}
We first introduce some notation and conventions that will be used in this paper.
Let $\mathcal{X}$ be a (real) Banach space with the norm $\left\|\cdot\right\|_{\mathcal{X}}$. We denote by $\mathcal{X}^*$ its dual space and by $\langle \cdot,\cdot \rangle_{\mathcal{X}^*,\mathcal{X}}$ the associated duality pairing. Given a (real) Hilbert space $\mathcal{H}$, its inner product will be denoted by $(\cdot,\cdot)_{\mathcal{H}}$.
Throughout this paper, we assume that $\Omega \subset \mathbb{R}^d$ ($d\in\{2,3\}$) is a bounded domain with a sufficiently smooth boundary $\partial \Omega$.
For the standard Lebesgue and Sobolev spaces in $\Omega$, we use the notation $L^{p}(\Omega)$, $W^{k,p}(\Omega)$ for any $p \in [1,\infty]$ and $k\in \mathbb{N}$,
equipped with the corresponding norms
$\|\cdot\|_{L^{p}(\Omega)}$, $\|\cdot\|_{W^{k,p}(\Omega)}$, respectively. The space $\wkpzero{k}{p}$ denotes the closure of $C^{\infty}_0 \bracket{\Omega}$ in $\wkp{k}{p}$. The dual spaces are denoted by $ \wkp{-k}{p} \stackrel{\rm{def}}{=} \bigbracket{ \wkpzero{k}{p'} }^*$, $\wkpzero{-k}{p} \stackrel{\rm{def}}{=} \bigbracket{ \wkp{k}{p'} }^*$, where $\frac{1}{p}+\frac{1}{p'}=1$, $p\in (1,\infty)$.
When $p = 2$, these spaces are Hilbert spaces and we use the standard convention $H^{k}(\Omega) \stackrel{\rm{def}}{=} W^{k,2}(\Omega)$.
The $L^2$-Bessel potential spaces are denoted by $\hs{s}$, $s\in\r$, which are defined by the restriction of distributions in $H^s\bracket{\rd{d}}$ to $\Omega$. For $s >0$, $\hszero{s}$ denotes the closure of $C^{\infty}_0\bracket{\Omega}$ in $\hs{s}$.
For the sake of simplicity, the norm and the inner product of $L^2(\Omega)$ will be denoted by $\|\cdot\|$ and $(\cdot,\cdot)$, respectively, while the duality pairing between $H^1(\Omega)$ and $H^1(\Omega)^*$ will be denoted by $\langle \cdot,\cdot \rangle$.
In order to avoid confusion, we also use the notation $(\cdot,\cdot)_{\mathcal{D}}$ to point out the integration domain $\mathcal{D}$.
Bold letters will be used for vector-valued
spaces, for example, $\bm{L}^p(\Omega)= L^p(\Omega;\mathbb{R}^d)$, $p \in [1,\infty]$. When there is no ambiguity, for a vector-valued function $\bm{f}:\Omega\to \rd{d}$, the expression $\bm{f} \in \lp{p}$ should be understood as $\bm{f} \in \bm{L}^p(\Omega) $.

Given a measurable set $I$ of $\mathbb{R}$, we introduce the function space $L^p(I;\mathcal{X})$ with $p\in [1,\infty]$, which consists of Bochner measurable $p$-integrable functions (if $p \in [1,\infty)$) or essentially bounded functions (if $p =\infty$) with values in a given Banach space $\mathcal{X}$. If $I=(a,b)$, we simply write $L^p(a,b;\mathcal{X})$. In addition, $f \in \lptbigloc{p}{[0,\infty);\mathcal{X}}$ if and only if $f \in \lptbig{p}{0,T;\mathcal{X}}$ for every $T>0$. The space $L^p\ul([0,\infty);\mathcal{X})$ denotes the uniformly local variant of $L^p(0,\infty;\mathcal{X})$ consisting of all strongly measurable $f:[0,\infty)\to \mathcal{X}$ such that
\begin{equation*}
\|f\|_{L^p\ul([0,\infty);\mathcal{X})}\stackrel{\rm{def}}{=} \mathop\mathrm{sup}\limits_{t\ge0} \|f\|_{L^p(t,t+1;\mathcal{X})} < \infty.
\end{equation*}
If $T\in(0,\infty)$, we set $L^p\ul([0,T);\mathcal{X})\stackrel{\rm{def}}{=}L^p(0,T;\mathcal{X})$.
For every $k\in \mathbb{N}_+$, $f\in \wkpt{k}{p}{0,T;\mathcal{X}}$ if and only if $ \set{(\frac{\mathrm{d}}{\mathrm{d}t})^jf}_{j=0}^k \subset \lpt{p}{0,T;\mathcal{X}}$, where $\dif{f}{t}$ denotes the vector-valued distributional derivative of $f$ and $(\frac{\mathrm{d}}{\mathrm{d}t})^0 f$  just means $f$.
Furthermore, $\wkptbiguloc{k}{p}{[0,\infty);\mathcal{X}}$ is defined by replacing $\lpt{p}{0,T;\mathcal{X}}$ with $\lptbiguloc{p}{[0,\infty);\mathcal{X}}$, and we set
$$ \hst{k}{0,T;\mathcal{X}}=\wkpt{k}{2}{0,T;\mathcal{X}},
\quad  \hstuloc{k}{[0,\infty);\mathcal{X}} = \wkptbiguloc{k}{2}{[0,\infty);\mathcal{X}}.
$$
Let $I=[0,T]$ if $T \in (0,\infty)$ or $I=[0,\infty)$ if $T=\infty$.
Then $BC\bracket{I;\mathcal{X}}$ denotes the Banach space of all bounded and continuous functions $f:I \to \mathcal{X}$ equipped with the supremum norm and $BUC\bracket{I;\mathcal{X}}$ is the subspace of all bounded and uniformly continuous functions.
Furthermore, for every $k \in \mathbb{N_+}$, $BC^k\bracket{I;\mathcal{X}}$ (resp. $BUC^k\bracket{I;\mathcal{X}}$) denotes all functions $f \in BC\bracket{I;\mathcal{X}}$ (resp. $f \in BUC\bracket{I;\mathcal{X}}$), whose Fr\'echet derivatives of order no larger than $k$ exist and belong to $BC\bracket{I;\mathcal{X}}$ (resp. $BUC\bracket{I;\mathcal{X}}$).
Finally, we denote by $BC_{w}\bracket{I;\mathcal{X}}$ the topological vector space of all bounded and weakly continuous functions $f:I \to \mathcal{X}$.

In the subsequent analysis, the following shorthands will be used
$$
H \stackrel{\rm{def}}{=} L^2\left(\Omega\right),
\quad V \stackrel{\rm{def}}{=} H^1\left(\Omega\right),
\quad W \stackrel{\rm{def}}{=}\big\{u\in H^2(\Omega)\ |\ \partial_\mathbf{n}u=0\ \text{a.e. on}\ \partial\Omega\big\}.
$$
As usual, $H$ is identified with its dual, and we have the following continuous, dense, and compact embeddings:
\begin{equation*}
    W \hookrightarrow V  \hookrightarrow H \hookrightarrow V^*.
\end{equation*}
Besides, we recall the interpolation inequality
\begin{equation*}
    \left\|f\right\|^2 \le \xi \left\|\nabla f\right\|^2 + C\left(\xi\right) \left\|f\right\|^2_{V^*}, \quad \forall f \in V,
\end{equation*}
where $\xi \in (0,1)$ is arbitrary and $C\left(\xi\right)$ is a positive constant only depending on $\xi$ and $\Omega$.
For every $f\in V^*$, $\overline{f}$ denotes its generalized mean value over $\Omega$, given by
$$\overline{f}=|\Omega|^{-1}\langle f,1\rangle.$$
If $f\in L^1(\Omega)$, its mean value is given by $$\overline{f}=|\Omega|^{-1}\int_\Omega f \,\mathrm{d}x.$$
Then we recall the well-known Poincar\'{e}--Wirtinger inequality:
\begin{equation*}
\left\|f-\overline{f}\right\|\leq C_P \|\nabla f\|,\quad \forall\,
f\in V,
\end{equation*}
where the positive constant $C_P$ depends only on $\Omega$.

Let $A_N: V \to V^*$ denote the extension of the Laplace operator $-\Delta$ subject to the homogeneous Neumann boundary condition, namely,
\begin{equation*}
    \left\langle A_N f,g \right\rangle = \int_{\Omega} \nabla f \cdot \nabla g \,\mathrm{d}x, \quad \forall\,f,g \in V.
\end{equation*}
Setting
\begin{align*}
    & H_{(0)}\stackrel{\rm{def}}{=}\{f\in H\ |\ \overline{f} =0\},\quad
    V_{(0)} \stackrel{\rm{def}}{=} V \cap H_{(0)},  \quad V^*_{(0)} \stackrel{\rm{def}}{=} \left\{ L \in V^* \ |\   \left\langle L,1 \right\rangle = 0 \right\},
\end{align*}
it is easy to verify that the restriction of $A_N$ to $V_{(0)}$ is a linear isomorphism between $V_{(0)}$ and $V^*_{(0)}$. Thus, we can define the inverse operator $\mathcal{N} \stackrel{\rm{def}}{=} \left(A_N|_{V_{(0)}}\right)^{-1}:V_{(0)}^* \to V_{(0)}$
and the following identities hold
\begin{align*}
&\left\langle A_Nu,\mathcal{N}L \right\rangle = \left\langle L,u \right\rangle,\quad \forall\, u \in V_{(0)},\ L \in V_{(0)}^*,
\\
&\left\langle L_1,\mathcal{N}L_2 \right\rangle = \inner{\nabla \left(\mathcal{N}L_1\right)}{\nabla \left(\mathcal{N}L_2\right)},
\quad \forall\, L_1,\,L_2\ \in V_{(0)}^*.
\end{align*}
Define
\begin{align*}
    &\left\|L\right\|_* \stackrel{\rm{def}}{=} \left\| \nabla \left(\mathcal{N}L\right)\right\| = \sqrt{\left\langle L,\mathcal{N}L \right\rangle},\quad \forall\, L\in V^*_{(0)}, \\
    &\left\|L\right\|_{-1}^2 \stackrel{\rm{def}}{=} \left\|L-\overline{L}\right\|_*^2 + |\overline{L}|^2, \quad \forall\, L\in V^*.
\end{align*}
We observe that $\left\|\cdot\right\|_*$ and $\left\|\cdot\right\|_{-1}$ are equivalent norms in $V^*_{(0)}$ and $V^*$ with respect to the usual dual norms.

For a bounded domain $\Omega \subset \r^d$, $d\in \{2,3\}$, we define
$$C^{\infty}_{0,\sigma}(\Omega) \stackrel{\rm{def}}{=} \{ \bm{f} \in C^{\infty}_0(\Omega)^d \ | \ \div \bm{f} =0 \}$$
and
\begin{align*}
    &\lpsigma{p} \stackrel{\rm{def}}{=} \overline{C^{\infty}_{0,\sigma}(\Omega)}^{\lp{p}^d},
    \quad
    \hssigma{s} \stackrel{\rm{def}}{=} \overline{C^{\infty}_{0,\sigma}(\Omega)}^{\hs{s}^d},
    \quad
    \wkpsigma{k}{p} \stackrel{\rm{def}}{=} \overline{C^{\infty}_{0,\sigma}(\Omega)}^{\wkp{k}{p}^d},
\end{align*}
for all $p \in (1,\infty)$, $k \in \mathbb{N}_+$ and $s > 0$.
The Helmholtz projection corresponding to $\ls$ is denoted by $\mathbb{P}_{\sigma}$ (see, for example, \cite{S2001}). We note that $\psigma \bm{f} = \bm{f} - \nabla p$, where $p \in V_{(0)}$ is the solution to the (weak) Neumann problem
$$ \inner{\nabla p}{\nabla \phi}_{\Omega}
= \inner{\bm{f}}{\nabla \phi}_{\Omega},
\quad  \forall\, \phi \in C^{\infty} \bigbracket{\overline{\Omega}}.
$$
In addition, if $\p \Omega$ is sufficiently smooth (for example, of $C^2$ class), then we can also define the Helmholtz--Weyl decomposition of $\bm{L}^p(\Omega)$ for all $p \in (1,\infty)$ thanks to the well-posedness of corresponding (homogeneous) Neumann boundary value problem (see \cite[Section III.1]{G2011}).
When the Helmholtz--Weyl decomposition is valid and there is no confusion, $\psigma$ also denotes the unique projection operator from $\bm{L}^p(\Omega)$ to $\lpsigma{p}$, whose null space is
$$ \set{ \bm{\omega} \in \bm{L}^p(\Omega) \ |\  \bm{\omega} = \nabla p, \ \mathrm{for} \ \mathrm{some} \ p \in \wkp{1}{p}  }. $$

In the subsequent analysis, the capital letter $C$ will denote a generic positive constant that depends on the structural data of the problem. Its meaning may change from line to line and even within the same chain of computations. Specific dependence will be pointed out if necessary.
Finally, for simplicity, we set
$$Q_{(s,t)}=\Omega \times (s,t), \quad Q_{t}=Q_{(0,t)},\quad Q=Q_{(0,\infty)},$$
and analogously
$$S_{(s,t)} = \p \Omega \times (s,t),\quad  S_t = S_{(0,t)}, \quad  S=S_{(0,\infty)}.$$

\subsection{Basic assumptions.}\label{sectionassumptions}
Let us introduce the assumptions that will be used throughout this paper (cf. \cite{ADG2013ndm,ADG2013dm,AGG2012,B1999}).
\medskip

\ba \label{assumptions} \rm
We assume that $\Omega \subset \rd{d}\ (d\in \{ 2,3\})$ is a bounded domain with a smooth boundary $\partial\Omega$ (here it is enough to assume that $\p \Omega$ is of class $C^2$) and $\beta>0$, $\sigma_2\geq 0$, $c\in (-1,1)$ are given constants.
In addition, we impose the following conditions.
\medskip
\begin{itemize}
    \item[$(\mathbf{H0})$] The average density is given by
    $$
    \rho(\phi) = \halfof{\rho_1 - \rho_2}\phi + \halfof{\rho_1 + \rho_2},
    $$
    where $\rho_1,\rho_2$ are given positive constants. \medskip

    \item[$(\mathbf{H1})$] The potential function
    $F_{\phi}\in C([-1,1]) \cap C^2((-1,1))$ satisfies
    \begin{equation*}
        \limit{s \to -1^+} F_{\phi}'(s) = -\infty,\quad
        \limit{s \to 1^-} F_{\phi}'(s) = + \infty, \quad
        F_{\phi}''(s) \ge \theta_{\phi} > 0, \quad \forall\, s \in (-1,1),
    \end{equation*}
    where $\theta_\phi$ is a given constant.
    Similarly, $F_{\psi}\in C([0,1]) \cap C^2((0,1))$ satisfies
    \begin{equation*}
        \limit{s \to 0^+} F_{\psi}'(s) = -\infty,\quad
        \limit{s \to 1^-} F_{\psi}'(s) = + \infty, \quad
        F_{\psi}''(s) \ge \theta_{\psi} > 0, \quad \forall\, s \in (0,1).
    \end{equation*}
     We set $F_{\psi}(s) = +\infty $ whenever $|s| > 1$ and $F_{\psi}(s) = +\infty $ whenever $s \in (-\infty,0) \cup (1,+\infty)$. Without loss of generality, we take
     \begin{equation*}
         \fhi(0)= \fhi'(0) = \fsi\left(\half\right) = \fsi'\left(\half\right) = 0.
     \end{equation*}
   Besides, we assume that $G(\cdot,\cdot)\in BUC^2(\rd{2})$.\medskip

    \item[$(\mathbf{H2})$] The viscosity coefficient satisfies $\nu \in BUC\bracket{\r}$ with
    $$ \underline{\nu} \stackrel{\rm{def}}{=} \inf\limits_{s \in \r} \set{\nu(s)} >0.$$
    Its upper bound is denoted by $\nu^* \stackrel{\rm{def}}{=} \sup\limits_{s \in \r} \set{\nu(s)}$.
\end{itemize}

Next, we present assumptions on the mobilities $\mbhi$ and $\mbsi$ as well as the function $\sigma_1$ in the case of non-degenerate mobilities.
\medskip
\begin{itemize}
    \item[$(\mathbf{H3})$] (Non-degenerate mobilities). We assume
    $m_{\phi}(\cdot),m_{\psi}(\cdot) \in BC(\mathbb{R}) \cap Lip(\mathbb{R})$. They are bounded from below by
    $$ \underline{\mbhi} \stackrel{\rm{def}}{=} \inf\limits_{s \in \r} \set{\mbhi(s)} >0 \quad  \mathrm{and} \quad \underline{\mbsi} \stackrel{\rm{def}}{=} \inf\limits_{s \in \r} \set{\mbsi(s)} >0. $$
    Their upper bounds are denoted by $$\mbhi^* \stackrel{\rm{def}}{=} \sup\limits_{s \in \r} \set{\mbhi(s)}\quad \text{and}\quad \mbsi^* \stackrel{\rm{def}}{=} \sup\limits_{s \in \r} \set{\mbsi(s)},$$
    respectively. Moreover, we assume
     $$
     \sigma_1(\cdot) \in BUC(\r),\quad \text{with}\quad \sigma_1(s) \ge 0,\quad \forall\, s\in \r.
     $$
\end{itemize}

In the case of degenerate mobilities, instead of $(\mathbf{H3})$, we make the following hypothesis $(\mathbf{H3*})$ on the mobilities $\mbhi$, $\mbsi$ and the function $\sigma_1$.
For convenience, we introduce two auxiliary functions (see \cite{B1999,EG1996})
$$
\whi \in C^2\bigbracket{ (-1,1);\r^+_0 }\quad \text{and}\quad \wsi \in C^2\bigbracket{ (0,1);\r^+_0 }\quad \text{with}\quad \r^+_0=[0,\infty),
$$
defined by
    \begin{align*}
        & \whi''(s) = \frac{1}{\mbhi(s)}, \quad s\in (-1,1),\quad \text{with}\ \ \whi(0)=\whi'(0) = 0, \\
        & \wsi''(s) = \frac{1}{\mbsi(s)}, \quad s\in (0,1),\quad \text{with}\ \ \wsi\left(\half\right)=\wsi'\left(\half\right) = 0.
    \end{align*}

\begin{itemize}
    \item[$(\mathbf{H3*})$] (Degenerate mobilities). Let
    $$m_{\phi}(\cdot) \in C^1\bigbracket{[-1,1]},\quad m_{\psi}(\cdot) \in C^1\bigbracket{[0,1]},\quad \sigma_1(\cdot) \in C\bigbracket{[-1,1]}$$
    be nonnegative functions.
    Moreover, $\mbhi(s)=0$ if and only if $s=\pm 1$, while $\mbsi(s)=0$ if and only if $s=0,1$.
    We assume that
    $$ F_{\phi}'' \mbhi\in C([-1,1]),\quad F_{\psi}''\mbsi\in C([0,1]),$$
    and there exists some positive constant $\alpha$ such that
    \begin{equation}
        F_{\phi}''(s) \mbhi(s) \le \alpha, \ \ \forall\, s \in (-1,1) \quad \text{and}\quad F_{\psi}''(s) \mbsi(s) \le \alpha, \ \  \forall\, s \in (0,1),
        \label{ass:F-m}
    \end{equation}
    \begin{equation}\label{ass:sig1W}
        \sigma_1(s) \bigabs{\whi'(s)} \le \alpha, \quad \forall\, s \in (-1,1).
    \end{equation}
\end{itemize}
\ea

\begin{remark} \rm
\label{mob}
 For convenience in the subsequent analysis, in $(\mathbf{H1})$ we have absorbed the smooth (possibly non-convex) part of $\fhi$, $\fsi$ into $G$. For example, $\fhi$ is usually assumed to be $\lambda$-convex such that $\fhi''(s)+\lambda\geq 0$ for some $\lambda>0$ and all $s\in (-1,1)$. Then we take $$\widetilde{\fhi}(s)=\fhi(s)+\frac{\lambda+1}{2}s^2$$
 as the potential function and absorb the smooth concave part $-\frac{\lambda+1}{2}s^2$ into $G$. The same argument applies to $\fsi$.
 Since the solution $(\phi,\psi)$ takes values in $[-1,1]\times[0,1]$, we only need the functions $\nu(\cdot)$, $m_{\phi}(\cdot)$ in $[-1,1]$, $m_{\psi}(\cdot)$ in $[0,1]$ and $G(\cdot,\cdot)$ in $[-1,1]\times[0,1]$. We can easily extend these functions to $\mathbb{R}$ or $\rd{2}$ such that the conditions in Assumption \ref{assumptions} are fulfilled. In particular, any bi-variable polynomial is admissible for the interaction energy density $G$.
\end{remark}
\medskip

\begin{remark} \rm
    The viscosity function $\nu$ may also depend on $\psi$ (see, e.g., \cite{DGW2023}). Because the concentration of the surfactant is usually very small, its influence on the viscosity of the mixed flow is negligible. For simplicity, we assume here that $\nu$ depends only on $\phi$ (cf. \cite{ZKYWYS2019}). Nevertheless, if $\nu$ is taken as a smooth and strictly positive function of $\phi$ and $\psi$, the results obtained in this study are still valid with minor modifications to the proof.
\end{remark}
\medskip

\begin{remark} \label{entropy_logarithmic_singularity&mobility_monotonicity} \rm
The assumptions $(\mathbf{H3}*)$ yield at least linear degeneracy of the mobility functions near their corresponding end points, namely, there exist some constants $C_{\phi}>0$ and $C_{\psi}>0$ such that
$$ 0\leq \mbhi(s_1) \le C_{\phi} (1-s_1^2), \quad 0\leq \mbsi(s_2) \le C_{\psi} s_2(1-s_2), $$
for all $s_1 \in [-1,1]$ and $s_2 \in [0,1]$, respectively. Indeed, from the following facts
\begin{align*}
    \mbhi (s) & = \mbhi (-1)+ \intinterval{\mbhi'(\tau)}{\tau}{-1}{s}  \le \norm{\mbhi'}_{C([-1,1])} (1+s),\quad \forall\, s\in[-1,0],\\
    \mbhi (s) & = \mbhi (1)- \intinterval{\mbhi'(\tau)}{\tau}{s}{1}  \le \norm{\mbhi'}_{C([-1,1])} (1-s),\quad \forall\, s\in[0,1],
\end{align*}
we can deduce that
$$\mbhi(s) \le \norm{\mbhi'}_{C([-1,1])} (1-s^2),\quad \forall\, s\in[-1,1].$$
In a similar manner, we get
$$\mbsi(s) \le 2\norm{\mbsi'}_{C([-1,1])} s(1-s),\quad \forall\,s\in[0,1].$$
As a consequence, it holds
\begin{equation}\label{entropy_control_log}
    \abs{\ln(1+s_1) - \ln (1-s_1)} \le 2 C_{\phi} \whi'(s_1), \quad  \abs{\ln(s_2) - \ln (1-s_2)}\le C_{\psi} \wsi'(s_2),
\end{equation}
for all $s_1 \in (-1,1)$ and $s_2 \in (0,1)$. The above controls along with \eqref{ass:sig1W} enable us to derive uniform estimates for solutions to the approximate problem in Section \ref{wpdm}.

On the other hand, we recall that in \cite{B1999} the second order derivative of the singular potential was required to be monotonic near the pure phases. Take the singular potential function $\fhi$ as an example. If we approximate it with quadratic functions near $\pm 1$ like in \cite{B1999}, the requirement on the monotonicity of $\fhi''$ seems to be indispensable in order to guarantee the following pointwise control:
$$
\fhi^{\epsilon}(s) \le \fhi(s),\quad \forall\,s\in (-1,1).
$$
In this study, we are able to get rid of the assumption on the monotonicity of $\fhi''$ by applying a different way to approximate the singular potential $\fhi$. Indeed, the main part of the approximate potential $\tfhi^{\epsilon}$ introduced in Section \ref{wpdm} naturally satisfies
$$ \bigbracket{\tfhi^{\epsilon}}''(s) \le \fhi''(s),\quad \forall\,s\in(-1,1).
$$
A similar property holds for the potential $\tfsi^{\epsilon}$ as well. Differently from \cite{B1999}, here we construct the approximate solutions by using the result for the case of non-degenerate mobilities and singular potentials, then the (almost everywhere) pointwise convergence of phase-field variables ensures that their limits also stay in the corresponding physical intervals.
\end{remark}

\subsection{Statement of results.}\label{state_results}
We first introduce the notion of global weak solution to problem \eqref{eq:nsch1}--\eqref{eq:nschi} with non-degenerate mobilities.
\medskip

\bd \label{wsnschndmd} \rm
Let $T \in (0,\infty]$ and set $I=[0,\infty)$ if $T= \infty$ or $I = [0,T]$ if $T<\infty$. Suppose that the assumptions $(\mathbf{H0})$, $(\mathbf{H1})$, $(\mathbf{H2})$ and $(\mathbf{H3})$ are satisfied.
For any initial data $\sp_0 \in \ls$, $\psi_0, \phi_0 \in V$ with $(\phi_0,\psi_0) \in [-1,1]\times[0,1]$ almost everywhere in $\Omega$ and $\overline{\phi_0} \in (-1,1)$ and $\overline{\psi_0} \in (0,1)$, we call $(\sp,\phi,\psi,\mhi,\msi)$ with the regularity properties
\begin{align*}
        & \sp \in C_{w} \bigbracket{I;\ls}\cap L^2\loc \bigbracket{I;\hssigma{1}},\\
        & \phi,\psi \in C_{w} \bigbracket{I;V} \cap L^2\loc (I;W) \cap L^{\infty} (Q_T), \\
        &\phi(x,t) \in (-1,1) \ \ \mathrm{and} \ \ \psi(x,t) \in (0,1) \ \ \ae \ \text{in} \ \ Q_T,\\
        & \mhi,\msi \in L^2\loc \bigbracket{I;V},\\
        & F_{\phi}'(\phi),F_{\psi}'(\psi) \in L^2\loc (I;H),
\end{align*}
a (finite energy) weak solution to problem \eqref{eq:nsch1}--\eqref{eq:nschi} on $I$, if it satisfies
\begin{align}
        & - \inner{\rho \sp}{\pdif{\bmtheta}{t}}_{Q_T} - \inner{\rho \sp \otimes \sp}{\nabla \bmtheta}_{Q_T}
        - \biginner{\sp \otimes \mathbf{J}}{\nabla \bmtheta}_{Q_T}
        + \inner{\nu(\phi) D \sp}{D \bmtheta}_{Q_T}
        \nonumber\\
        & \quad = \biginner{\mhi \nabla \phi}{\bmtheta}_{Q_T}+\biginner{\msi \nabla \psi}{\bmtheta}_{Q_T},
         \label{testns}
\end{align}
for all $\bmtheta \in C^{\infty}_0(Q_T)^d$ with $\div \bmtheta = 0$, and
\begin{align}\label{testchphi}
   & - \biginner{\phi}{\pdif{\zeta}{t}}_{Q_T} + \biginner{\sp \cdot \nabla \phi}{\zeta}_{Q_T} + \intinterval{\bigbracket{\overline{\phi} -c} \inner{ \sigma_1 \bracket{\phi} }{ \zeta }}{t}{0}{T} = - \biginner{\mbhi(\phi)\nabla \mhi}{\nabla \zeta}_{Q_T} ,\\
\label{testchpsi}
    & - \biginner{\psi}{\pdif{\zeta}{t}}_{Q_T} + \biginner{\sp \cdot \nabla \psi}{\zeta}_{Q_T} = - \biginner{\mbsi(\psi)\nabla \msi}{\nabla \zeta}_{Q_T},
\end{align}
for all $\zeta \in C_0^{\infty}\big((0,T);C^1\big(\overline{\Omega}\big)\big)$,
equations \eqref{eq:nsch4} and \eqref{eq:nsch6} for chemical potentials $\mu_\phi$, $\mu_\psi$ hold almost everywhere in $Q_T$, the initial conditions
$(\sp,\phi,\psi)|_{t=0} =(\sp_0,\phi_0,\psi_0)$ are satisfied.
Moreover, $\overline{\phi} \in (-1,1)$ and $\overline{\psi}=\overline{\psi_0} \in (0,1)$ and the following energy inequality holds
\begin{align}
    & \eto \bigbracket{\sp(t),\phi(t),\psi(t)} + \intomega{\nu (\phi) |D \sp|^2+ m_{\phi}(\phi)\abs{\gradmhi}^2 +m_{\psi}(\psi) \abs{\gradmsi}^2}{(x,\tau)}{Q_{(s,t)}}\nonumber\\
    & \quad +  \intomega{ \sigma_1 (\phi) \bigbracket{ \overline{\phi} - c } \Bigbracket{ \mhi - \frac{\rho_1 - \rho_2}{4} \abs{\sp}^2 } }{\bracket{x,\tau}}{Q_{\bracket{s,t}}}
    \le \eto \bigbracket{\sp(s),\phi(s),\psi(s)}, \label{energyineq}
\end{align}
for almost all $s \in [0,\infty) \cap I$ (including $s=0$) and all $t \in [s,\infty) \cap I$.
\ed
\medskip

Our first result reads as follows.
\medskip
\bt\label{wsnschndme} Suppose that the assumptions $(\mathbf{H0})$, $(\mathbf{H1})$, $(\mathbf{H2})$ and $(\mathbf{H3})$ are satisfied. For any initial data $\sp_0 \in \ls$, $\psi_0, \phi_0 \in V$ with $(\phi_0,\psi_0) \in [-1,1]\times[0,1]$ almost everywhere in $\Omega$, $\overline{\phi_0} \in (-1,1)$ and $\overline{\psi_0} \in (0,1)$,
problem \eqref{eq:nsch1}--\eqref{eq:nschi} admits a global weak solution $\bracket{\sp,\phi,\psi,\mhi,\msi}$ on $[0,\infty)$ in the sense of Definition \ref{wsnschndmd}.
Moreover, if it holds
    \begin{equation} \label{additional_requirements_uniform_boundedness}
        \norm{\sigma_1(\phi)}_{\lp{\infty}} \bigbracket{\overline{\phi}-c} \in \lpt{1}{0,\infty} \ \ \mathrm{and} \ \  \norm{\sigma_1(\phi(t))}_{\lp{\infty}} \bigbracket{\overline{\phi}(t)-c} \to 0 \ \mathrm{as} \ t \to \infty,
    \end{equation}
    then the weak solution $\bracket{\sp,\phi,\psi,\mhi,\msi}$  satisfies
    \begin{align*}
        & \sp \in BC_{w} \bigbracket{[0,\infty);\ls} \cap L^2 \bigbracket{0,\infty;\hssigma{1}},\\
        & \phi,\psi \in BC_{w} \bigbracket{[0,\infty);V} \cap L^2\ul ([0,\infty);W),\\
        & \mhi,\msi \in L^2\ul \bigbracket{[0,\infty);V} \ \mathrm{with} \ \gradmhi,\gradmsi \in L^2 \bigbracket{0,\infty;\bm{L}^2(\Omega)}, \\
        & F_{\phi}'(\phi),F_{\psi}'(\psi) \in L^2\ul ([0,\infty);H).
\end{align*}
\et

\begin{remark}\label{add-sigphi} \rm
From the mass relation
\begin{equation}
\overline{\phi}(t) -c = \bigbracket{\overline{\phi_0} -c } \mathrm{exp}  \Bigbracket{ -\intinterval{\overline{\sigma_1(\phi)}}{\tau}{0}{t} }, \quad  \forall \, t \ge 0,
\label{mass-rel}
\end{equation}
it is straightforward to check that the condition \eqref{additional_requirements_uniform_boundedness} for $\phi$ can be fulfilled, provided that one of the following assumptions holds:
\smallskip
    \begin{itemize}
        \item the function $\sigma_1(s)$ is identically equal to a constant $\sigma \ge 0$, \smallskip
        \item $\overline{\phi_0} = c$, \smallskip
        \item there exists a positive constant $\sigma_*$ such that $ \min\limits_{s \in [-1,1]} \set{\sigma_1(s)} \ge \sigma_*$.
    \end{itemize}
\end{remark}
\medskip

Let us proceed to investigate the case of degenerate mobilities. We define weak solutions for problem \eqref{eq:nsch1}--\eqref{eq:nschi} in a manner analogous to the framework used in \cite{ADG2013dm}.

Introduce the fluxes
\begin{equation*}
    \jhi = \mbhi (\phi) \nabla \mhi,
    \quad \jsi = \mbsi (\psi) \nabla \msi.
\end{equation*}
We shall work with the new variables $(\jhi,\jsi)$ instead of the chemical potentials $(\mu_\phi,\mu_\psi)$, as estimates for $(\mu_\phi,\mu_\psi)$ and their gradients seem hard to find in the case of degenerate mobilities.
\medskip

\bd\label{wsnschdmd} \rm
Let $T \in (0,\infty]$ and set $I=[0,\infty)$ if $T= \infty$ or $I = [0,T]$ if $T<\infty$.
Suppose that the assumptions $(\mathbf{H0})$, $(\mathbf{H1})$, $(\mathbf{H2})$ and $(\mathbf{H3}*)$ are satisfied.
For any initial data $\sp_0 \in \ls$, $\psi_0, \phi_0 \in V$ with $(\phi_0,\psi_0) \in [-1,1]\times[0,1]$ almost everywhere in $\Omega$, $\overline{\phi_0} \in (-1,1)$,  $\overline{\psi_0} \in (0,1)$, and
$F_{\phi}(\phi_0)$, $\whi(\phi_0)$, $F_{\psi}(\psi_0)$, $\wsi(\psi_0) \in \lp{1}$, we call $(\sp,\phi,\psi,\jhi,\jsi)$ with the regularity properties
    \begin{align*}
        & \sp \in C_{w} \bigbracket{I;\ls}\cap L^2\loc \big( I;\hssigma{1} \big),\\
        & \phi,\psi \in C_{w} \bigbracket{I;V} \cap L^2\loc (I;W) \cap L^{\infty}(Q_T), \\
        &\phi(x,t) \in [-1,1] \ \ \mathrm{and} \ \ \psi(x,t) \in [0,1] \ \ \ae \ \mathrm{in} \ \ Q_T,\\
        & \jhi, \jsi \in L^2\loc \bigbracket{I;\bm{L}^2(\Omega)}
    \end{align*}
a weak solution to problem \eqref{eq:nsch1}--\eqref{eq:nschi} on $I$, if the following conditions are satisfied:
\begin{align*}
        & - \inner{\rho \sp}{\pdif{\bmtheta}{t}}_{Q_T}
        - \inner{\rho \sp \otimes \sp}{\nabla \bmtheta}_{Q_T}
        - \inner{\sp \otimes \mathbf{J}}{\nabla \bmtheta}_{Q_T}
        + \inner{\nu(\phi) D \sp}{D \bmtheta}_{Q_T}
        \nonumber\\
        & \quad = -\inner{ \Delta \phi \nabla \phi}{\bmtheta}_{Q_T} + \sigma_2 \Biginner{\mathcal{N} \bigbracket{ \phi - \overline{\phi} } \nabla \phi}{ \bmtheta }
        -\beta \inner{ \Delta \psi \nabla \psi}{ \bmtheta}_{Q_T},
\end{align*}
with $\mathbf{J}= -\frac{\rho_1-\rho_2}{2}\mathbf{J}_\phi$, for all $\bmtheta \in C_0^{\infty}(Q_T)^d$ with $\div \bmtheta = 0$, and
\begin{align*}
    &- \inner{\phi}{\pdif{\zeta}{t}}_{Q_T}
    + \inner{\sp \cdot \nabla \phi}{\zeta}_{Q_T} + \intinterval{\bigbracket{\overline{\phi} -c} \inner{ \sigma_1 \bracket{\phi} }{ \zeta }}{t}{0}{T}
    = - \inner{\jhi}{\nabla \zeta}_{Q_T} ,\\
    &- \inner{\psi}{\pdif{\zeta}{t}}_{Q_T}
    + \inner{\sp \cdot \nabla \psi}{\zeta}_{Q_T}
    = - \inner{\jsi}{\nabla \zeta}_{Q_T},
\end{align*}
for all $\zeta \in C_0^{\infty}\big((0,T);C^1\big(\overline{\Omega}\big)\big)$,
where for all $T_0 \in I$, it holds
\begin{align*}
    & \inner{\jhi}{\bmeta}_{Q_{T_0}} =
    \intomega{\Bigbracket{ F_{\phi}''(\phi)\nabla \phi + \frac{\p^2 G}{\p \phi^2}(\phi,\psi) \nabla \phi + \frac{\p^2 G}{\p \phi \p \psi}(\phi,\psi) \nabla \psi}\cdot (\mbhi(\phi)\bmeta )}{(x,t)}{Q_{T_0}} \nonumber\\
    &  \qquad \qquad \qquad + \inner{\Delta \phi}{\div (\mbhi(\phi) \bmeta)}_{Q_{T_0}} + \sigma_2 \biginner{\nabla \mathcal{N} \bigbracket{\phi - \overline{\phi}}}{ \mbhi(\phi) \bmeta}_{Q_{T_0}}
\end{align*}
and
\begin{align*}
    & \inner{\jsi}{\bmeta}_{Q_{T_0}} =
    \intomega{ \Bigbracket{F_{\psi}''(\psi)\nabla \psi + \frac{\p^2 G}{\p \psi^2}(\phi,\psi) \nabla \psi + \frac{\p^2 G}{\p \phi \p \psi}(\phi,\psi) \nabla \phi}\cdot (\mbsi(\psi)\bmeta )}{(x,t)}{Q_{T_0}} \nonumber\\
    &  \qquad \qquad \qquad + \inner{\beta \Delta \psi}{\div (\mbsi(\psi) \bmeta)}_{Q_{T_0}}
\end{align*}
for all $\bmeta \in L^2\bigbracket{0,T_0;\bm{H}^1(\Omega)} \cap \bm{L}^{\infty}(Q_{T_0})$ with $\bmeta \cdot \mathbf{n} = 0$ on $S_{T_0}$, and the initial conditions $(\sp,\phi,\psi)|_{t=0} =(\sp_0,\phi_0,\psi_0)$ are satisfied.
\ed
\medskip

Then we state the second result of this paper.
\medskip
\bt\label{wsnschdme}
Suppose that the assumptions $(\mathbf{H0})$, $(\mathbf{H1})$, $(\mathbf{H2})$ and $(\mathbf{H3}*)$ are satisfied. For any initial data $\sp_0 \in \ls$, $\psi_0, \phi_0 \in V$ with $(\phi_0,\psi_0) \in [-1,1]\times[0,1]$ almost everywhere in $\Omega$, $\overline{\phi_0} \in (-1,1)$, $\overline{\psi_0} \in (0,1)$ and
$F_{\phi}(\phi_0)$, $\whi(\phi_0)$, $F_{\psi}(\psi_0)$, $\wsi(\psi_0) \in \lp{1}$, problem \eqref{eq:nsch1}--\eqref{eq:nschi} admits a global weak solution $(\sp,\phi,\psi,\jhi,\jsi)$ on $[0,\infty)$ in the sense of Definition \ref{wsnschdmd}. Moreover, there exist $\mathbf{\widehat{J}}_{\phi},\mathbf{\widehat{J}}_{\psi} \in L^2\loc(Q)^d$ such that
\begin{equation*}
    \jhi = \sqrt{\mbhi(\phi)}\, \mathbf{\widehat{J}}_{\phi}, \quad  \jsi = \sqrt{\mbsi(\psi)}\, \mathbf{\widehat{J}}_{\psi}.
\end{equation*}
Assume in addition, $\sigma_1= 0$ or $\overline{\phi_0}=c$, then we can construct a global weak solution that satisfies the following energy inequality
\begin{align}\label{energyinequalitydm}
    &\eto \bigbracket{\sp(t),\phi(t),\psi(t)} + \intomega{\left(\nu (\phi) |D \sp|^2+ \bigabs{\mathbf{\widehat{J}}_{\phi}}^2 + \bigabs{\mathbf{\widehat{J}}_{\psi}}^2 \right)}{(x,\tau)}{Q_{(s,t)}} \nonumber \\
    & \quad \le \eto \bigbracket{\sp(s),\phi(s),\psi(s)},
\end{align}
for almost all $s \in [0,\infty) \cap I$ (including $s=0$) and all $t \in [s,\infty) \cap I$, which further implies
    \begin{align*}
        & \sp \in BC_{w} \bigbracket{[0,\infty);\ls}\cap L^2 \big( 0,\infty;\hssigma{1} \big),\\
        & \phi,\psi \in BC_{w} \bigbracket{[0,\infty);V}, \\
        & \jhi, \jsi \in L^2 \bigbracket{0,\infty;\bm{L}^2(\Omega)}.
    \end{align*}
\et

\begin{remark}\label{rem:mo-pot}\rm
In Theorem \ref{wsnschdme} we focus on degenerate mobilities and singular potentials in a general form \cite{EG1996,B1999}. In particular, condition \eqref{ass:F-m} indicates that the singularity of $F_\phi$, $F_\psi$ must be balanced by the degeneracy of $m_\phi$, $m_\psi$, respectively. On the other hand, regularity potentials could also be considered under certain specific choices of degenerate mobilities. For instance, with minor modifications of the argument in \cite{ADG2013dm}, it is straightforward to check that the conclusion of Theorem \ref{wsnschdme} still holds for
$F_\phi, F_\psi\in C^1(\mathbb{R})$
and
$$
m_\phi(s)
=\begin{cases}
1-s^2, &\quad  s\in [-1,1],\\
0,&\quad \text{else},
\end{cases}
\qquad
m_\psi(s)
=\begin{cases}
s(1-s), &\quad  s\in [0,1],\\
0,&\quad \text{else}.
\end{cases}
$$
Other mobilities with linear degeneracy near the degenerate points are also possible.
\end{remark}
\medskip

\begin{remark}\rm
Since estimates for the chemical potential $\mu_\phi$ are unavailable when the mobility $m_\phi$ is degenerate, we cannot obtain an energy inequality in the general case with mass transfer (cf. \eqref{energyineq}). If $\sigma_1 = 0$ or $\overline{\phi_0}=c$, we have
$$\sigma_1 \bracket{\phi} \bracket{\overline{\phi}-c} \equiv 0,$$
then the energy inequality \eqref{energyinequalitydm} holds. With the aid of \eqref{energyinequalitydm}, we can derive uniform-in-time estimates for global weak solutions on $[0,\infty)$. Otherwise, the additional term due to mass transfer
$$
\intomega{ \sigma_1 \bracket{\phi^\epsilon}  \bigbracket{\overline{\phi^\epsilon} - c} \mhi^\epsilon}{\bracket{x,\tau}}{Q_{(s,t)}}
$$
appearing in the energy inequality for approximate solutions yields terms like
$$\norm{\Delta \phi^\epsilon}_{\lpt{2}{s,t;\ltwo}}^2 + C_1(t-s).$$
Then we need to apply entropy estimates to control $\norm{\Delta \phi^\epsilon}_{\lpt{2}{s,t;\ltwo}}^2$, which leads to another term of the form $C_2(t-s)$ (see the proof of Lemma \ref{energyestimateapdm} for further details). This only gives time-dependent estimates of the solution.
\end{remark}
\medskip

\begin{remark}\rm
Different assumptions have been imposed in the literature to handle the mass source \cite{CFG2024,M2019,GLRS2022,JWZ2015}.
In particular, if we replace the assumption \eqref{ass:sig1W} on $\sigma_1$ by those as in \cite{CFG2024}, the conclusion of Theorem \ref{wsnschdme} still holds.
More precisely, instead of \eqref{ass:sig1W}, let us assume that $\mbhi$ is monotonic close to the pure phases and
$$
\sigma_1(s) \le \alpha \mbhi(s), \quad \forall\, s \in \r,
$$
and denote the corresponding set of assumptions by $(\mathbf{H3**})$.
We observe that for some typical choices of $\mbhi$, such as
$$\mbhi(s) = (1-s^2)^k,\quad k \ge 1,$$
it holds
$$
\mbhi(s)\whi'(s) \to 0 \quad \mathrm{as} \quad s \to \pm1,
$$
which implies that if $(\mathbf{H3**})$ is satisfied then $(\mathbf{H3*})$ holds as well.
Actually, this property holds in more general settings. To this end, we recall that $\mbhi$ is assumed to be monotonically decreasing (resp. increasing) in $[1-2^{-n_0},1]$ (resp. $[-1,-1+2^{-n_0}]$) for some sufficiently large $n_0 \in \mathbb{N}_+$. This yields
\begin{align*}
            0\leq \whi'(s) \mbhi(s) & = \intinterval{\frac{\mbhi(s)}{\mbhi(\tau)}}{\tau}{0}{1-2^{-n_0}}
            + \intinterval{\frac{\mbhi(s)}{\mbhi(\tau)}}{\tau}{1-2^{-n_0}}{1-2^{-m}}
            + \intinterval{\frac{\mbhi(s)}{\mbhi(\tau)}}{\tau}{1-2^{-m}}{s} \\
            & \le C\mbhi(s) + \frac{\mbhi(s)}{\mbhi(1-2^{-m})} + 2^{-m},
\end{align*}
for all $m \ge n_0+1$ and $ s \in [1 - 2^{-m},1)$, where $C>0$ is independent of $m$ and $s$. Since
$$ \limsup\limits_{m \to \infty}\, \limsup\limits_{s \to 1^-} \Bigbracket{C\mbhi(s)+\frac{\mbhi(s)}{\mbhi(1-2^{-m})} + 2^{-m}} = 0, $$
we find that
$$ \mbhi(s)\whi'(s) \to 0\quad \text{as}\ s \to 1^-.$$
Analogously, it holds $ \mbhi(s)\whi'(s) \to 0$ as $s \to -1^+$.
\end{remark}
\medskip

\begin{remark} \rm
    We note that in \cite{EG1996,B1999}, the potential function $\fhi$ was not required to be continuous at the pure phases. Here, this assumption is imposed for the derivation of the energy inequality \eqref{energyinequalitydm} (the validity of an energy inequality was not addressed in \cite{EG1996,B1999}). Indeed, it requires that $\eto^{\epsilon}\bigbracket{\sp^{\epsilon}(t),\phi^{\epsilon}(t),\psi^{\epsilon}(t)}$ converges to $\eto\bigbracket{\sp(t),\phi(t),\psi(t)}$ for almost all $t \in \r_+$ and $t=0$ as $\epsilon \to 0$ (cf. \cite[Section 5.3]{ADG2013ndm}). Thanks to the boundedness of $\fhi$ and $\fsi$, we can apply Lebesgue's dominated convergence theorem to show that as $\epsilon \to 0$, it holds
    $$
    \intomega{\tfhi^{\epsilon} \bracket{\phi^{\epsilon}(t)}}{x}{\Omega} \to \intomega{\fhi \bracket{\phi(t)}}{x}{\Omega},
    \quad \intomega{\tfsi^{\epsilon} \bracket{\psi^{\epsilon}(t)}}{x}{\Omega} \to \intomega{\fsi \bracket{\psi(t)}}{x}{\Omega},
    $$
    for almost all $t \in \r_+$ and $t=0$.
    On the other hand, if we only assume $\fhi \in C^2\bigbracket{(-1,1)}$ and $\fsi \in C^2\bigbracket{(0,1)}$, then in the case of $\sigma_1 =0$ or $\overline{\phi_0}=c$ we can obtain a weaker version of the energy inequality for the global weak solution $\bigbracket{\sp,\phi,\psi,\jhi,\jsi}$, that is,
    $$ \eto \bigbracket{\sp(t),\phi(t),\psi(t)}
    + \intomega{\left(\nu (\phi) |D \sp|^2+ \bigabs{\mathbf{\widehat{J}}_{\phi}}^2 + \bigabs{\mathbf{\widehat{J}}_{\psi}}^2 \right)}{(x,\tau)}{Q_t}
    \le \eto \bigbracket{\sp_0,\phi_0,\psi_0}, $$
    for all $t \in [0,\infty)$. In this case, we only need
    $$ \liminf\limits_{\epsilon \to 0} \eto^{\epsilon}\bigbracket{\sp^{\epsilon}(t),\phi^{\epsilon}(t),\psi^{\epsilon}(t)} \ge \eto\bigbracket{\sp(t),\phi(t),\psi(t)}, $$
    which can be guaranteed by Fatou's Lemma.
\end{remark}
\medskip

\begin{remark} \rm
    In the case of degenerate mobilities, the property that $\phi$ and $\psi$ take their values in the corresponding physical intervals $[-1,1]$ and $[0,1]$ is guaranteed by the degeneracy of $\mbhi$ and $\mbsi$ near the pure states, instead of the singularity of $\fhi'\bracket{\phi}$ and $\fsi'\bracket{\psi}$. When the order of degeneracy is sufficiently high for the entropy function $\whi$ (resp. $\wsi$) to possess singularity at the pure phases, it follows from the argument in \cite[Section 3.4]{B1999} and the integrability of $\whi\bracket{\phi}$ (resp. $\wsi\bracket{\psi}$) that $\phi \in (-1,1)$ (resp. $\psi \in (0,1)$) almost everywhere in $Q_T$.
    For example, if $\mbhi(s) = (1-s^2)^2$, $s \in (-1,1)$, then we have
    $$\whi(s) = \left(s-\frac{3}{2}\right) \ln{(1-s)} - \left(s+\frac{3}{2}\right) \ln{(1+s)},$$
    which tends to $+\infty$ as $s\to1^-$ or $s \to -1^+$.
\end{remark}
\medskip

\begin{remark}\rm
As in \cite{ADG2013ndm} (see also Section \ref{proofexistencendm}), for both non-degenerate mobilities and degenerate mobilities, we can show that
$$\pdif{ \psigma \bracket{\rho \sp} }{t} \in \lptbigloc{\frac{8}{7}}{[0,\infty);\wkpsigma{1}{4}^*}$$
by comparison in the weak formulation of \eqref{eq:nsch1}. Then applying an argument similar to that in \cite{BF2013}, we can prove the existence of a unique pressure function $\pi \in \wkptbigloc{-1}{\infty}{[0,\infty);\lp{\frac{4}{3}}}$ with $\overline{\pi} = 0$ such that
        $$ \pdif{ \bracket{\rho \sp} }{t} + \div \bigbracket{ \sp \otimes \bigbracket{\rho \sp +\mathbf{J}} } - \div \bracket{\nu\bracket{\phi} D \sp} + \nabla \pi = \bm{f} $$
holds in the sense of distributions, where
\begin{equation*}
            \bm{f} =\left\{
            \begin{aligned}
                & \mhi \nabla \phi + \msi \nabla \psi,
                &&\quad  \text{if} \ \mbhi \ \text{and} \ \mbsi \ \text{are non-degenerate}, \\
                & -\Delta \phi \nabla \phi - \beta \Delta \psi \nabla \psi + \sigma_2 \mathcal{N} \bracket{\phi-\overline{\phi}}\nabla \phi, && \quad \text{if} \ \mbhi \ \text{and} \ \mbsi \ \text{are degenerate}.
            \end{aligned}
            \right.
\end{equation*}
Further details can be found in the Appendix.
\end{remark}
\medskip

\section{Proof of Theorem \ref{wsnschndme}: the Case of Non-degenerate Mobilities}\label{wpndm}\medskip

In this section, we prove Theorem \ref{wsnschndme} on the existence of global weak solutions to problem \eqref{eq:nsch1}--\eqref{eq:nschi} in the case of non-degenerate mobilities. Thus, we assume that assumptions $(\mathbf{H0})$, $(\mathbf{H1})$, $(\mathbf{H2})$ and $(\mathbf{H3})$ are satisfied.

The proof of Theorem \ref{wsnschndme} is divided into two steps. First, we construct a family of approximate solutions through a suitable implicit time-discretization scheme by adapting the idea in \cite{ADG2013ndm}. Then we derive uniform estimates and pass to the limit to establish the existence of a weak solution.
\smallskip

\subsection{Implicit time discretization.}\label{implicit_time_discret}
In this subsection, we work with a slightly stronger assumption such that $\mbhi, \mbsi \in BUC^1(\r)$.
As in \cite[Theorem 4.3]{AW2007}, for the energy functional $\widetilde{E}_\phi:H\to \mathbb{R}$ with the effective domain $$\mathcal{D}(\widetilde{E}_\phi)=\big\{\phi\in V\ |\  \phi\in [-1,1]\ \text{a.e. in}\ \Omega\big\}$$
such that
$$
\widetilde{E}_\phi(\phi)=
\begin{cases}
\displaystyle{\int_\Omega\left(\dfrac12|\nabla \phi|^2+  \fhi(\phi)\right)\mathrm{d}x},
&\quad \text{if}\ \phi\in \mathcal{D}(\widetilde{E}_\phi),\\
+\infty,&\quad \text{else},
\end{cases}
$$
we define the domain of its
subgradient as
$$
\mathcal{D}(\partial \widetilde{E}_\phi)=\big\{ \phi\in W\ |\
\fhi'(\phi)\in H,\ \fhi''(\phi)|\nabla \phi|^2\in L^1(\Omega)\big\}.
$$
Analogously, for
$$
\widetilde{E}_\psi(\psi)=
\begin{cases}
\displaystyle{\int_\Omega\left(\dfrac{\beta}{2}|\nabla \psi|^2+  \fsi(\psi)\right)\mathrm{d}x},
&\quad \text{if}\ \psi\in \mathcal{D}(\widetilde{E}_\psi),\\
+\infty,&\quad \text{else},
\end{cases}
$$
with
$$\mathcal{D}(\widetilde{E}_\psi)=\big\{\psi\in V\ |\  \phi\in [0,1]\ \text{a.e. in}\ \Omega\big\},$$
we define
$$
\mathcal{D}(\partial \widetilde{E}_\psi)=\big\{ \psi\in W\ |\
\fsi'(\psi)\in H,\ \fsi''(\psi)|\nabla \psi|^2\in L^1(\Omega)\big\}.
$$
Applying the same argument for \cite[Lemma 4.1]{AW2007}, we can verify that both
$\widetilde{E}_\phi$ and $\widetilde{E}_\psi$
are proper, lower semi-continuous, convex functionals.

Similar results hold for the case with a prescribed mean value. For any given $a\in \mathbb{R}$, set
$$
H_{(a)}=\{f\in H\ |\ \overline{f}=a\},\quad V_{(a)}=V\cap H_{(a)}, \quad W_{(a)}= W\cap H_{(a)}.
$$
Then for $a \in (-1,1)$, we define the functional $\widetilde{E}_{\phi,a} : H_{(a)} \to \mathbb{R}$ as $ \widetilde{E}_\phi|_{H_{(a)}}$ with the effective domain $$\mathcal{D}(\widetilde{E}_{\phi,a})=\big\{\phi\in V_{(a)}\ |\  \phi\in [-1,1]\ \text{a.e. in}\ \Omega\big\}$$
such that
$$
\widetilde{E}_{\phi,a} (\phi)=
\begin{cases}
\displaystyle{\int_\Omega\left(\dfrac12|\nabla \phi|^2+  \fhi(\phi)\right)\mathrm{d}x},
&\quad \text{if}\ \phi\in \mathcal{D}(\widetilde{E}_{\phi,a}),\\
+\infty,&\quad \text{else}.
\end{cases}
$$
According to \cite[Lemma 4.1]{AW2007}, we find that $\widetilde{E}_{\phi,a}$ is proper, lower semi-continuous and convex. Using \cite[Theorem 4.3]{AW2007}, we also see that the domain of its subgradient $\partial \widetilde{E}_{\phi,a}$ is
$$
\mathcal{D}(\partial \widetilde{E}_{\phi,a})=\big\{ \phi\in W_{(a)} \ |\
\fhi'(\phi)\in H,\ \fhi''(\phi)|\nabla \phi|^2\in L^1(\Omega)\big\},
$$
where
$$\partial \widetilde{E}_{\phi,a} (\phi) = -\Delta \phi + P_0 F'_{\phi} (\phi).$$
Here, $P_0$ is the orthogonal projection from $H$ to $H_{(0)}$ given by
$$ P_0 f = f - \overline{f}, \quad  \forall f \in H. $$
For any $a\in (-1,1)$, $\partial \widetilde{E}_{\phi,a}$ is a maximal monotone operator. Moreover, it holds
\begin{align}
\norm{\phi}^2_W + \norm{F'_\phi (\phi)}^2 + \intomega{ F''_\phi (\phi) | \nabla \phi |^2 }{x}{\Omega} \le C \bigbracket{ \norm{\partial \widetilde{E}_{\phi,a} (\phi)}^2 + \norm{\phi}^2 },
\notag
\end{align}
where the constant $C>0$ is independent of $\phi \in \mathcal{D} (\partial \widetilde{E}_{\phi,a})$.  An analogous consideration applies to $\widetilde{E}_\psi$. For any given $b \in (0,1)$, we define $\widetilde{E}_{\psi,b} \stackrel{\rm{def}}{=} \widetilde{E}_\psi |_{H_{(b)}}: H_{(b)} \to \mathbb{R}$ as
$$
\widetilde{E}_{\psi,b} (\psi)=
\begin{cases}
\displaystyle{\int_\Omega\left(\dfrac{\beta}{2}|\nabla \psi|^2+  \fsi(\psi)\right)\mathrm{d}x},
&\quad \text{if}\ \psi\in \mathcal{D}(\widetilde{E}_{\psi,b}),\\
+\infty,&\quad \text{else},
\end{cases}
$$
with
$$\mathcal{D}(\widetilde{E}_{\psi,b})=\big\{\psi\in V_{(b)}\ |\  \psi\in [0,1]\ \text{a.e. in}\ \Omega\big\}.$$
Then we have
$$
\mathcal{D}(\partial \widetilde{E}_{\psi,b})=\big\{ \psi\in W_{(b)} \ |\
\fsi'(\psi)\in H,\ \fsi''(\psi)|\nabla \psi|^2\in L^1(\Omega)\big\},
$$
with
$$\partial \widetilde{E}_{\psi,b} (\psi) = -\beta \Delta \psi + P_0 F'_{\psi} (\psi),$$
which is a maximal monotone operator. Besides, the following estimate holds
$$ \norm{\psi}^2_W + \norm{F'_\psi (\psi)}^2 + \intomega{ F''_\psi (\psi) | \nabla \psi |^2 }{x}{\Omega} \le C \bigbracket{ \norm{\partial \widetilde{E}_{\psi,b} (\psi)}^2 + \norm{\psi}^2 }, $$
where the constant $C>0$ is independent of $\psi \in \mathcal{D} (\partial \widetilde{E}_{\psi,b})$.

Now we present the implicit time discretization.
Let $h=\frac{1}{N}$ for $N \in \mathbb{N}_+$.
Given
$$
\sp^k \in \ls,\ \ \  \phi^k,\psi^k \in V   \ \text{with}\ F_{\phi}' \bigbracket{ \phi^k },\,F_{\psi}' \bigbracket{ \psi^k } \in H\ \ \ \text{and}\ \rho^k = \rho \bigbracket{ \phi^k },
$$
we look for functions
$$
\big(\sp^{k+1},\phi^{k+1},\psi^{k+1},\widehat{\mhi}^{k+1},\widehat{\msi}^{k+1} \big) = (\sp,\phi,\psi,\widehat{\mhi},\widehat{\msi}),
$$
satisfying
    \begin{equation*}
        \sp \in \hssigma{1},\quad
        \phi \in \mathcal{D} \bigbracket{ \mathrm{\p}\widetilde{E}_{\phi,a}},\quad
        \psi \in \mathcal{D} \bigbracket{ \mathrm{\p}\widetilde{E}_{\psi,b}} ,\quad
         \widehat{\mhi}\in W_{(0)},\quad  \widehat{\msi} \in W_{(0)},
    \end{equation*}
with prescribed mean values
$$
a= - h\overline{\sigma_1(\phi^k)}(\overline{\phi^k}-c)+ \overline{\phi^k},\qquad b=\overline{\psi^k},
$$
as a solution to the following nonlinear discrete system:
\begin{align}
    & \Biginner{\frac{\rho \sp - \rho^k \sp^k}{h}}{\bmtheta} + \biginner{\div \bigbracket{\rho^k \sp \otimes \sp}}{\bmtheta} + \biginner{\nu \big(\phi^k \big) D \sp}{\nabla \bmtheta} +\biginner{\div \bracket{\sp \otimes \mathbf{J}}}{\bmtheta}
    \nonumber \\
    &\quad = \biginner{ \widehat{\mhi} \nabla \phi^k }{\bmtheta}+\biginner{\widehat{\msi} \nabla \psi^k}{\bmtheta},
    \qquad \forall\, \bmtheta \in C^{\infty}_{0,\sigma}(\Omega),
    \label{nsdiscrete1}\\
    &\mathbf{J} = \mathbf{J}^{k+1} = -\gamma \mbhi \big(\phi^k \big) \nabla \widehat{\mhi},  \qquad \mathrm{a.e.\ in} \ \Omega,
    \nonumber\\
    & \frac{\phi -\phi^k}{h} + \sp \cdot \nabla \phi^k + \sigma_1\bigbracket{\phi^k}
    \Bigbracket{\overline{\phi^k} -c} = \div \bigbracket{\mbhi \bigbracket{\phi^k} \nabla \widehat{\mhi}}, \qquad \mathrm{a.e.\ in} \ \Omega,
    \label{chdiscrete1}\\
    & \widehat{\mhi} - P_0 G_{\phi} \bigbracket{\phi,\phi^k,\psi} - \sigma_2 \mathcal{N} \bigbracket{\phi - \overline{\phi}} = -\Delta \phi + P_0 F_{\phi}'(\phi),
    \qquad \quad \mathrm{a.e.\ in} \ \Omega,
    \label{chdiscrete2}\\
    & \frac{\psi -\psi^k}{h} + \sp \cdot \nabla \psi^k = \div \bigbracket{\mbsi \bigbracket{\psi^k} \nabla \widehat{\msi}}, \qquad \quad \ \ \mathrm{a.e.\ in} \ \Omega,
    \label{chdiscrete3}\\
    & \widehat{\msi} - P_0 G_{\psi} \bigbracket{\phi^k,\psi,\psi^k} = -\beta \Delta \psi + P_0 F_{\psi}'(\psi),
    \qquad \mathrm{a.e.\ in} \ \Omega.
    \label{chdiscrete4}
\end{align}
Here, we denote
$$\gamma \stackrel{\rm{def}}{=} \frac{\rho_1 - \rho_2}{2},$$
and the functions $G_\phi$, $G_\psi$ are given by
\begin{equation}\notag
    G_{\phi}(a,b,c)\stackrel{\rm{def}}{=}\left\{
    \begin{aligned}
        &\frac{G(a,c)-G(b,c)}{a-b} , \ &  \mathrm{if} \ a \neq b,\\
        &\pdif{G}{\phi} (a,c), \ &  \mathrm{if} \ a=b,
    \end{aligned}
    \right.
\end{equation}
\begin{equation}\notag
    G_{\psi}(c,a,b)\stackrel{\rm{def}}{=}\left\{
    \begin{aligned}
        &\frac{G(c,a)-G(c,b)}{a-b} , \ &  \mathrm{if} \ a \neq b,\\
        &\pdif{G}{\psi} (c,a), \ &  \mathrm{if} \ a=b,
    \end{aligned}
    \right.
\end{equation}
for $(a,b,c) \in \rd{3}$. Then it is straightforward to check that $G_{\phi}(\cdot,\cdot,\cdot)$, $G_{\psi}(\cdot,\cdot,\cdot) $ belong to $ BUC^1(\rd{3})$ and
$$\norm{G_{\phi}}_{{BUC}^1(\rd{3})}\le \norm{G}_{BUC^2(\rd{2})},
\quad\norm{G_{\psi}}_{BUC^1(\rd{3})} \le \norm{G}_{BUC^2(\rd{2})}.
$$
\begin{remark} \rm
    Since
    $$\rho(\phi) = \gamma \phi + \frac{ \rho_1 + \rho_2 }{2}\quad \text{and}\quad \mathbf{ J } = -\gamma \mbhi \bigbracket{\phi^k} \nabla \widehat{\mhi},$$
    it follows from \eqref{chdiscrete1} that
    $$ \frac{\rho - \rho^k }{h} + \sp \cdot \nabla \rho^k + \gamma \sigma_1 \bigbracket{\phi^k} \Bigbracket{ \overline{\phi^k} -c } = - \div \mathbf{ J } .$$
    Substituting the above equation into \eqref{nsdiscrete1}, we obtain the following equivalent form:
\begin{align}
    & \Biginner{\frac{\rho \sp - \rho^k \sp^k}{h}}{\bmtheta} + \Biginner{ \Bigbracket{\div \mathbf{J} - \frac{\rho-\rho^k}{h}-\sp \cdot \nabla \rho^k - \gamma  \sigma_1\bigbracket{\phi^k} \Bigbracket{\overline{\phi^k} -c}}\halfof{\sp}}{\bmtheta} \nonumber\\
    & \qquad + \biginner{\div \bigbracket{\rho^k \sp \otimes \sp}}{\bmtheta}  + \biginner{\nu \bigbracket{\phi^k} D \sp}{D \bmtheta} +\inner{\bracket{\mathbf{J}\cdot \nabla}\sp}{\bmtheta} \nonumber \\
    & \quad = \biginner{ \widehat{\mhi} \nabla \phi^k }{\bmtheta}+\biginner{ \widehat{\msi} \nabla \psi^k }{\bmtheta},\label{nsdiscrete2}
\end{align}
for all $\bmtheta \in C^{\infty}_{0,\sigma}(\Omega)$.
\end{remark}
\medskip

\begin{remark} \label{averagephi} \rm
Assume $\overline{\phi^k}\in (-1,1)$, $\overline{\psi^k}\in (0,1)$.
Taking the spatial average on both sides of equations \eqref{chdiscrete1} and \eqref{chdiscrete3}, we have
\begin{equation*}
    \overline{\phi^{k+1}} -c
    = \prod\limits_{n=0}^k
    \Bigbracket{1-h \overline{\sigma_1\bigbracket{\phi^n}}}
    \bigbracket{\overline{\phi_0} -c}
    \quad
    \text{and}
    \quad
    \overline{\psi^{k+1}} = \overline{\psi^k},
    \qquad \forall\, k\in\mathbb{N}.
\end{equation*}
As a consequence, if $h\sigma_1^* \le \frac12$, where $$\sigma_1^* \stackrel{\rm{def}}{=} \mathop{\mathrm{max}}\limits_{s \in [-1,1]} \set{\sigma_1(s)},$$
then the sequence $\bigabs{\overline{\phi^k}-c}$ is monotonically non-increasing with respect to $k$.
In particular, $b=\overline{\psi^k}\in (0,1)$, and under the assumption $h\sigma_1^* \le \frac12$, we find that
$$
a= - h\overline{\sigma_1(\phi^k)}(\overline{\phi^k}-c)+ \overline{\phi^k} \in (-1,1).
$$
\end{remark}

In what follows, we prove the existence of a solution to the discrete problem \eqref{nsdiscrete1}--\eqref{chdiscrete4}.
To this end, the following lemma will be useful.
\medskip
\bl \label{apriorichemical}
Assume that $\phi \in \mathcal{D}(\partial \widetilde{E}_\phi)$, $ \psi \in \mathcal{D}(\partial \widetilde{E}_\psi)$ and $\widehat{\mhi}$, $\widehat{\msi} \in V_{(0)}$ solve \eqref{chdiscrete2},  \eqref{chdiscrete4} with given functions $\phi^k,\psi^k \in H^2(\Omega)$ satisfying $ \bigbracket{\phi^k,\psi^k} \in [-1,1]\times[0,1]$ in $\Omega$.
Then, for any $\delta \in (0,1)$, there exist two positive constants $K_{\phi} = K_{\phi}(\delta)$, $K_{\psi} = K_{\psi}(\delta)$ such that if
$$(\overline{\phi},\overline{\psi}) \in [-1+\delta,1-\delta] \times [\delta,1-\delta],$$
then it holds
\begin{align*}
    & \norm{F_{\phi}'(\phi)} \le K_{\phi}\bracket{\norm{\widehat{\mhi}} + 1},\quad    &&\norm{F_{\psi}'(\psi)} \le K_{\psi}\bracket{\norm{\widehat{\msi}} + 1},\\
    &\norm{\p \widetilde{E}_\phi(\phi)} \le K_{\phi} \bracket{\norm{\widehat{\mhi}}+1},\quad
    &&\norm{\p \widetilde{E}_\psi(\psi)} \le K_{\psi} \bracket{\norm{\widehat{\msi}}+1}.
\end{align*}
\el
\bpf
We keep in mind that due to $\phi \in \mathcal{D}(\partial \widetilde{E}_\phi)$, $ \psi \in \mathcal{D}(\partial \widetilde{E}_\psi)$, it holds $(\phi,\psi) \in [-1,1]\times[0,1]$ almost everywhere in $\Omega$. Therefore, the coupling terms $G_{\phi}\bigbracket{\phi,\phi^k,\psi}$ and $G_{\psi}\bigbracket{\phi^k,\psi,\psi^k}$ are bounded by $\norm{G}_{C^1([-1,1]\times[0,1])}$. Testing Equation \eqref{chdiscrete2} with $\bigbracket{\phi - \overline{\phi}}$, we obtain
\begin{align}
    & \biginner{\widehat{\mhi}}{  \phi - \overline{\phi} }_{\Omega} - \biginner{G_{\phi}\bigbracket{\phi,\phi^k,\psi}}{  \phi - \overline{\phi} }_{\Omega} - \sigma_2 \bignorm{ \phi - \overline{\phi} }^2_* \nonumber \\
    & \quad = \norm{ \nabla \phi }^2 + \biginner{\fhi' \bracket{\phi}}{\phi - \overline{\phi}}_{\Omega} \ . \label{discretetestphi-avephi}
\end{align}
According to \cite{MZ2004}, there exists some positive constant $C_{\delta}$ depending on $\delta$ such that
$$ \bignorm{ \fhi' \bracket{\phi} }_{\lp{1}} \le C_{\delta} \biginner{\fhi' \bracket{\phi}}{\phi - \overline{\phi}}_{\Omega} + C_{\delta}.
$$
Substituting this estimate into \eqref{discretetestphi-avephi}, we find
\begin{align*}
    \bignorm{ \fhi' \bracket{\phi} }_{\lp{1}} & \le C_{\delta}  \biginner{\widehat{\mhi}}{ \phi - \overline{\phi} } + 2 C_{\delta} \norm{G_{\phi}}_{C([-1,1]^2\times [0,1])} + C_{\delta} \\
    & \le C_\delta (\norm{\widehat{\mhi}} +1).
\end{align*}
Here, we have used the facts $\phi,\phi^k \in [-1,1]$ and $\psi, \psi^k\in [0,1]$ almost everywhere in $\Omega$. On the other hand, from \eqref{chdiscrete2}, we infer that
\begin{align}
\|\phi\|_{H^2(\Omega)}^2 + \|\fhi'(\phi)\|^2
& \leq C(\| -\Delta \phi + P_0 F_{\phi}'(\phi)\|^2 + \|\phi\|^2)
\notag\\
&\leq C(\|\widehat{\mhi} - P_0 G_{\phi} \bigbracket{\phi,\phi^k,\psi} - \sigma_2 \mathcal{N} \bigbracket{\phi - \overline{\phi}}\|^2+ 1)
\notag\\
&\leq C(\|\widehat{\mhi}\|^2+1),
\notag
\end{align}
which further implies
$$
\norm{\p \widetilde{E}_\phi(\phi)} = \|-\Delta \phi + F_{\phi}'(\phi)\|
\leq \|\Delta\phi\|+ \|F_{\phi}'(\phi)\| \leq C(\|\widehat{\mhi}\|+1).
$$
Estimates for $\fsi'(\psi)$ and $\p \widetilde{E}_{\psi}$ can be obtained in a similar manner.
\epf
\medskip

Let us set (cf. \cite{ADG2013ndm})
\begin{align*}
        & X \stackrel{\rm{def}}{=}  \hssigma{1}  \times
        \mathcal{D}(\p \widetilde{E}_{\phi,a}) \times \mathcal{D}(\p \widetilde{E}_{\psi,b}) \times
        W_{(0)} \times W_{(0)},
        \\
        & Y \stackrel{\rm{def}}{=} \bigbracket{\hssigma{1}}^* \times
        H_{(0)} \times H_{(0)} \times H_{(0)} \times H_{(0)},
        \\
        & \widetilde{X}\stackrel{\rm{def}}{=}  \hssigma{1}  \times
          (\hs{2-s}\cap H_{(a)}) \times (\hs{2-s}\cap H_{(b)}) \times W_{(0)} \times W_{(0)}, \quad  s \in \Bigbracket{0,\frac{1}{4}},\\
        & \widetilde{Y} \stackrel{\rm{def}}{=}  \bm{L}^\frac{3}{2}(\Omega)  \times  (\wkp{1}{\frac{3}{2}}\cap H_{(0)}) \times  V_{(0)}  \times  (\wkp{1}{\frac{3}{2}} \cap H_{(0)}) \times V_{(0)},
\end{align*}
for some given $a\in (-1,1)$, $b\in (0,1)$.
Then we have\medskip
\bl \label{sdiscretee}
Let $\sp^k \in \ls$ and $\phi^k, \psi^k \in \hs{2}$ be given, with $(\phi^k, \psi^k) \in [-1,1] \times [0,1]$ almost everywhere in $\Omega$ and
$\overline{\phi^k} \in (-1,1)$, $\overline{\psi^k} \in (0,1)$.
If $h>0$ is sufficiently small or $\sigma_1 =0$, then there is some $\bigbracket{\sp,\phi,\psi,\widehat{\mhi},\widehat{\msi}} \in X$ (with $a= -h\overline{\sigma_1(\phi^k)}(\overline{\phi^k}-c)+ \overline{\phi^k}$, $b=\overline{\psi^k}$)
that solves the time-discrete problem \eqref{nsdiscrete1}--\eqref{chdiscrete4} and in addition, it satisfies the following discrete energy inequality:
\begin{align}
    &\eto \bigbracket{\sp,\phi,\psi} + \intomega{\rho^k \frac{\bigabs{\sp - \sp^k}^2}{2}}{x}{\Omega} + \half \Bigbracket{\bignorm{\nabla \phi -\nabla \phi^k}^2+ \beta \bignorm{\nabla \psi -\nabla \psi^k}^2} \nonumber \\
    & \qquad +\halfof{\sigma_2} \bignorm{\phi - \phi^k - \big(\overline{\phi} - \overline{\phi^k} \big)}_*^2 \nonumber \\
    & \qquad + h \intomega{\left(\nu \bigbracket{\phi^k} |D\sp|^2 + \mbhi \bigbracket{\phi^k} \abs{\nabla \widehat{\mhi}}^2 + \mbsi\bigbracket{\psi^k} \abs{\nabla \widehat{\msi}}^2\right)}{x}{\Omega} \nonumber \\
    & \qquad + h \bigbracket{\overline{\phi^k} -c} \intomega{\Bigbracket{\widehat{\mhi} - \frac{\rho_1 - \rho_2}{4} \abs{\sp}^2} \sigma_1 \bigbracket{\phi^k}}{x}{\Omega} \nonumber\\
    & \qquad - |\Omega| \bigbracket{ \overline{\phi}-\overline{\phi^k} } \bigbracket{ \overline{F'_\phi (\phi)} + \overline{G_\phi (\phi,\phi^k,\psi)} }  \nonumber \\
    & \quad \le \eto \bigbracket{\sp^k,\phi^k,\psi^k}.
    \label{energy_projected_discrete}
\end{align}
\el

\bpf
First, we show that every $(\sp,\phi,\psi,\widehat{\mhi},\widehat{\msi}) \in X$ solving problem \eqref{nsdiscrete1}--\eqref{chdiscrete4} satisfies the energy inequality
\eqref{energy_projected_discrete}.
By the same calculations as in \cite{ADG2013ndm}, we have
\begin{align*}
    & \intomega{\Bigbracket{(\div \bmj) \halfof{\sp}+ \bracket{\bmj \cdot \nabla}\sp}\cdot \sp}{x}{\Omega} = 0,\\
    & \intomega{\Bigbracket{\div (\rho^k \sp \otimes \sp)- (\nabla \rho^k \cdot \sp)\halfof{\sp}}\cdot \sp}{x}{\Omega}=0,\\
    & \bigbracket{\rho \sp - \rho^k \sp^k}\cdot \sp = \Bigbracket{\rho \halfof{|\sp|^2} - \rho^k \halfof{|\sp^k|^2}} + (\rho -\rho^k)\halfof{|\sp|^2} + \rho^k \halfof{|\sp - \sp^k|^2}.
\end{align*}
Testing \eqref{nsdiscrete2} by $\bmtheta=\mathbf{u}$ and using the above identities, we can deduce that
\begin{align}
    & \intomega{\frac{\rho |\sp|^2-\rho^k |\sp^k|^2}{2h}}{x}{\Omega} + \intomega{\rho^k \frac{\bigabs{\sp - \sp^k}^2}{2h}}{x}{\Omega} + \intomega{\nu\bigbracket{\phi^k}|D \sp|^2}{x}{\Omega} \nonumber \\
    & \quad - \intomega{ \widehat{\mhi} \bigbracket{\nabla \phi^k \cdot \sp} }{x}{\Omega}- \intomega{ \widehat{\msi} \bigbracket{\nabla \psi^k \cdot \sp} }{x}{\Omega}=  \gamma \Bigbracket{\overline{\phi^k} -c} \intomega{\frac{\abs{\sp}^2}{2} \sigma_1\bigbracket{\phi^k}}{x}{\Omega}.
    \label{onestepnstest}
\end{align}
Testing \eqref{chdiscrete1} with $\widehat{\mhi}$, we obtain
\begin{align*}
    &\intomega{\frac{\phi - \phi^k}{h} \widehat{\mhi}}{x}{\Omega} + \intomega{\bigbracket{\sp\cdot \nabla \phi^k}\widehat{\mhi}}{x}{\Omega} +\intomega{\mbhi \bigbracket{\phi^k} \abs{\nabla \widehat{\mhi}}^2}{x}{\Omega} \\
    & \quad = - \Bigbracket{\overline{\phi^k} -c}  \intomega{\sigma_1\bigbracket{\phi^k} \widehat{\mhi}}{x}{\Omega}.
\end{align*}
Similarly, it holds
\begin{equation}\notag
    \intomega{\frac{\psi - \psi^k}{h} \widehat{\msi}}{x}{\Omega} + \intomega{\bigbracket{\sp\cdot \nabla \psi^k}\widehat{\msi}}{x}{\Omega} + \intomega{\mbsi \bigbracket{\psi^k} \abs{\nabla \widehat{\msi}}^2}{x}{\Omega} = 0.
\end{equation}
Next, testing \eqref{chdiscrete2} and \eqref{chdiscrete4} by $\frac{\phi-\phi^k}{h}$ and $\frac{\psi-\psi^k}{h}$, respectively, adding the resultants together, we obtain
\begin{align}
    &\frac{1}{h}\bigbracket{\biginner{\nabla \phi}{\nabla \phi - \nabla \phi^k}+\beta \biginner{\nabla \psi}{\nabla \psi - \nabla \psi^k}}- \intomega{\Bigbracket{ \widehat{\mhi} \frac{\phi-\phi^k}{h} + \widehat{\msi} \frac{\psi-\psi^k}{h}} }{x}{\Omega} \nonumber \\
    & \qquad + \intomega{\frac{G(\phi,\psi)-G(\phi^k,\psi^k)}{h}}{x}{\Omega} + \intomega{\Bigbracket{F_{\phi}'(\phi)\frac{\phi-\phi^k}{h} + F_{\psi}'(\psi)\frac{\psi-\psi^k}{h}}}{x}{\Omega} \nonumber \\
    & \qquad + \sigma_2 \Biginner{\mathcal{N}\bigbracket{\phi - \overline{\phi}}}{\frac{\phi - \phi^k}{h}} \notag \\
    & \quad = \frac{1}{h} |\Omega| \bigbracket{\overline{\phi}-\overline{\phi^k} } \bigbracket{ \overline{F'_\phi (\phi)} + \overline{G_\phi (\phi,\phi^k,\psi)}},
    \label{onestepchemicaltest}
\end{align}
%
Collecting the above identities, we arrive at the energy inequality \eqref{energy_projected_discrete} by using the following facts due to the convexity of $\fhi$, $\fsi$:
\begin{align*}
    & F_{\phi}'(\phi)\bigbracket{\phi-\phi^k} \ge F_{\phi}(\phi) - F_{\phi}\bigbracket{\phi^k},\quad  F_{\psi}'(\psi)\bigbracket{\psi-\psi^k} \ge F_{\psi}(\psi) - F_{\psi}\bigbracket{\psi^k},
\end{align*}
and
\begin{equation*}
    \biginner{\mathcal{N}(\phi - \overline{\phi})}{\phi - \phi^k} = \half (\bignorm{\phi - \overline{\phi}}_*^2 - \bignorm{\phi^k - \overline{\phi^k}}_*^2 + \bignorm{\phi - \phi^k - \overline{\phi} + \overline{\phi^k}}_*^2).
\end{equation*}

We now proceed to show the existence of a solution to the time-discrete problem \eqref{nsdiscrete1}--\eqref{chdiscrete4} by the Leray--Schauder principle.

For $\mathbf{w} = (\sp,\phi,\psi,\widehat{\mhi},\widehat{\msi})\in X$, we define $\mathcal{L}_k:X\to Y$ by setting
    \begin{eqnarray*}
        \mathcal{L}_k(\mathbf{w}) &=&
        \left(
        \begin{array}{c}
            L_k(\sp)\\
            - \div \bigbracket{ \mbhi \bigbracket{\phi^k} \nabla \widehat{\mhi}} \\
            - \Delta \phi + P_0 F_{\phi}'(\phi) \\
            - \div \bigbracket{ \mbsi \bigbracket{\psi^k} \nabla \widehat{\msi}} \\
            - \beta \Delta \psi + P_0 F_{\psi}'(\psi)
        \end{array}
        \right),
    \end{eqnarray*}
where
$$\inner{L_k(\sp)}{\bmtheta} = \intomega{\nu\bigbracket{\phi^k}D \sp : \nabla \bmtheta}{x}{\Omega},\quad \forall\, \bmtheta \in \hssigma{1}.$$
Next, we define $\mathcal{F}_k: \widetilde{X} \to Y$ as
    \begin{eqnarray*}
        \mathcal{F}_k(\mathbf{w})&=&
        \left(
        \begin{array}{c}
            \Big[ -\dfrac{\rho \sp - \rho^k \sp^k}{h}-\div \bracket{\rho^k \sp \otimes \sp}+ \widehat{\mhi} \nabla \phi^k+ \widehat{\msi} \nabla \psi^k
            \\
            -\Bigbracket{ \div \mathbf{J} - \dfrac{\rho-\rho^k}{h}-\sp \cdot \nabla \rho^k - \gamma  \sigma_1\bigbracket{\phi^k}\Bigbracket{\overline{\phi^k} -c} } \dfrac{\sp}{2} - \bracket{\mathbf{J}\cdot \nabla}\sp \Big]
            \\
             -\dfrac{\phi -\phi^k}{h} - \sp \cdot \nabla \phi^k - \sigma_1 (\phi^k) \Bigbracket{\overline{\phi^k} -c}
            \\
            \widehat{\mhi} - P_0 G_{\phi} \bigbracket{\phi,\phi^k,\psi} - \sigma_2 \mathcal{N} \bigbracket{ \phi - \overline{\phi}}
            \\
            -\dfrac{\psi -\psi^k}{h} - \sp \cdot \nabla \psi^k
            \\
            \widehat{\msi} - P_0 G_{\psi} \bigbracket{\phi^k,\psi,\psi^k}
        \end{array}
        \right),
    \end{eqnarray*}
where the first element (the first and second lines) should be understood in the sense of weak formulation (see \eqref{nsdiscrete2}).
Then $\mathbf{w}\in X$ is a weak solution of the time-discrete problem \eqref{nsdiscrete1}--\eqref{chdiscrete4} if and only if
\begin{equation}
\label{FP}
    \mathcal{L}_k(\mathbf{w}) - \mathcal{F}_k(\mathbf{w})=0.
\end{equation}

Applying the same argument as in \cite[Section 4]{ADG2013ndm} with minor modifications, we can show that $\mathcal{L}_k:X\to Y$ is invertible with the inverse $\mathcal{L}_k^{-1}:Y \to X$,
which is continuous from $Y$ to $\widetilde{X}$.
Since the embedding $i_{\widetilde{Y}}:\widetilde{Y} \hookrightarrow Y$ is compact, it follows that $\mathcal{L}_k^{-1} \circ i_{\widetilde{Y}}$ is a compact mapping from $\widetilde{Y}$ to $\widetilde{X}$, still denoted by $\mathcal{L}_k^{-1}$ in the subsequent analysis.

Next, since $G_{\phi}$ and $G_{\psi}$ belong to $BUC^1(\rd{3})$, we can recover the estimates in \cite[Section 4]{ADG2013ndm} to conclude that $\mathcal{F}_k:\widetilde{X}\to \widetilde{Y}$ is continuous and maps bounded sets in $\widetilde{X}$ into bounded sets in $\widetilde{Y}$. The related details are omitted here.

Define $\mathbf{f}=\mathcal{L}_k(\mathbf{w})$, we can rewrite the equation \eqref{FP} as
\begin{equation*}
    \mathbf{f} - \mathcal{F}_k \circ \mathcal{L}_k^{-1}(\mathbf{f}) = \mathbf{0}.
\end{equation*}
Thus, the problem reduces to finding a fixed point of the compact mapping
$$\mathcal{K}_k\stackrel{\rm{def}}{=}\mathcal{F}_k \circ \mathcal{L}_k^{-1}:\widetilde{Y}\to  \widetilde{Y}.$$
By the Leray--Schauder principle, it suffices to show that
\begin{equation*}
        \exists\, R >0 \ \mathrm{such \ that, \ if} \ \mathbf{f}\in \widetilde{Y} \ \mathrm{and} \ 0\le \lambda \le 1 \ \mathrm{fulfills} \ \mathbf{f}=\lambda \mathcal{K}_k \mathbf{f}, \ \mathrm{then} \ \norm{\mathbf{f}}_{\widetilde{Y}} \le R.
\end{equation*}
Recalling that $\mathbf{w} = \mathcal{L}_k^{-1}(\mathbf{f})$, we have
\begin{equation*}
    \mathbf{f}=\lambda \mathcal{K}_k \mathbf{f} \quad \Longleftrightarrow
    \quad \mathcal{L}_k(\mathbf{w}) - \lambda \mathcal{F}_k(\mathbf{w})=\mathbf{0},
\end{equation*}
which is equivalent to the following weak formulation
\begin{align}
    & \lambda\Biginner{\frac{\rho \sp - \rho^k \sp^k}{h}}{\bmtheta} + \lambda\biginner{\div \bigbracket{\rho^k \sp \otimes \sp}}{\bmtheta}
    + \biginner{\nu\bigbracket{\phi^k} D \sp}{\nabla  \bmtheta} +\lambda\inner{\bracket{\mathbf{J}\cdot \nabla}\sp}{\bmtheta}
    \nonumber \\
    & \qquad + \lambda \Biginner{\Bigbracket{\div \mathbf{J} - \frac{\rho-\rho^k}{h}-\sp \cdot \nabla \rho^k - \gamma  \sigma_1 \bigbracket{\phi^k} \Bigbracket{\overline{\phi^k} -c} }\halfof{\sp}}{\bmtheta}
    \nonumber \\
    & \quad = \lambda \bigbracket{\biginner{\widehat{\mhi} \nabla \phi^k}{\bmtheta}+\biginner{\widehat{\msi} \nabla \psi^k}{\bmtheta}},\quad
    \forall\, \bmtheta \in C^{\infty}_{0,\sigma}(\Omega),
    \label{lambdansdiscrete2} \\
    & \mathbf{J} = -\gamma \mbhi \big(\phi^k \big) \nabla \widehat{\mhi},  \qquad \mathrm{a.e.\ in} \ \Omega, \\
    & \lambda \frac{\phi -\phi^k}{h} + \lambda \sp \cdot \nabla \phi^k + \lambda  \sigma_1 \bigbracket{\phi^k} \Bigbracket{ \overline{\phi^k} - c } = \div \bigbracket{ \mbhi\bigbracket{\phi^k} \nabla \widehat{\mhi} }, \qquad  \mathrm{a.e.\ in} \ \Omega,
    \label{lambdachdiscrete1}\\
    & \lambda \widehat{\mhi} - \lambda P_0 G_{\phi}\bigbracket{\phi,\phi^k,\psi} - \lambda \sigma_2 \mathcal{N} \bigbracket{ \phi - \overline{\phi} } = -\Delta \phi + P_0 F_{\phi}'(\phi), \qquad \quad \mathrm{a.e.\ in} \ \Omega,
    \label{lambdachdiscrete2}\\
    & \lambda \frac{\psi -\psi^k}{h} + \lambda \sp \cdot \nabla \psi^k = \div \bigbracket{ \mbsi \bigbracket{\psi^k} \nabla \widehat{\msi} },
    \qquad\quad \ \   \mathrm{a.e.\ in} \ \Omega,
    \label{lambdachdiscrete3}\\
    & \lambda \widehat{\msi} - \lambda P_0 G_{\psi}\bigbracket{\phi^k,\psi,\psi^k} = -\beta \Delta \psi + P_0 F_{\psi}'(\psi), \qquad  \mathrm{a.e.\ in} \ \Omega.
    \label{lambdachdiscrete4}
\end{align}
By calculations similar to those for the energy inequality \eqref{energy_projected_discrete}, we obtain
\begin{align*}
    & h \intomega{\left(\nu\bigbracket{\phi^k}\abs{D \sp}^2 + \mbhi\bigbracket{\phi^k} \abs{\nabla \widehat{\mhi}}^2 + \mbsi\bigbracket{\psi^k} \abs{\nabla \widehat{\msi}}^2\right) }{x}{\Omega} + \half \bigbracket{\norm{\gradphi}^2 + \beta \norm{\gradpsi}^2} \\
    & \qquad + \intomega{\bracket{F_{\phi}(\phi)+F_{\psi}(\psi)}}{x}{\Omega} \\
    & \quad \le \intomega{\halfof{\rho^k \bigabs{\sp^k}^2}}{x}{\Omega}  + \half \Bigbracket{\bignorm{\gradphi^k}^2+\beta \bignorm{\gradpsi^k}^2} \\
    & \qquad  + \intomega{\bigbracket{F_{\phi}\bigbracket{\phi^k}+F_{\psi}\bigbracket{\psi^k} +
    \lambda \abs{G(\phi,\psi)} + \lambda \bigabs{G \bigbracket{\phi^k,\psi^k}}}}{x}{\Omega} + \frac{\lambda\sigma_2}{2} \bignorm{\phi^k - \overline{\phi^k}}_*^2
    \nonumber \\
    & \qquad - \lambda \intomega{\halfof{\rho \abs{\sp}^2}}{x}{\Omega}  - \lambda h \Bigbracket{\overline{\phi^k}-c} \intomega{\sigma_1\bigbracket{\phi^k}\Bigbracket{\widehat{\mhi} - \frac{\rho_1 - \rho_2}{4} \abs{\sp}^2}}{x}{\Omega} \\
    & \qquad + |\Omega| \bigbracket{ \overline{\phi} - \overline{\phi^k} } \bigbracket{ \overline{F'_\phi (\phi)} + \lambda \overline{G_\phi (\phi,\phi^k,\psi)} } .
\end{align*}
Let us take the time step $h\in (0,1)$ sufficiently small satisfying (cf. Remark \ref{averagephi})
\begin{equation*}
h\sigma_1^* \le \frac12 \quad \text{and}\quad  h\sigma_1^* | \gamma| \le \frac{1}{2}\min \{ \rho_1,\rho_2 \}.
\end{equation*}
When $\sigma_1=0$, the above restrictions on $h$ are not necessary.
It follows from $\phi \in \mathcal{D}(\p \widetilde{E}_{\phi,a})$ and $\psi \in \mathcal{D}(\p \widetilde{E}_{\psi,b})$ with $a= -h\overline{\sigma_1(\phi^k)}(\overline{\phi^k}-c)+ \overline{\phi^k}\in (-1,1)$, $b=\overline{\psi^k}\in (0,1)$ that $(\phi,\psi) \in [-1,1]\times [0,1]$ almost everywhere in $\Omega$. Besides, as in the proof of Lemma \ref{apriorichemical}, there exists a positive constant $K_{\phi} = K_{\phi}(a, \overline{\phi^k})$ being uniform with respect to $\lambda \in [0,1]$ such that (cf. \eqref{lambdachdiscrete2})
$$
\bigabs{ \overline{F'_\phi(\phi)} } \le K_\phi (1 + \norm{\mu_\phi}).
$$
Then we can infer from H\"{o}lder's inequality, the Poincar\'e--Wirtinger inequality and Young's inequality that
\begin{align*}
     & - \lambda \intomega{\halfof{\rho \abs{\sp}^2}}{x}{\Omega}  - \lambda h \Bigbracket{\overline{\phi^k}-c} \intomega{\sigma_1\bigbracket{\phi^k}\Bigbracket{\widehat{\mhi} - \frac{\rho_1 - \rho_2}{4} \abs{\sp}^2}}{x}{\Omega}
     \\
     &\quad \leq 2h \abs{\Omega}^{\frac{1}{2}} \sigma_1^* \norm{\widehat{\mu_\phi}}
     \\
     &\quad \leq \frac14 h\underline{m_\phi} \|\nabla \widehat{\mu_\phi}\|^2 + \dfrac{4h \abs{\Omega} (\sigma_1^*)^2 C_P^2}{\underline{m_\phi}}.
\end{align*}
Using Young's inequality again, we get
$$ |\Omega| \bigabs{ \overline{\phi}-\overline{\phi^k} } \bigabs{ \overline{F'_\phi (\phi)} } \le \frac14 h \underline{m_\phi} \|\nabla \widehat{\mu_\phi}\|^2 + \dfrac{ 4 |\Omega|^2 C_P^2 K_\phi^2 }{h\underline{m_\phi}} + 2 |\Omega| K_\phi. $$
Therefore, with the help of Korn's inequality for $\sp \in \hssigma{1}$, we have
\begin{equation}\label{estimateschauder}
    \norm{\sp}_{\bm{H}^1(\Omega)} + \norm{\nabla \widehat{\mhi}}+\norm{\nabla \widehat{\msi}} + \norm{\phi}_{\hone}+ \norm{\psi}_{\hone} \le C_k,
\end{equation}
where $C_k$ is a positive constant depending on $h$ but is independent of $\lambda$.
\smallskip

Consider the following Neumann problems (cf. \cite{ADG2013ndm}):
    \begin{equation*}
    \begin{cases}
    -\div \bigbracket{ \mbhi \bigbracket{\phi^k} \nabla \widehat{\mhi}} =g_\phi,\quad & \text{in} \ \Omega,
    \\
    \p_{\mathbf{n}} \widehat{\mhi} =0,\quad & \text{on}\ \partial \Omega,
    \end{cases}
    \end{equation*}
    \begin{equation*}
    \begin{cases}
    -\div \bigbracket{ \mbsi \bigbracket{\psi^k} \nabla \widehat{\msi}} =g_\psi,\quad & \text{in} \ \Omega,
    \\
    \p_{\mathbf{n}} \widehat{\msi} =0,\quad & \text{on}\ \partial \Omega,
    \end{cases}
    \end{equation*}
    with $\overline{\widehat{\mu_\phi}} = \overline{\widehat{\mu_\psi}} = 0$ and $g_\phi, g_\psi$ given by (cf. \eqref{lambdachdiscrete1}, \eqref{lambdachdiscrete3})
    \begin{align*}
    g_\phi &= -\lambda \frac{\phi -\phi^k}{h} - \lambda \sp \cdot \nabla \phi^k - \lambda  \sigma_1 \bigbracket{\phi^k} \Bigbracket{ \overline{\phi^k} - c } \in H_{(0)},
    \\
    g_\psi &= -\lambda \frac{\psi -\psi^k}{h} - \lambda \sp \cdot \nabla \psi^k  \in H_{(0)}.
    \end{align*}
    With \eqref{estimateschauder}, we can apply the bootstrap argument in \cite{ADG2013ndm} to deduce that
   \begin{align*}
        & \norm{\Delta \widehat{\mu_\phi}}
         \le \dfrac{C}{ \underline{m_\phi} } \left(\dfrac{ \norm{\nabla \phi^k}_{\bm{L}^6(\Omega)}^2 }{ \underline{m_\phi}^2 } +1  \right) \norm{g_\phi},
        \\
        & \norm{\Delta \widehat{\mu_\psi}} \le \dfrac{C}{ \underline{m_\psi} } \left(\dfrac{ \norm{\nabla \psi^k}_{\bm{L}^6(\Omega)}^2 }{ \underline{m_\psi}^2 } +1  \right) \norm{g_\psi}.
    \end{align*}
Hence, applying the elliptic regularity theory, we get
\begin{equation*}
     \norm{\widehat{\mhi}}_{\hs{2}}+ \norm{\widehat{\msi}}_{\hs{2}} \le C_k.
\end{equation*}
The above estimate combined with \eqref{lambdachdiscrete2} and \eqref{lambdachdiscrete4} yields $\norm{\p \widetilde{E}_{\phi,a}}+\norm{\p \widetilde{E}_{\psi,b}}\le C_k$.
As in the proof of Lemma \ref{apriorichemical}, we also get $\|\phi\|_{H^2(\Omega)} + \|\psi\|_{H^2(\Omega)} \le C_k$.
Summing up, we obtain
\begin{equation*}
    \norm{\mathbf{w}}_{\widetilde{X}} \le C_k,
\end{equation*}
and as a result, it holds
\begin{equation*}
   \norm{\mathbf{f}}_{\widetilde{Y}}=\norm{\lambda \mathcal{F}_k (\mathbf{w})}_{\widetilde{Y}} \leq C_k(\norm{\mathbf{w}}_{\widetilde{X}} +1)\le C_k,
\end{equation*}
where the constant $C_k>0$ is independent of $\lambda\in [0,1]$.

Hence, the Leray--Schauder principle applies. This establishes the existence of a weak solution to the time-discrete problem \eqref{nsdiscrete1}--\eqref{chdiscrete4} that satisfies the discrete energy inequality \eqref{energy_projected_discrete}. The proof is complete. \epf
\medskip

\begin{remark} \label{discrete_nonproject_wse} \rm
    Let  $\mathbf{w} = (\sp,\phi,\psi,\widehat{\mu_\phi},\widehat{\mu_\psi}) \in X$ be a solution to the problem \eqref{nsdiscrete1}--\eqref{chdiscrete4}. Set
    $$ \mu_\phi = \widehat{\mu_\phi} + \overline{F'_\phi(\phi)} + \overline{G_\phi(\phi,\phi^k,\psi)},
    \quad
    \mu_\psi = \widehat{\mu_\psi} + \overline{F'_\psi(\psi)} + \overline{G_\psi(\phi^k,\psi,\psi^k)}. $$
    It is straightforward to check that $\widetilde{\mathbf{w}} = (\sp,\phi,\psi,\mu_\phi,\mu_\psi)$ solves the problem
    \begin{align}
        & \Biginner{\frac{\rho \sp - \rho^k \sp^k}{h}}{\bmtheta} + \biginner{\div \bigbracket{\rho^k \sp \otimes \sp}}{\bmtheta} + \biginner{\nu \big(\phi^k \big) D \sp}{\nabla \bmtheta} +\biginner{\div \bracket{\sp \otimes \mathbf{J}}}{\bmtheta}
        \nonumber \\
        &\quad = \biginner{\mhi \nabla \phi^k}{\bmtheta}+\biginner{\msi \nabla \psi^k}{\bmtheta}, \qquad \forall\, \bmtheta \in C^{\infty}_{0,\sigma}(\Omega),
        \label{nsdiscrete1_nonproject}\\
        &\mathbf{J} = -\gamma \mbhi \big(\phi^k \big) \nabla \mhi,  \qquad \mathrm{a.e.\ in} \ \Omega,
        \nonumber\\
        & \frac{\phi -\phi^k}{h} + \sp \cdot \nabla \phi^k + \sigma_1\bigbracket{\phi^k}
        \Bigbracket{\overline{\phi^k} -c} = \div \bigbracket{\mbhi \bigbracket{\phi^k} \nabla \mhi}, \qquad \mathrm{a.e.\ in} \ \Omega,
        \label{chdiscrete1_nonproject}\\
        & \mhi - G_{\phi} \bigbracket{\phi,\phi^k,\psi} - \sigma_2 \mathcal{N} \bigbracket{\phi - \overline{\phi}} = -\Delta \phi + F_{\phi}'(\phi), \qquad \mathrm{a.e.\ in} \ \Omega,
        \label{chdiscrete2_nonproject}\\
        & \frac{\psi -\psi^k}{h} + \sp \cdot \nabla \psi^k = \div \bigbracket{\mbsi \bigbracket{\psi^k} \nabla \msi}, \qquad \mathrm{a.e.\ in} \ \Omega,
        \label{chdiscrete3_nonproject}\\
        & \msi - G_{\psi} \bigbracket{\phi^k,\psi,\psi^k} = -\beta \Delta \psi + F_{\psi}'(\psi), \qquad \mathrm{a.e.\ in} \ \Omega,
        \label{chdiscrete4_nonproject}
    \end{align}
    in addition, the following discrete energy inequality holds
    \begin{align}
        &\eto \bigbracket{\sp,\phi,\psi} + \intomega{\rho^k \frac{\bigabs{\sp - \sp^k}^2}{2}}{x}{\Omega} + \half \Bigbracket{\bignorm{\nabla \phi -\nabla \phi^k}^2+ \beta \bignorm{\nabla \psi -\nabla \psi^k}^2} \nonumber \\
        & \qquad +\halfof{\sigma_2} \bignorm{\phi - \phi^k - \big(\overline{\phi} - \overline{\phi^k} \big)}_*^2 \nonumber \\
        & \qquad + h \intomega{\left(\nu \bigbracket{\phi^k} |D\sp|^2 + \mbhi \bigbracket{\phi^k} \abs{\nabla \mhi}^2 + \mbsi\bigbracket{\psi^k} \abs{\nabla \msi}^2\right)}{x}{\Omega} \nonumber \\
        & \qquad + h \bigbracket{\overline{\phi^k} -c} \intomega{\Bigbracket{\mhi - \frac{\rho_1 - \rho_2}{4} \abs{\sp}^2} \sigma_1 \bigbracket{\phi^k}}{x}{\Omega} \nonumber\\
        & \quad \le \eto \bigbracket{\sp^k,\phi^k,\psi^k}.
        \label{energy_discrete}
    \end{align}
\end{remark}

\subsection{Existence of a global weak solution.}\label{proofexistencendm}
\medskip

We are now in a position to prove Theorem \ref{wsnschndme}.
\medskip
\bpf{\bf (Proof of Theorem \ref{wsnschndme}.)}\,
The construction of a global weak solution to problem \eqref{eq:nsch1}--\eqref{eq:nschi} on $[0,\infty)$ can be done by adapting the arguments in \cite{ADG2013ndm}. Let $N\in \mathbb{N}_+$ be given. First of all, we approximate the mobilities $\mbhi$ and $\mbsi$ with
$$ \mbhi^N \stackrel{\rm{def}}{=} \eta^N * \mbhi, \quad \mbsi^N \stackrel{\rm{def}}{=} \eta^N * \mbsi,$$
where $\eta^N$ is a family of positive mollifiers in $\r$ and $*$ stands for the convolution in $\r$.
That is,
$$ \eta^N = N^d \eta(N x), \ x \in \rd{d}, $$
where $\eta \in C_c^{\infty} \bigbracket{\rd{d}} $ satisfies
$$
\eta \ge 0 \ \ \mathrm{and} \ \ \intomega{\eta \bracket{x}}{x}{\rd{d}} =1.
$$
Since $\mbhi$ and $\mbsi$ are globally Lipschitz continuous on $\r$, then $\bigbracket{\mbhi^N,\mbsi^N}$ converges to $(\mbhi,\mbsi)$ in $C(\r)^2$ as $N\to \infty$. Additionally, $\bignorm{\mbhi^N}_{C^1(\r)}$ (resp. $\bignorm{\mbsi^N}_{C^1(\r)}$) is bounded by $\norm{\mbhi}_{W^{1,\infty}(\r)}$ (resp. $\norm{\mbsi}_{W^{1,\infty}(\r)}$) uniformly with respect to $N$, and it holds
$$
\min\limits_{s \in \r} \{ \mbhi^N(s) \} \ge \underline{\mbhi}, \quad  \min\limits_{s \in \r} \{ \mbsi^N(s) \} \ge \underline{\mbsi}.
$$

Next, we approximate the initial data $(\phi_0,\psi_0)\in \hs{1} \times \hs{1}$ by $(\phi_0^N,\psi_0^N) \in W \times W$ in order to apply Lemma \ref{sdiscretee}. As shown in \cite{GP2022}, we know that there exists a sequence $\left\{(\phi_0^N,\psi_0^N)\right\}_{N=1}^{\infty} \subset H^3_N(\Omega) \times H^3_N(\Omega)$ satisfying
\begin{equation*}
    \left\|\phi_0^N \right\|_{L^{\infty}(\Omega)},\ \left\|\psi_0^N \right\|_{L^{\infty}(\Omega)} \le 1- \frac{1}{N}, \quad  \forall\, N \in \mathbb{N_+},
\end{equation*}
and
\begin{equation*}
    \bigbracket{\phi_0^N,\psi_0^N} \to (\phi_0,\psi_0) \ \ \mathrm{in} \ \ V \ \ \mathrm{as} \ \ N \to \infty,
\end{equation*}
where $H^3_N(\Omega)\stackrel{\rm{def}}{=} W \cap H^3(\Omega)$.
In addition, recall that in Remark \ref{averagephi} if $N \ge 1+ \sigma_1^*$, then the function $\phi^{k+1}$ constructed there satisfies
    \begin{align*}
        \Bigabs{ \overline{\phi^{k+1}} - c } & =  \prod\limits_{n=0}^k \Bigbracket{1-\frac{\overline{\sigma_1\bracket{\phi^n}}}{N}} \Bigabs{\overline{\phi^N_0} - c} \\
        & \le e^{- \frac{1}{N} \sum\limits_{n=0}^k \overline{\sigma_1\bracket{\phi^n}}} \Bigabs{\overline{\phi^N_0} - c} \le 2 e^{- \frac{1}{N} \sum\limits_{n=0}^k \overline{\sigma_1\bracket{\phi^n}}}.
    \end{align*}
As a consequence, for the time-discrete problem \eqref{nsdiscrete1_nonproject}--\eqref{chdiscrete4_nonproject} with $h=\frac{1}{N}$, $(\mbhi^N, \mbsi^N)$ as the mobility functions and $(\sp_0,\phi_0^N,\psi_0^N)$ as the initial data, we can apply Lemma \ref{sdiscretee} combined with Remark \ref{discrete_nonproject_wse} successively to construct a sequence of solutions $\big\{\sp^{k+1,N},\phi^{k+1,N},\psi^{k+1,N},\mhi^{k+1,N},\msi^{k+1,N}\big\}_{k=0}^{\infty}$.
\smallskip

In the subsequent analysis, we use the following convention for the sake of simplicity (see \cite{ADG2013ndm}). For a given sequence $\set{f^{k,N}}_{k=0}^{\infty}$ that may depend on $N$, we denote by $f^N$ the function on $[-h,\infty)$ such that
$$f^N(t) = f^{k,N}, \ \ \  t \in [(k-1)h,kh),\quad \forall\, k \in \mathbb{N}.$$
In addition, we use $\widetilde{f}^N$ to denote the piecewise linear interpolant of $f^N(t^{k,N})$, where $t^{k,N} = kh, \ k \in \mathbb{N}$, that is,
$$ \widetilde{f}^N = \frac{1}{h} \chi_{[0,h]} *_t f^N. $$
Besides, for a given function $f=f(t)$, we set
\begin{align*}
    & (\Delta^+_h f)(t) \stackrel{\rm{def}}{=} f(t+h)-f(t), \qquad (\Delta^-_h f)(t) \stackrel{\rm{def}}{=} f(t)-f(t-h),\\
    & \p^+_{t,h} f(t) \stackrel{\rm{def}}{=} \frac{1}{h} (\Delta^+_h f)(t), \qquad\quad \ \ \ \p^-_{t,h} f(t) \stackrel{\rm{def}}{=} \frac{1}{h} (\Delta^-_h f)(t),\\
    & f_h \stackrel{\rm{def}}{=} (\tau^*_h f)(t) = f(t-h).
\end{align*}
For any $\bmtheta \in C^{\infty}_0(Q)^d$ satisfying $\div \bmtheta =0$, we choose
$$\bmtheta^k \stackrel{\rm{def}}{=} \intinterval{\bmtheta(t)}{t}{kh}{(k+1)h}$$
as the test function in \eqref{nsdiscrete1_nonproject} and sum over $k \in \mathbb{N}$ to obtain
\begin{align}
    & \intinterval{\biginner{\pn (\rho^N \sp^N)}{\bmtheta}}{t}{0}{\infty} - \intinterval{\biginner{ \rho^N_h \sp^N \otimes \sp^N}{\nabla \bmtheta}}{t}{0}{\infty} \nonumber \\
    &\qquad + \intinterval{\biginner{\nu \bigbracket{\phi^N_h} D \sp^N}{\nabla  \bmtheta}}{t}{0}{\infty} -\intinterval{\biginner{\bigbracket{\sp^N \otimes \mathbf{J}^N}}{\nabla \bmtheta}}{t}{0}{\infty} \nonumber \\
    & \quad = \intinterval{\biginner{\mhi^N \nabla \phi^N_h}{\bmtheta}+\biginner{\msi^N \nabla \psi^N_h}{\bmtheta}}{t}{0}{\infty},\label{nsglobaldiscrete}
\end{align}
where
$$
\rho^N= \gamma \phi^N + \frac{ \rho_1 + \rho_2 }{2}\quad \text{and}\quad \bmj^N = - \gamma \mbhi^N \bigbracket{\phi^N_h} \nabla \mhi^N.
$$
Analogously, we obtain
\begin{align}
    & -\intinterval{\biginner{\pn \phi^N}{\zeta}}{t}{0}{\infty} + \intinterval{\biginner{\sp^N \phi_h^N}{\nabla \zeta} }{t}{0}{\infty} - \intinterval{\Bigbracket{ \overline{\phi^N_h} - c } \biginner{\sigma_1\bigbracket{\phi^N_h}}{\zeta} }{t}{0}{\infty} \nonumber \\
    & \quad = \intinterval{\biginner{\mbhi^N \bigbracket{\phi^N_h} \nabla \mhi^N}{\nabla \zeta}}{t}{0}{\infty},\label{chglobaldiscrete1}\\
    & -\intinterval{\biginner{\pn \psi^N}{\zeta}}{t}{0}{\infty} + \intinterval{\biginner{\sp^N \psi_h^N}{\nabla \zeta}}{t}{0}{\infty}= \intinterval{\biginner{\mbsi^N \bigbracket{\psi^N_h} \nabla \msi^N}{\nabla \zeta}}{t}{0}{\infty},\label{chglobaldiscrete3}
\end{align}
for all $\zeta \in C^{\infty}_0\big((0,\infty);C^1\big(\overline{\Omega} \big)\big)$, and
\begin{align}
    & \mhi^N - G_{\phi} \bigbracket{\phi^N,\phi^N_h,\psi^N}
    -\sigma_2 \mathcal{N} \bigbracket{ \phi^N - \overline{\phi^N} }
    = -\Delta \phi^N + F_{\phi}' \bigbracket{\phi^N}, \label{chglobaldiscrete2}\\
    & \msi^N - G_{\psi} \bigbracket{\phi^N_h,\psi^N,\psi^N_h} = -\beta \Delta \psi^N + F_{\psi}' \bigbracket{\psi^N}, \label{chglobaldiscrete4}
\end{align}
almost everywhere in $Q$.
Let us define, for all $t \in (t^{k,N},t^{k+1,N}), \ k \in \mathbb{N}$, the following quantities
\begin{align*}
    Diss^N(t) & \stackrel{\rm{def}}{=} \intomega{\left( \mbhi^N\bigbracket{\phi^{k,N}} \bigabs{\nabla \mhi^{k+1,N}}^2 + \mbsi^N\bigbracket{\psi^{k,N}} \bigabs{\nabla \msi^{k+1,N}}^2 \right)}{x}{\Omega} \\
    & \quad + \intomega{ \nu\bigbracket{\phi^{k,N}} \bigabs{D \sp^{k+1,N}}^2 }{x}{\Omega}, \\
    R^N(t) & \stackrel{\rm{def}}{=} \bigbracket{\overline{\phi^{k,N}} -c} \intomega{ \Bigbracket{\mhi^{k+1,N} - \frac{\gamma}{2} \abs{\sp^{k+1,N}}^2 } \sigma_1\bigbracket{\phi^{k,N}} }{x}{\Omega}.
\end{align*}%
Setting $E^{k,N}= \eto \bigbracket{\sp^{k,N},\phi^{k,N},\psi^{k,N}}$, $k \in \mathbb{N}$, and using $\widetilde{E}^N$ to denote the corresponding linear interpolant of $E^N(t) = \eto \bigbracket{\sp^N(t),\phi^N(t),\psi^N(t)}$ at $t^{k,N}$, from the discrete energy inequality \eqref{energy_discrete} we find that
\begin{equation}\label{diffdiss}
    -\frac{\mathrm{d}}{\mathrm{dt}} \widetilde{E}^N (t) \ge Diss^N(t) + R^N(t),
\end{equation}
for all $t \in (t^{k,N},t^{k+1,N}), \ k \in \mathbb{N} $.

Multiplying the above inequality by $g \in W^{1,1}(0,\infty) \cap C_c ([0,\infty))$ with $g \ge 0$, integrating over $(0,\infty)$, and applying integration by parts, we obtain
\begin{equation*}
    \eto \bigbracket{\sp_0,\phi^N_0,\psi^N_0} g(0) + \intinterval{\widetilde{E}^N (t) g'(t) }{t}{0}{\infty} \ge \intinterval{ D^N(t)g(t) }{t}{0}{\infty},
\end{equation*}
where we denote
$$D^N(t) \stackrel{\rm{def}}{=} Diss^N(t) + R^N(t).$$
Integrating \eqref{diffdiss} over $(s,t)$, $s,t \in h \mathbb{N}$, gives
\begin{align*}
    &\eto \bigbracket{\sp^N(t),\phi^N(t),\psi^N(t)} + \intomega{\left(\mbhi^N \bigbracket{\phi^N_h} \bigabs{\nabla \mhi^N }^2 + \mbsi^N \bigbracket{\psi^N_h} \bigabs{\nabla \msi^N }^2\right)}{(x,\tau)}{Q_{(s,t)}} \\
    & \qquad + \intomega{\nu(\phi^N_h) \bigabs{D \sp^N}^2}{(x,\tau)}{Q_{(s,t)}} +
        \intomega{ \bigbracket{\overline{\phi^N_h} -c} \Bigbracket{\mhi^N - \frac{\gamma}{2} \abs{\sp^N}^2 } \sigma_1\bigbracket{\phi^N_h}}{\bracket{x,\tau}}{Q_{(s,t)}}
     \\
    & \quad \le \eto \bigbracket{\sp^N(s),\phi^N(s),\psi^N(s)}.
\end{align*}
From Lemma \ref{apriorichemical}, it holds
\begin{equation*}
    \bigabs{R^N} \le \half \intomega{\mbhi^N \bigbracket{\phi^N_h} \bigabs{\nabla \mhi^N }^2}{x}{\Omega} + C \bigbracket{ 1 + \eto \bigbracket{\sp^N,\phi^N,\psi^N} + \eto \bigbracket{\sp^N_h,\phi^N_h,\psi^N_h} },
\end{equation*}
where $C>0$ depends on $\sigma_1$, $\Omega$, $K_{\phi}(\delta)$, $\gamma$, $\mbhi$ and coefficients of the system \eqref{eq:nsch1}--\eqref{eq:nsch6} but is independent of $N$. Here $\delta \in (0,1)$ is a constant such that $c,\overline{\phi^N_0} \in (-1+\delta,1-\delta)$ for all $ N \in \mathbb{N}_+$. Therefore, $E^N(t)$ satisfies
\begin{equation*}
    E^N(t) \le E^N(s) + C\intinterval{\big(1 + E^N(\tau) + E^N(\tau -h)\big)}{\tau}{s}{t},
\end{equation*}
for all $0\le s \le t < \infty$ with $s,t \in h \mathbb{N}$. In particular, we have
\begin{equation*}
    (1-Ch)E^N((k+1)h) \le Ch + (1+Ch)E^N(kh), \quad \forall\, k \in \mathbb{N}.
\end{equation*}
Without loss of generality, we set $h>0$ sufficiently small that $Ch < 1/2$, then it follows from an iteration of the above inequality that
\begin{equation*}
    E^N(t) \le e^{2Ct} \Bigbracket{E^N(0)+\frac{1}{2}}, \quad  \forall\, t \ge 0.
\end{equation*}
From the uniform boundedness of $\eto \bigbracket{\sp_0,\phi^N_0,\psi^N_0}$ with respect to $N$, it follows that, for any given $T>0$, the following terms are uniformly bounded with respect to $N$,
in the corresponding spaces:
\begin{align*}
    & \sp^N \in L^2\bigbracket{0,T;\hssigma{1}} \cap L^{\infty} \bigbracket{0,T;\ls},\\
    & \nabla \mhi^N, \nabla \msi^N \in L^2 \bigbracket{0,T;\bm{L}^2(\Omega)},\\
    & \phi^N, \psi^N \in L^{\infty} \bigbracket{0,T;V}.
\end{align*}
In addition, we have
$$
\intinterval{\Big( \Bigabs{\overline{\mhi^N}} + \Bigabs{\overline{\msi^N}} \Big)}{t}{0}{T} \le C(T),
$$
for a monotonic function $C:\r^+ \to \r^+$.

Using a diagonal argument as in \cite{BF2013}, we are able to extract a convergent subsequence such that\smallskip
\begin{itemize}
    \item $\sp^N$ converges weakly in $L^2 \bigbracket{0,L;\hssigma{1}}$ and weakly-$^*$ in $L^{\infty}\bigbracket{0,L;\ls}$,\smallskip
    \item $\phi^N,\phi^N$ converge weakly-$^*$  in $L^{\infty} \bigbracket{0,L;V}$,\smallskip
    \item $\mhi^N,\msi^N$ converge weakly in $L^2 \bigbracket{0,L;V}$,\smallskip
\end{itemize}
for all $L \in \mathbb{N}_+$.
Hence, there exists a limit $(\sp,\phi,\psi,\mhi,\msi)$ on $[0,\infty)$ such that
\begin{align*}
    &\sp^N \rightharpoonup \sp \ &&\mathrm{in} \ L^2\bigbracket{0,T;\hssigma{1}},\\
    &\sp^N \mathop{\rightharpoonup}\limits^* \sp \ &&\mathrm{in} \ L^{\infty}\bigbracket{0,T;\ls},\\
    &\phi^N \mathop{\rightharpoonup}\limits^* \phi, \ \psi^N \mathop{\rightharpoonup}\limits^* \psi \ &&\mathrm{in} \ L^{\infty}\bigbracket{0,T;V},\\
    &\mhi^N \rightharpoonup \mhi, \ \msi^N \rightharpoonup \msi \ &&\mathrm{in} \ L^2\bigbracket{0,T;V},
\end{align*}
for any $T>0$. Hereafter, all limits are meant to be for suitable subsequences $N_k \to \infty$ (resp. $h_k \to 0$) for $k \to \infty$ unless otherwise stated.

For $\big(\widetilde{\phi}^N,\widetilde{\psi}^N\big)$, we have $\pdif{\widetilde{\phi}^N}{t} = \pn \phi^N$, $\pdif{\widetilde{\psi}^N}{t} = \pn \psi^N$ and
\begin{align}
    \bignorm{\widetilde{\phi}^N - \phi^N }_{V^*} \le h \bignorm{\pdif{\widetilde{\phi}^N}{t}}_{V^*},  \qquad
    &\bignorm{\widetilde{\psi}^N - \psi^N }_{V^*} \le h \bignorm{\pdif{\widetilde{\psi}^N}{t}}_{V^*}.
    \label{interpolationerror}
\end{align}
On account of the boundedness of $\sp^N \phi^N_h$, $\sp^N \psi^N_h$, $\mbhi^N \bigbracket{\phi^N_h} \nabla \mhi^N$, $\mbsi^N \bigbracket{\psi^N_h} \nabla \msi^N$
in $L^2 \bigbracket{0,T; \bm{L}^2(\Omega)}$, we infer from \eqref{chglobaldiscrete1} and \eqref{chglobaldiscrete3} that $\pdif{\widetilde{\phi}^N}{t}, \pdif{\widetilde{\psi}^N}{t} \in L^2 \bigbracket{0,T;V^*} $ are uniformly bounded with respect to $N$.
Together with the boundedness of $\widetilde{\phi}^N,\widetilde{\psi}^N$ in $L^{\infty}(0,T;V)$, the Aubin--Lions--Simon lemma (see \cite{Si1987}) enables us to obtain
the strong convergence
\begin{equation*}
    \widetilde{\phi}^N \to \widetilde{\phi}, \quad \widetilde{\psi}^N \to \widetilde{\psi} \ \ \mathrm{in} \ \ C \bigbracket{[0,T];H},
\end{equation*}
for any $T>0$ for some $\tphi,\tpsi \in C \bigbracket{[0,\infty);H}$.
In particular, up to a subsequence, we also have $\Bigbracket{\tphi^N,\tpsi^N} \to \Bigbracket{\tphi,\tpsi}$ almost everywhere in $\Omega \times (0,\infty)$.
On the other hand, thanks to \eqref{interpolationerror}, it holds
\begin{equation*}
    \tphi^N - \phi^N \to 0, \quad \tpsi^N - \psi^N \to 0 \ \ \mathrm{in} \ \ L^2 \bigbracket{0,T;V^*},
\end{equation*}
which implies $\Bigbracket{\tphi,\tpsi} = (\phi,\psi)$. Furthermore, we observe that
\begin{equation*}
    \tphi^N,\tpsi^N \in H^1 \bigbracket{0,T;V^*} \cap L^2 \bigbracket{0,T;V} \hookrightarrow BUC \bigbracket{[0,T];H}
\end{equation*}
and
\begin{equation*}
    \tphi^N,\tpsi^N \in L^{\infty}\bigbracket{0,T;V}
\end{equation*}
are uniformly bounded with respective to $N$. Therefore, applying \cite[Lemma 2.1]{ADG2013ndm} we can conclude
$$\phi,\psi \in C_w \bigbracket{[0,\infty);V}.$$ From the convergence $\widetilde{\phi}^N \to \phi$ and $\widetilde{\psi}^N \to \psi$ in $C \big( [0,T];H \big)$ for any $T \in (0,\infty)$,
it follows that $$\tphi^N(0) \to \phi(0),\quad \tpsi^N(0) \to \psi(0)\ \ \ \text{in}\ H.$$
Since $\Bigbracket{\tphi^N(0),\tpsi^N(0)} = \bigbracket{\phi^N_0,\psi^N_0}$, we find that $(\phi(0),\psi(0)) = (\phi_0,\psi_0)$ after passing to the limit as $N\to \infty$.

To show the weak convergence of \eqref{chglobaldiscrete2} and \eqref{chglobaldiscrete4}, we observe that the left-hand side of \eqref{chglobaldiscrete2} (resp. \eqref{chglobaldiscrete4}) converges weakly in $L^2 \bigbracket{0,T;H}$ to
\begin{equation*}
    f_{\phi} \stackrel{\rm{def}}{=} \mhi - \pdif{G}{\phi}(\phi,\psi) -\sigma_2 \mathcal{N} \bigbracket{ \phi - \overline{\phi} }\ \ \ (\text{resp.}\ f_{\psi} \stackrel{\rm{def}}{=} \msi - \pdif{G}{\psi}(\phi,\psi) ),
\end{equation*}
which implies
$$\p \widetilde{E}_\phi(\phi^N) \rightharpoonup f_{\phi}\quad (\text{resp.}\ \p \widetilde{E}_\psi(\psi^N) \rightharpoonup f_{\psi}).$$
Hence, it suffices to show that
\begin{align}
    &\mathop{\mathrm{limsup}}\limits_{N \to \infty} \biginner{\p \widetilde{E}_\phi\bigbracket{\phi^N}}{\phi^N}_{Q_T} \le \inner{f_{\phi}}{\phi}_{Q_T},\label{limsupsubdifphi} \\
    &\mathop{\mathrm{limsup}}\limits_{N \to \infty} \biginner{\p \widetilde{E}_\psi \bigbracket{\psi^N}}{\psi^N}_{Q_T} \le \inner{f_{\psi}}{\psi}_{Q_T} , \label{limsupsubdifpsi}
\end{align}
to conclude $\p \widetilde{E}_\phi(\phi) = f_{\phi}$ (resp. $\p \widetilde{E}_\psi(\psi) = f_{\psi}$), since $\p \widetilde{E}_\phi$, $\p \widetilde{E}_\psi$ are maximal monotone operators.
Due to the strong convergence of $\phi^N$, $\psi^N$ in $L^2 \bigbracket{0,T;H}$, \eqref{limsupsubdifphi} and \eqref{limsupsubdifpsi} indeed hold with equal sign, because
\begin{align*}
    & \biginner{\p \widetilde{E}_\phi \bigbracket{\phi^N}}{\phi^N}_{Q_T} = \biginner{f^N_{\phi}}{\phi^N}_{Q_T} \to \inner{f_{\phi}}{\phi}_{Q_T}, \\
    &\biginner{\p \widetilde{E}_\psi \bigbracket{\psi^N}}{\psi^N}_{Q_T} = \biginner{f^N_{\psi}}{\psi^N}_{Q_T} \to \inner{f_{\psi}}{\psi}_{Q_T},
\end{align*}
with the natural definition of $f^N_{\phi},f^N_{\psi}$ as the left-hand side of \eqref{chglobaldiscrete2} and \eqref{chglobaldiscrete4}, respectively.
In addition, using Lemma \ref{apriorichemical}, we can derive uniform estimates of $\phi^N,\psi^N$ in $L^2 \bigbracket{0,T;W}$. Then, from the strong convergence of $\phi^N,\psi^N$ in $L^2 \bigbracket{0,T;V^*}$ and the interpolation, we can obtain the strong convergence
\begin{equation*}
    \phi^N \to \phi, \ \psi^N \to \psi  \ \ \mathrm{in} \ L^2 \bigbracket{0,T;V}.
\end{equation*}
Moreover, by \cite[Lemma 7.4]{GGW2018} and the equations \eqref{chglobaldiscrete2}, \eqref{chglobaldiscrete4}, we can derive the uniform boundedness of $\phi^N$ and $\psi^N$ in $L^2 \bigbracket{0,T;\wkp{2}{p}}$, with $p=6$ if $d=3$ and for any $p\ge2$ if $d=2$, so that
$\phi, \psi\in L^2 \bigbracket{0,T;\wkp{2}{p}}$.

Following the same argument as in \cite{ADG2013ndm}, we find that the following terms are uniformly bounded in corresponding spaces:
\begin{align*}
    \rho^N_h \sp^N \otimes \sp^N \ \ &\mathrm{in} \ L^2 \bigbracket{0,T;\lp{\frac{3}{2}}}, \\
    D \sp^N \ \ &\mathrm{in} \ \lptbig{2}{0,T;\ltwo},\\
    \sp^N \otimes \nabla \mhi^N \ \ & \mathrm{in} \ \lptbig{\frac{8}{7}}{0,T;\lp{\frac{4}{3}}}, \\
    \mhi^N \nabla \phi^N_h, \msi^N \nabla \psi^N_h \ \ & \mathrm{in} \ \lptbig{2}{0,T;\lp{\frac{3}{2}}}.
\end{align*}
By comparison in \eqref{nsglobaldiscrete}, we can deduce that the time derivative $\pdif{\mathbb{P}_{\sigma} \bigbracket{\widetilde{\rho \sp}^N}}{t} $ is uniformly bounded in $\lptbig{\frac{8}{7}}{0,T;\wkpsigma{1}{4}^*}$ with respect to $N$. From the uniform boundedness of $\sp^N$ in $L^{\infty} \bigbracket{ 0,T;\ls } \cap L^2 \bigbracket{ 0,T;\bm{H}^1(\Omega)}$ and the uniform boundedness of $\phi^N$ in $L^2 \bigbracket{ 0,T;\wkp{2}{6} } \cap L^{\infty} \bracket{Q_T}$, we also have the uniform boundedness of $\mathbb{P}_{\sigma} \bigbracket{\widetilde{\rho \sp}^N} $ in $ \lptbig{2}{0,T;\bm{H}^1(\Omega)}$. Then it follows from the Aubin--Lions--Simon lemma that
\begin{equation*}
    \mathbb{P}_{\sigma} (\widetilde{\rho \sp}^N) \to \mathbf{w} \ \mathrm{in} \ \lptbig{2}{0,T;\bm{L}^2(\Omega)},
\end{equation*}
for all $0<T<\infty$, where the distribution $\mathbf{w}$ defined on $Q$ belongs to $L^\infty \bigbracket{0,t;\bm{L}^2(\Omega)}$ for all $t \in \r_+$. Furthermore, we can show that $\mathbf{w} \in C_w \bigbracket{[0,\infty);\bm{L}^2(\Omega)}$. The weak continuity property of $\mathbf{w}$ with respect to time is a direct consequence of \cite[Lemma 2.1]{ADG2013ndm}, because it also belongs to $C \bigbracket{[0,\infty);\wkpsigma{1}{4}^*}$. Since the projection
\begin{equation*}
    \psigma : \lptbig{2}{0,T;\bm{L}^2(\Omega)} \to \lptbig{2}{0,T;\ls}
\end{equation*}
is weakly continuous, it follows from the weak convergence $\widetilde{\rho \sp}^N \rightharpoonup \rho \sp$ in $\lptbig{2}{0,T;\bm{L}^2(\Omega)}$ that $$\mathbf{w} = \psigma (\rho \sp).$$
In a similar manner, with the uniform boundedness of $\psigma \bigbracket{\rho^N \sp^N}$ in $ \lptbig{2}{0,T;\bm{H}^1(\Omega)}$ and $\pn \psigma \bigbracket{\rho^N \sp^N}$ in $\lptbig{\frac{8}{7}}{0,T;\wkpsigma{1}{4}^*}$ with respect to $N$, we can apply a time discrete version of the Aubin--Lions lemma (see \cite[Theorem 1]{DJ12}) to conclude that $\psigma \bigbracket{\rho^N \sp^N} \to \widehat{\mathbf{w}}$ in $\lptbig{2}{0,T;\bm{L}^2(\Omega)}$ for some $\widehat{\mathbf{w}}\in \lptbig{2}{0,T;\bm{L}^2(\Omega)}$. This together with the weak convergence $\rho^N \sp^N \rightharpoonup \rho \sp$ in $\lptbig{2}{0,T;\bm{L}^2(\Omega)}$ yields $ \widehat{\mathbf{w}}=\psigma (\rho \sp)$ so that $\psigma \bigbracket{\rho^N \sp^N} \to \psigma (\rho \sp)$ in $\lptbig{2}{0,T;\bm{L}^2(\Omega)}$.
Combining the strong convergence of $\psigma \bigbracket{\rho^N \sp^N}$ in $L^2 \bigbracket{ 0,T;\bm{L}^2(\Omega)}$ and the weak convergence of $\sp^N$ in $L^2 \bigbracket{ 0,T;\ls }$, using the same arguments as in \cite[Section 5.1]{ADG2013ndm}, we can deduce the strong convergence
$$\sp^N \to \sp\ \ \text{in}\ \lptbig{2}{0,T;\bm{L}^2(\Omega)}$$
and thus the pointwise convergence almost everywhere in $Q$.
For readers' convenience, we provide some details here.
Firstly, we have
\begin{align*}
    \intomega{ \rho^N \bigabs{ \sp^N }^2 }{(x,\tau)}{Q_T} & = \Biginner{  \psigma \bigbracket{\rho^N \sp^N}  }{ \sp^N }_{Q_T} \\
    & \to \biginner{  \psigma (\rho \sp) }{ \sp }_{Q_T} = \intomega{ \rho \abs{ \sp }^2 }{(x,\tau)}{Q_T},\quad \forall\, T >0.
\end{align*}
On the other hand, the strong convergence of $\phi^N$ in $\lptbig{2}{0,T;V}$ and the weak convergence of $\sp^N$ in $\lptbig{2}{0,T;\ls}$ ensure the weak convergence
$$\sqrt{ \rho^N } \sp^N \rightharpoonup \sqrt{ \rho } \sp \ \ \mathrm{in} \ \ \lptbig{2}{0,T;\bm{L}^2(\Omega)}.$$
Therefore, the uniform convexity of $\lptbig{2}{0,T;H}$ yields the strong convergence $$\sqrt{ \rho^N } \sp^N \to \sqrt{ \rho } \sp \ \  \mathrm{in} \ \ \lptbig{2}{0,T;\bm{L}^2(\Omega)}.$$
Then, using again the strong convergence of $\phi^N$ in $\lptbig{2}{0,T;V}$ and the uniform boundedness of $\phi^N$ with respect to $N$ in $\lptbig{\infty}{0,T;V}$, we can conclude
$$ \sp^N = \frac{1}{\sqrt{\rho^N}} \left(\sqrt{\rho^N} \sp^N\right) \to \frac{1}{\sqrt{\rho}} \left(\sqrt{\rho} \sp\right) = \sp \ \  \mathrm{in} \ \ \lptbig{2}{0,T;\bm{L}^2(\Omega)}. $$
Here, we have also used the fact that $\abs{\rho^N} \ge \underline{\rho} = \min\limits_{s \in [-1,1]} \rho(s)>0.$

With the above convergence results, we are allowed to pass to the limit in \eqref{nsglobaldiscrete}, \eqref{chglobaldiscrete1} and \eqref{chglobaldiscrete3} as $N\to\infty$ (up to a suitable subsequence) to conclude that the limit functions $(\sp,\phi,\psi,\mhi,\msi)$ satisfy equations \eqref{testns}, \eqref{testchphi} and \eqref{testchpsi}.
To verify that $(\sp,\phi,\psi,\mhi,\msi)$ is a weak solution to problem \eqref{eq:nsch1}--\eqref{eq:nschi} in the sense of Definition \ref{wsnschndmd}, it remains to prove that the initial value $\sp_0$ is achieved by $\sp$ and  $(\sp,\phi,\psi,\mhi,\msi)$ satisfies the energy inequality \eqref{energyineq}.

Like in \cite[Section 5.2]{ADG2013ndm}, we see that $\sp$ belongs to $BC_w \bigbracket{[0,T]; \ls }$ and the initial value $\sp_0$ can be attained in the sense that $\sp(0) = \sp_0$ in $\bm{L}^2(\Omega)$.
For the proof of the energy inequality \eqref{energyineq}, since we now have an additional term $R^N$ in $D^N$ compared to \cite[Section 5.3]{ADG2013ndm}, we provide the details here.
Since, for any $T>0$, we have
\begin{align*}
&\sigma_1\bracket{\phi^N} \mathop{\rightharpoonup}\limits^*  \sigma_1\bracket{\phi}\ \ \text{in}\ \lptbig{\infty}{0,T;\lp{\infty}},
\quad
\sigma_1\bracket{\phi^N} \to \sigma_1\bracket{\phi}\ \ \text{in}\ \lptbig{2}{0,T;H},\\
&  \sp^N \to \sp \ \ \text{in} \ L^2 \bigbracket{ 0,T;\bm{L}^2(\Omega) },\quad \mhi^N \rightharpoonup \mhi\ \ \text{in}\ L^2\bigbracket{ 0,T;\bm{H}^1(\Omega) },
\end{align*}
then it holds
$$ \intomega{\sigma_1\bracket{\phi^N} \abs{\sp^N}^2}{x}{\Omega} \to \intomega{\sigma_1\bracket{\phi} \abs{\sp}^2}{x}{\Omega}$$
and
$$\Biginner{\mhi^N}{\sigma_1\bracket{\phi^N}}_{\Omega} \to \Biginner{\mhi}{\sigma_1\bracket{\phi}}_{\Omega} \ \ \mathrm{in} \ L^1(0,T). $$
Combined with the strong convergence of $\overline{\phi^N}$ in $L^{\infty} (0,\infty)$, we can conclude
$$ \intinterval{R^N(\tau) g(\tau)}{\tau}{0}{T} \to \intinterval{R(\tau) g(\tau)}{\tau}{0}{T},
$$
with
$$ R \stackrel{\rm{def}}{=}   \bigbracket{\overline{\phi} -c} \intomega{ \Bigbracket{ \mhi - \frac{\gamma}{2} \abs{\sp}^2 } \sigma_1 \bracket{\phi} }{x}{\Omega},
$$
for any $T>0$ and $g \in W^{1,1}(0,\infty) \cap C_c ([0,\infty))$ with $g\ge 0$. Combining the weak convergence of $\sp^N$ and $\mhi^N, \ \msi^N$ in $L^2 \bigbracket{ 0,T;\bm{H}^1(\Omega) }$ and the almost everywhere point-wise convergence of $\phi^N$, we obtain the following weak convergence results
\begin{align*}
    \sqrt{\nu \bigbracket{\phi^N}} \nabla \sp^N &\rightharpoonup \sqrt{\nu \bigbracket{\phi}} \nabla \sp,
\end{align*}
and
\begin{align*}
    \sqrt{\mbhi^N \bigbracket{\phi^N}} \nabla \mhi^N &\rightharpoonup \sqrt{\mbhi \bigbracket{\phi}} \nabla \mhi,\quad
    \sqrt{\mbsi^N \bigbracket{\psi^N}} \nabla \msi^N \rightharpoonup \sqrt{\mbsi \bigbracket{\psi}} \nabla \msi,
\end{align*}
in $L^2 \bigbracket{ 0,T;\bm{L}^2(\Omega) }$. Then by the weak lower semi-continuity of norms and approximating the continuous functions with the step functions, we can conclude
$$ \liminf\limits_{N \to \infty} \intinterval{D^N(\tau) g(\tau)}{\tau}{0}{\infty} \ge \intinterval{D(\tau) g(\tau)}{\tau}{0}{\infty},$$
for all $g \in W^{1,1}(0,\infty) \cap C_c ([0,\infty))$ with $g\ge 0$, where
$$D(t)\stackrel{\rm{def}}{=} \intomega{\nu\bigbracket{\phi(t)} \bigabs{D \sp(t)}^2 + \mbhi\bigbracket{\phi(t)} \bigabs{\nabla \mhi(t)}^2 + \mbsi\bigbracket{\psi(t)} \bigabs{\nabla \msi(t)}^2 }{x}{\Omega} + R(t).
$$
Additionally, from the strong convergence $\sp^N(t) \to \sp(t)$ in $\bm{L}^2(\Omega)$, $\phi^N(t) \to \phi(t)$ and $\psi^N(t) \to \psi(t)$ in $V$ for almost all $t \in (0,\infty)$, it holds
\begin{equation*}
    \eto^N(t) \to \eto(\sp(t),\phi(t),\psi(t)),
\end{equation*}
for almost all $t \in (0,\infty)$. Thus, it holds
$$\eto(\sp_0,\phi_0,\psi_0)g(0) + \intinterval{\eto(\sp(\tau),\phi(\tau),\psi(\tau)) g'(\tau)}{\tau}{0}{\infty} \ge \intinterval{D(\tau) g(\tau)}{\tau}{0}{\infty},$$
for all $g \in W^{1,1}(0,\infty) \cap C_c ([0,\infty))$ with $g\ge 0$. As a consequence, the energy inequality \eqref{energyineq} follows from \cite[Lemma 2.2]{ADG2013ndm}.

With the energy inequality \eqref{energyineq}, we can derive further bounds of the global weak solution in $[0,\infty)$, under the additional assumption \eqref{additional_requirements_uniform_boundedness}. Since the mass relation \eqref{mass-rel} holds for all $t \ge 0$, then there exists $\delta\in (0,1)$ such that $\overline{\phi}(t)\in [-1+\delta,1-\delta]$ for $t\geq 0$.
From \eqref{eq:nsch4}, we can apply similar arguments as in the proof of Lemma \ref{apriorichemical} to get
\begin{equation*}
    |\overline{\mhi}| \le K_{\phi}(\delta)\bigbracket{\norm{\nabla \mhi} +1}.
\end{equation*}
By the Poincar\'e--Wirtinger inequality and $(\mathbf{H3})$, we obtain
    \begin{equation}\label{averagemhi}
        \norm{\mhi}^2 \le K_{\phi}(\delta) \Bigbracket{ \intomega{\mbhi(\phi) \abs{\nabla \mhi}^2 }{x}{\Omega} +1 }.
    \end{equation}
Estimates for $\mu_\psi$ can be obtained in a similar manner.
In addition, owing to Korn's inequality and Poincar\'e's inequality, there exists a constant $C_3>0$ depending on $\Omega$ and $\underline{\nu}$ such that
\begin{equation}\label{kornforsp}
    \norm{\sp}^2 \le C_3 \intomega{\nu(\phi) \abs{ D \sp }^2}{x}{\Omega}.
\end{equation}
Since it is assumed in \eqref{additional_requirements_uniform_boundedness} that $\norm{\sigma_1(\phi)}_{\lp{\infty}} \bigbracket{\overline{\phi}-c} \to 0 \ \mathrm{as} \ t \to \infty$, then there exists $t_0>0$ sufficiently large such that
    \begin{equation*}
        |\gamma |\norm{\sigma_1(\phi(t))}_{\lp{\infty}} \bigabs{\overline{\phi}(t)-c} \le \frac{1}{C_3},\quad \forall\, t\geq t_0.
    \end{equation*}
In view of the interpolation
    \begin{align*}
        \norm{\sigma_1\bracket{\phi}} &\le \norm{\sigma_1\bracket{\phi}}_{\lp{1}}^{\half} \norm{\sigma_1\bracket{\phi}}_{\lp{\infty}}^{\half}
        \le \abs{\Omega}^\half \big|\overline{\sigma_1\bracket{\phi}}\big|^\half  \norm{\sigma_1\bracket{\phi}}_{\lp{\infty}}^{\half},
    \end{align*}
  we infer from \eqref{mass-rel} and \eqref{averagemhi} that
    \begin{align*}
        & \Bigabs{ \bigbracket{\overline{\phi(t)}-c} \intomega{\sigma_1\bracket{\phi(t)} \mhi(t)}{x}{\Omega} } \\
        & \quad \le \bigabs{\overline{\phi(t)} -c} \abs{\Omega}^\half \big|\overline{\sigma_1\bracket{\phi}}\big|^\half  \norm{\sigma_1\bracket{\phi}}_{\lp{\infty}}^{\half} \norm{\mhi(t)} \\
        & \quad \le \epsilon \norm{\sigma_1\bracket{\phi}}_{\lp{\infty}} \abs{\overline{\phi(t)} -c} \norm{\mhi(t)}^2 + \frac{\abs{\Omega}}{4 \epsilon} \big|\overline{\sigma_1\bracket{\phi(t)}}\big| \bigabs{\overline{\phi_0} -c}e^{-\intinterval{\overline{\sigma_1\bracket{\phi(\tau)}}}{\tau}{0}{t}}  \\
        & \quad \le \epsilon K_{\phi}(\delta) \norm{\sigma_1\bracket{\phi}}_{\lp{\infty}} \abs{\overline{\phi(t)} -c}  + \epsilon \sigma_1^* K_{\phi}(\delta)  \intomega{\mbhi(\phi) \abs{\nabla \mhi}^2 }{x}{\Omega} - \frac{\abs{\Omega}}{4 \epsilon} \frac{\mathrm{d}}{\mathrm{d}t} \abs{\overline{\phi}-c},
    \end{align*}
for any $\epsilon>0$.
Taking $$\epsilon = \bigbracket{4 \bracket{1 + \sigma_1^*} K_{\phi}(\delta)}^{-1} $$
and substituting the above estimate as well as \eqref{kornforsp} into \eqref{energyineq}, we obtain
\begin{align*}
    & \eto (\sp(t),\phi(t),\psi(t)) + \half \intomega{\left(\nu (\phi) |D \sp|^2+ m_{\phi}(\phi)\abs{\gradmhi}^2 +m_{\psi}(\psi) \abs{\gradmsi}^2\right)}{(x,\tau)}{Q_{(t_0,t)}} \\
    & \quad \le \eto (\sp(t_0),\phi(t_0),\psi(t_0)) + \frac{\abs{\Omega}}{2 \epsilon} + \epsilon K_{\phi}(\delta) \intinterval{\norm{\sigma_1\bracket{\phi}}_{\lp{\infty}} \abs{\overline{\phi(t)} -c}}{\tau}{0}{\infty},
\end{align*}
for all $t\ge t_0$. This energy inequality combined with the assumption \eqref{additional_requirements_uniform_boundedness} and the estimates on the finite interval $[0,t_0]$ entails that
\begin{align*}
    & \sp \in L^{\infty} \bigbracket{0,\infty;\ls}\cap L^2 \bigbracket{0,\infty;\hssigma{1}},\\
    & \phi,\psi \in L^{\infty} \bigbracket{0,\infty;V},\\
    & \mhi,\msi \in L^2\ul \bigbracket{[0,\infty);V} \ \mathrm{with} \ \gradmhi,\gradmsi \in L^2 \bigbracket{0,\infty;\bm{L}^2(\Omega)}.
\end{align*}
Finally, by a standard argument as in \cite[Section 3]{GGW2018}, we can further obtain
\begin{align*}
    & \phi,\psi \in L^4\ul \bigbracket{[0,\infty);W} \cap L^2\ul \bigbracket{[0,\infty);\wkp{2}{p}}, \\
    & F_{\phi}'(\phi),F_{\psi}'(\psi) \in L^2\ul \bigbracket{[0,\infty);\lp{p}},
\end{align*}
for $p=6$ if $d=3$ and any $p \in [2,\infty)$ if $d=2$.

The proof of Theorem \ref{wsnschndme} is complete.
\epf


\section{Proof of Theorem \ref{wsnschdme}: the Case of Degenerate Mobilities }\label{wpdm}
\smallskip

In this section, we prove Theorem \ref{wsnschdme} on the existence of a global weak solution to problem \eqref{eq:nsch1}--\eqref{eq:nschi} in the case of degenerate mobilities. There are several approaches to managing the degeneracy of mobility functions \cite{ADG2013dm,CFG2024,EG1996}. Here we shall adapt the strategy in \cite{ADG2013dm}, taking advantage of the results obtained in Theorem \ref{wsnschndme}.

In view of Remark \ref{rem:mo-pot}, we assume below that the assumptions $(\mathbf{H0})$, $(\mathbf{H1})$, $(\mathbf{H2})$ and $(\mathbf{H3}*)$ are satisfied. Then the proof for Theorem \ref{wsnschdme} consists of two steps. First, we regularize the degenerate mobilities by non-degenerate functions and construct approximate solutions using Theorem \ref{wsnschndme}. Second, we derive uniform estimates and pass to the limit to obtain a global weak solution to the original problem.

\subsection{The approximate problem.} \label{approximate_degenerate}
We regularize the degenerate mobility functions $\mbhi$, $\mbsi$ in a different way compared to \cite{ADG2013dm}. More precisely, we introduce the following strictly positive functions $\mbhi^{\epsilon}$,  $\mbsi^{\epsilon}$ such that
\begin{equation}\notag
    \mbhi^{\epsilon}(s)
    = \left\{
        \begin{aligned}
            & 2 \epsilon, &&\ s \le -1+\delta_{\phi,1}^{\epsilon}, \\
            & \epsilon + \mbhi(s), &&\ s \in \bigbracket{ -1+\delta_{\phi,1}^{\epsilon},1-\delta_{\phi,2}^{\epsilon} },\\
            & 2 \epsilon, &&\ s \ge 1-\delta_{\phi,2}^{\epsilon},
        \end{aligned}
        \right.
\end{equation}
\begin{equation}\notag
    \mbsi^{\epsilon}(s)
    = \left\{
        \begin{aligned}
            & 2 \epsilon, &&\ s \le \delta_{\psi,1}^{\epsilon}, \\
            & \epsilon + \mbsi(s), &&\ s \in \bigbracket{\delta_{\psi,1}^{\epsilon}, 1-\delta_{\psi,2}^{\epsilon}},\\
            & 2 \epsilon, &&\ s \ge 1-\delta_{\psi,2}^{\epsilon},
        \end{aligned}
        \right.
\end{equation}
where
$$
\epsilon \in I_M\stackrel{\rm{def}}{=} \left(0,\,\half \mbhi(0)\right] \ \bigcap\ \left(0,\,\half \mbsi\left(\half\right)\right],
$$
and $\delta_{\phi,k}^{\epsilon}$, $\delta_{\psi,k}^{\epsilon}$, $k \in \set{1,2}$ are given by
\begin{align*}
    & \delta_{\phi,1}^{\epsilon} \stackrel{\rm{def}}{=} \inf \set{s \in(0,1):\mbhi(-1+s) \ge \epsilon}, \quad  \delta_{\phi,2}^{\epsilon} \stackrel{\rm{def}}{=} \inf \set{s \in(0,1):\mbhi(1-s) \ge \epsilon}, \\
    & \delta_{\psi,1}^{\epsilon} \stackrel{\rm{def}}{=} \inf \left\{s \in\left(0,\half\right):\mbsi(s) \ge \epsilon\right\}, \quad \ \,  \delta_{\psi,2}^{\epsilon} \stackrel{\rm{def}}{=} \inf \left\{s \in\left(0,\half\right):\mbsi(1-s) \ge \epsilon\right\}.
\end{align*}
Obviously, $\delta_{\phi,k}^{\epsilon}$, $\delta_{\psi,k}^{\epsilon}$, $k \in \set{1,2}$, are strictly increasing with respect to $\epsilon\in I_M$, and tend to $0$ as $\epsilon \to 0$.

Corresponding to $\mbhi^{\epsilon}$ and $\mbsi^{\epsilon}$, we define two functions $\whi^{\epsilon} \in C^2(\r;\r^+_0)$ and $ \wsi^{\epsilon} \in C^2(\r;\r^+_0)$ as follows
\begin{align*}
    & (\whi^{\epsilon})''(s) = \frac{1}{\mbhi^{\epsilon}(s)}, \ \ \whi^{\epsilon}(0)=(\whi^{\epsilon})'(0) = 0, \ \ s\in \r,\\
    & (\wsi^{\epsilon})''(s) = \frac{1}{\mbsi^{\epsilon}(s)}, \ \ \wsi^{\epsilon} \Bigbracket{\half}=(\wsi^{\epsilon})' \Bigbracket{\half} = 0, \ \ s\in \r.
\end{align*}
Since
$$ \mbhi^{\epsilon}(s_1) \ge \mbhi(s_1) + \epsilon, \quad \mbsi^{\epsilon}(s_2) \ge \mbsi(s_2) + \epsilon, $$
for all $s_1\in[-1,1]$ and $s_2 \in [0,1]$, respectively, we find
\begin{align}
    &\whi^{\epsilon}(s) \le \whi(s), \ \ \forall\, s \in (-1,1), \qquad
    \wsi^{\epsilon}(s) \le \wsi(s), \ \ \forall\, s \in (0,1).
    \label{upper-W}
\end{align}

Next, we introduce the corresponding regular parts of the approximate potential functions:
\begin{align*}
    & \bigbracket{ \tfhi^{\epsilon} }'' = \frac{\mbhi}{\mbhi^{\epsilon}} \fhi'' \in C\bigbracket{[-1,1]}, \quad
    \bigbracket{ \tfsi^{\epsilon} }'' = \frac{\mbsi}{\mbsi^{\epsilon}} \fsi'' \in C\bigbracket{[0,1]},
\end{align*}
which satisfy
\begin{equation*}
    \tfhi^{\epsilon}(0) = \bigbracket{\tfhi^{\epsilon}}'(0)=\tfsi^{\epsilon}\Bigbracket{\half} = \bigbracket{\tfsi^{\epsilon}}'\Bigbracket{\half}=0.
\end{equation*}
Then it is straightforward to check that
\begin{equation} \label{upper-F}
    \bigbracket{ \tfhi^{\epsilon} }'' \le \fhi'', \quad  \bigbracket{ \tfsi^{\epsilon} }'' \le \fsi'',
\end{equation}
and
\begin{align}
    & 0 < \bigbracket{\tfhi^{\epsilon}}''\mbhi^{\epsilon} = \fhi'' \mbhi \le \alpha, \quad  0 < \bigbracket{\tfsi^{\epsilon}}''\mbsi^{\epsilon} = \fsi'' \mbsi \le \alpha.
    \label{ass:sig1W-e}
\end{align}
In addition, $\tfhi^{\epsilon}$ (resp. $\tfsi^{\epsilon}$) converges uniformly to $\fhi$ in $C\bigbracket{[-1,1]}$ (resp. to $\fsi$ in $C\bigbracket{[0,1]}$). To see this, we first note that they are monotonously increasing functions with respect to $\epsilon^{-1}$ as $\epsilon$ tends to $0$, and $\tfhi^{\epsilon}$ (resp. $\tfsi^{\epsilon
}$) converges to $\fhi$ (resp. to $\fsi$) in $(-1,1)$ (resp. in $(0,1)$) pointwisely.
So the limit $\limit{\epsilon \to 0} \tfhi^{\epsilon}(1)$ exists and is less than or equal to $\fhi(1)$.
On the other hand, it follows from the convexity of $\tfhi^{\epsilon}$ that
$$\tfhi^{\epsilon}(1) \ge 2 \tfhi^{\epsilon}\Bigbracket{\halfof{1+s}} - \tfhi^{\epsilon}(s),
\quad \forall\, s \in [-1,1].
$$
Passing to the limit in the above inequality for any fixed $s \in (-1,1)$, we obtain
$$ \limit{\epsilon \to 0} \tfhi^{\epsilon}(1) \ge 2 \fhi \Bigbracket{\halfof{1+s}} - \fhi(s). $$
Then by the continuity of $\fhi$, we can conclude
$$ \limit{\epsilon \to 0} \tfhi^{\epsilon}(1) \ge 2 \fhi(1) - \fhi(1)=\fhi(1), $$
which gives
$$ \limit{\epsilon \to 0} \tfhi^{\epsilon}(1) = \fhi(1). $$
In a similar manner, we have
\begin{align*}
    & \limit{\epsilon \to 0} \tfhi^{\epsilon}(-1) = \fhi(-1), \quad\ \
     \limit{\epsilon \to 0} \tfsi^{\epsilon}(0) = \fsi(0),\quad\ \
     \limit{\epsilon \to 0} \tfsi^{\epsilon}(1) = \fsi(1).
\end{align*}
Finally, it follows from Dini's theorem that the above uniform convergence of $\tfhi^{\epsilon}$ (resp. $\tfsi^{\epsilon}$) on $[-1,1]$ (resp. on $[0,1]$) holds.

Let us define approximations of $\fhi$ and $\fsi$ as follows:
$$
\fhi^{\epsilon} = \tfhi^{\epsilon} + c_{\epsilon} \flnhi, \quad \fsi^{\epsilon} = \tfsi^{\epsilon} + c_{\epsilon} \flnsi,
$$
where
\begin{align*}
    & \flnhi(s) = (1+s)\ln(1+s) + (1-s)\ln(1-s), \quad  s \in (-1,1),\\
    & \flnsi(s) = s \ln(s) + (1-s)\ln(1-s), \quad  s \in (0,1),
\end{align*}
and
\begin{align*}
    & c_{\epsilon} = \min \left\{1,\  -\frac{\delta_{\phi,1}^{\epsilon}}{\flnhi'\bracket{-1+\delta_{\phi,1}^{\epsilon}}},\  \frac{\delta_{\phi,2}^{\epsilon}}{\flnhi'\bracket{1-\delta_{\phi,2}^{\epsilon}}},\   -\frac{\delta_{\psi,1}^{\epsilon}}{\flnsi'\bracket{\delta_{\psi,1}^{\epsilon}}},\  \frac{\delta_{\psi,2}^{\epsilon}}{\flnsi'\bracket{1-\delta_{\psi,2}^{\epsilon}}}\right\} \in (0,1].
\end{align*}
In particular, we observe that $c_\epsilon\to 0$ as $\epsilon\to 0$. In addition, it follows from Remark \ref{entropy_logarithmic_singularity&mobility_monotonicity}, \eqref{upper-W} and \eqref{ass:sig1W-e} that
\begin{align}
    & \abs{(\fhi^{\epsilon})'(s_1)} \le (\alpha + 2 C_{\phi}) \abs{\whi'(s_1)}, \quad\ \   \abs{(\fsi^{\epsilon})'(s_2)} \le (\alpha +  C_{\psi}) \abs{\wsi'(s_2)},
    \label{estimate_appontential_entropy}
\end{align}
for all $s_1 \in (-1,1)$ and $s_2 \in (0,1)$, respectively.

Consider now the following problem with regularized mobilities and approximated potentials:
\begin{numcases}{}
    \pdif{(\rho (\phi) \sp)}{t} + \div\Bigbracket{\sp \otimes \bigbracket{\rho (\phi) \sp - \gamma \mbhi^{\epsilon}(\phi) \nabla \mhi} } - \div \bigbracket{\nu (\phi) D\sp} + \nabla\pi \notag\\
    \quad = \mu_{\phi} \nabla \phi + \mu_{\psi} \nabla \psi \label{eq:nsapdmch1}, & in $Q$,\\
    \div \sp = 0 \label{eq:nsapdmch2}, & in $Q$,\\
    \pdif{\phi}{t} + \sp \cdot \nabla \phi + \sigma_1 \bracket{\phi} \bigbracket{\overline{\phi}-c} = \div \bigbracket{m_{\phi}^{\epsilon} (\phi) \nabla \mhi} \label{eq:nsapdmch3}, & in $Q$,\\
    \mu_{\phi} = -\Delta \phi + \sigma_2 \mathcal{N} \bigbracket{ \phi - \overline{\phi} } + \bigbracket{F_{\phi}^{\epsilon}}'(\phi) + \pdif{G}{\phi}(\phi,\psi) \label{eq:nsapdmch4}, & in $Q$,\\
    \pdif{\psi}{t} + \sp \cdot \nabla \psi = \div \bigbracket{m_{\psi}^{\epsilon} (\psi) \nabla \msi} \label{eq:nsapdmnch5}, & in $Q$,\\
    \mu_{\psi} = -\beta \Delta \psi + \bigbracket{F_{\psi}^{\epsilon}}'(\psi) + \pdif{G}{\psi}(\phi,\psi) \label{eq:nsapdmch6}, & in $Q$,\\
    \sp = \mathbf{0}, \ \pdif{\phi}{\nv} = \pdif{\mu_{\phi}}{\nv} = \pdif{\psi}{\nv} = \pdif{\mu_{\psi}}{\nv} = 0 \label{eq:nsapdmchb}, & on $S$,\\
    \sp|_{t=0} = \sp_0(x), \ \ \phi|_{t=0} = \phi_0(x), \ \  \psi|_{t=0} = \psi_0(x) \label{eq:nsapdmchi}, & in $\Omega$.
\end{numcases}
We are able to apply Theorem \ref{wsnschndme} to conclude that problem \eqref{eq:nsapdmch1}--\eqref{eq:nsapdmchi} admits a global weak solution $ \bigbracket{ \sp^{\epsilon},\phi^{\epsilon},\psi^{\epsilon},\mhi^{\epsilon},\msi^{\epsilon} }$ on $[0,\infty)$ in the sense of Definition \ref{wsnschndmd}.
Define
$$  \rho^{\epsilon}(\phi^\epsilon) = \gamma \phi^{\epsilon} + \halfof{\rho_1+\rho_2},\quad\ \
\jhi^{\epsilon} = \mbhi^{\epsilon}(\phi^{\epsilon}) \nabla \mhi^{\epsilon}, \quad \ \
\jsi^{\epsilon}= \mbsi^{\epsilon}(\psi^{\epsilon}) \nabla \msi^{\epsilon}.
$$
Using the identity
\begin{align*}
    & \left(\bigbracket{ \bigbracket{F_{\phi}^{\epsilon}}'(\phi^{\epsilon}) + \pdif{G}{\phi}(\phi^{\epsilon},\psi^{\epsilon})} \nabla \phi^{\epsilon}, \bmtheta\right)_Q
    + \left(\bigbracket{ \bigbracket{F_{\psi}^{\epsilon}}'(\psi^{\epsilon}) + \pdif{G}{\psi}(\phi^{\epsilon},\psi^{\epsilon})} \nabla \psi^{\epsilon},\bmtheta\right)_Q \\
    & \quad = \left(\nabla \bigbracket{F_{\phi}^{\epsilon}(\phi^{\epsilon}) + F_{\psi}^{\epsilon}(\psi^{\epsilon}) + G(\phi^{\epsilon},\psi^{\epsilon})},\bmtheta\right)_Q=0,
    \quad \forall\, \bmtheta\in   C^{\infty}_0(Q)^d\ \text{with}\ \div \bmtheta =0,
\end{align*}
we can write the weak formulation for $\sp^\epsilon$ as
\begin{align}
    & -\biginner{\rho^{\epsilon} \sp ^{\epsilon}}{\pdif{\bmtheta}{t}}_{Q} - \biginner{ \rho^{\epsilon} \sp^{\epsilon} \otimes \sp^{\epsilon}}{\nabla \bmtheta}_{Q} + \biginner{\nu(\phi^{\epsilon}) D \sp^{\epsilon}}{\nabla  \bmtheta}_{Q} + \biginner{\sp^{\epsilon} \otimes \gamma \jhi^{\epsilon}}{\nabla \bmtheta}\nonumber\\
    & \quad = \intomega{\bigbracket{ \sigma_2 \mathcal{N} \bigbracket{ \phi^{\epsilon} - \overline{\phi^{\epsilon}} }\nabla \phi^{\epsilon} -\Delta \phi^{\epsilon} \nabla \phi^{\epsilon} - \beta \Delta \psi^{\epsilon} \nabla \psi^{\epsilon}}\cdot \bmtheta}{(x,t)}{Q},\label{testapdmns}
\end{align}
for all $\bmtheta \in C^{\infty}_0(Q)^d$ with $\div \bmtheta =0$. Moreover,
\begin{align}
    &  -\biginner{\phi^{\epsilon}}{\pdif{\zeta}{t}}_{Q} + \intomega{ \bigbracket{ \sp^{\epsilon} \cdot \nabla \phi^{\epsilon} + \sigma_1(\phi^{\epsilon})\bigbracket{\overline{\phi^{\epsilon}}-c} } \zeta}{(x,t)}{Q} = -\biginner{\jhi^{\epsilon} }{ \nabla \zeta}_{Q}, \label{testapdmch1}
    \\
    & -\biginner{\psi^{\epsilon}}{\pdif{\zeta}{t}}_{Q} + \intomega{\sp^{\epsilon} \cdot \nabla \psi^{\epsilon} \zeta}{(x,t)}{Q} = -\biginner{\jsi^{\epsilon} }{ \nabla \zeta}_{Q}, \label{testapdmch2}
\end{align}
for all $\zeta \in C^{\infty}_0 \big((0,\infty);C^1 \big(\overline{\Omega} \big) \big)$, and
\begin{align}
    &\intomega{\Bigbracket{(\tfhi^{\epsilon})''(\phi^\epsilon)\nabla \phi^{\epsilon} + \frac{\p^2 G}{\p \phi^2}(\phi^{\epsilon},\psi^{\epsilon}) \nabla \phi^{\epsilon} + \frac{\p^2 G}{\p \phi \p \psi}(\phi^{\epsilon},\psi^{\epsilon}) \nabla \psi^{\epsilon}}\cdot \bigbracket{\mbhi^{\epsilon}(\phi^{\epsilon})\bmeta}}{(x,t)}{Q} \nonumber\\
    & \qquad + \biginner{\Delta \phi^{\epsilon} - c_\epsilon \flnhi'(\phi^{\epsilon}) }{\div \bigbracket{\mbhi^{\epsilon} (\phi^{\epsilon}) \bmeta}}_{Q} + \sigma_2 \biginner{ \nabla \mathcal{N} \bigbracket{ \phi^{\epsilon} - \overline{\phi^{\epsilon}} }}{\mbhi^{\epsilon}(\phi^{\epsilon})\bmeta}_{Q} \nonumber \\
    & \quad = \inner{\jhi^{\epsilon}}{\bmeta}_{Q}, \label{testapdmmhi}\\
    &\intomega{\Bigbracket{(\tfsi^{\epsilon})''(\psi^\epsilon)\nabla \psi^{\epsilon} + \frac{\p^2 G}{\p \psi^2}(\psi^{\epsilon},\psi^{\epsilon}) \nabla \psi^{\epsilon} + \frac{\p^2 G}{\p \phi \p \psi}(\phi^{\epsilon},\psi^{\epsilon}) \nabla \phi^{\epsilon}}\cdot \bigbracket{\mbsi^{\epsilon}(\psi^{\epsilon})\bmeta }}{(x,t)}{Q} \nonumber\\
    & \qquad + \biginner{\beta \Delta \psi^{\epsilon} - c_\epsilon \flnsi'(\psi^{\epsilon})}{\div \bigbracket{\mbsi^{\epsilon}(\psi^{\epsilon}) \bmeta}}_{Q} \nonumber \\
    & \quad = \inner{\jsi^{\epsilon}}{\bmeta}_{Q},\label{testapdmmsi}
\end{align}
for all $\bmeta \in C_c \bigbracket{[0,\infty);\bm{H}^1\bracket{\Omega}} \cap \bm{L}^{\infty}(Q)$ with $\bmeta \cdot \nv = 0$ on $S$.
\medskip

Next, we derive estimates for the approximate solution $ \bigbracket{\sp^{\epsilon},\phi^{\epsilon},\psi^{\epsilon},\mhi^{\epsilon},\msi^{\epsilon}}$ that are uniform with respect to $\epsilon$.
\medskip

\bl \label{energyestimateapdm}
Let $\bigbracket{\sp^{\epsilon},\phi^{\epsilon},\psi^{\epsilon},\mhi^{\epsilon},\msi^{\epsilon}}$ be a global weak solution to problem \eqref{eq:nsapdmch1}--\eqref{eq:nsapdmchi} in the sense of Definition \ref{wsnschndmd}. For any $T > 0$, we have the following estimates:
\medskip
\begin{enumerate}[(i)]
    \item Set $\widehat{\bmj}_{\phi}^{\epsilon} = \sqrt{\mbhi^{\epsilon}(\phi^{\epsilon})} \nabla \mhi^{\epsilon}, \ \widehat{\bmj}_{\phi}^{\epsilon} = \sqrt{\mbsi^{\epsilon}(\psi^{\epsilon})} \nabla \msi^{\epsilon}$. Then
        \begin{align}
           &\sup\limits_{0\le t \le T} \Big[ \intomega{\rho^{\epsilon}(t) \frac{ \abs{\sp^{\epsilon} (t) }^2}{2} + \half \abs{\nabla \phi^{\epsilon}(t)}^2 + \halfof{\beta} \abs{\nabla \psi^{\epsilon}(t)}^2 }{x}{\Omega} \nonumber\\
           & \qquad + \halfof{\sigma_2} \bignorm{ \phi^{\epsilon} (t) - \overline{ \phi^{\epsilon} } (t) }^2_* + \intomega{\fhi^{\epsilon}\bracket{\phi^{\epsilon}(t)} + \fsi^{\epsilon}\bracket{\psi^{\epsilon}(t)} + G \bigbracket{\phi^{\epsilon}(t),\psi^{\epsilon}(t)} }{x}{\Omega} \nonumber \\
           & \qquad + \intomega{\whi^{\epsilon}(\phi^{\epsilon})}{x}{\Omega} + \intomega{\wsi^{\epsilon}(\psi^{\epsilon})}{x}{\Omega} \Big]\nonumber\\
           & \qquad + \intomega{\nu (\phi^{\epsilon}) \abs{D \sp^{\epsilon}}^2 }{(x,\tau)}{Q_T} + \Bignorm{\widehat{\bmj}_{\phi}^{\epsilon}}^2_{L^2(Q_T)} + \Bignorm{\widehat{\bmj}_{\psi}^{\epsilon}}^2_{L^2(Q_T)}\nonumber\\
           & \qquad + \norm{\Delta \phi^{\epsilon}}^2_{L^2(Q_T)} + \intomega{(\fhi^{\epsilon})''(\phi^{\epsilon}) \abs{\nabla \phi^{\epsilon}}^2}{(x,\tau)}{Q_T} \nonumber \\
           & \qquad + \beta \norm{\Delta \psi^{\epsilon}}^2_{L^2(Q_T)} + \intomega{(\fsi^{\epsilon})''(\psi^{\epsilon}) \abs{\nabla \psi^{\epsilon}}^2}{(x,\tau)}{Q_T} \nonumber \\
           & \quad \le C .\label{unif-i}
        \end{align}
\item $\epsilon c_{\epsilon}^2 \norm{\flnhi'(\phi^{\epsilon})}_{L^2(Q_T)}^2\le C$ and $\epsilon c_{\epsilon}^2 \norm{\flnsi'(\psi^{\epsilon})}_{L^2(Q_T)}^2 \le C$.
\end{enumerate}
\medskip
Here, the positive constant $C$ is independent of $\epsilon$.
\el
\medskip

\bpf (i)
    First, we derive estimates from the energy inequality:
    \begin{align*}
        & \eto^{\epsilon} \bigbracket{\sp^{\epsilon}(t),\phi^{\epsilon}(t),\psi^{\epsilon}(t)}
        \\
        &\qquad + \intomega{\left(\nu (\phi^{\epsilon}) |D \sp^{\epsilon}|^2+ m_{\phi}^{\epsilon}(\phi^{\epsilon})\abs{\gradmhi^{\epsilon}}^2 +m_{\psi}^{\epsilon}(\psi^{\epsilon}) \abs{\gradmsi^{\epsilon}}^2\right)}{(x,\tau)}{Q_{(s,t)}}\nonumber\\
        & \qquad + \intomega{ \sigma_1 (\phi^{\epsilon}) \bigbracket{ \overline{\phi^{\epsilon}} - c } \Bigbracket{ \mhi^{\epsilon} - \halfof{\gamma} \abs{\sp^{\epsilon}}^2 } }{\bracket{x,\tau}}{Q_{\bracket{s,t}}}\\
        &\quad \le \eto^{\epsilon} \bigbracket{\sp^{\epsilon}(s),\phi^{\epsilon}(s),\psi^{\epsilon}(s)} ,
    \end{align*}
    where
    \begin{align*}
        \eto^{\epsilon}(\sp,\phi,\psi)
        &\stackrel{\rm{def}}{=} \intomega{ \left(\half \abs{\nabla \phi}^2 + \halfof{\beta} \abs{\nabla \psi}^2 + \halfof{\sigma_2} \bigabs{\nabla \mathcal{N} \bigbracket{\phi-\overline{\phi}}}^2\right) }{x}{\Omega}\\
        &\qquad + \intomega{ \left(\fhi^{\epsilon}(\phi) + \fsi^{\epsilon}(\psi) + G(\phi,\psi)\right) }{x}{\Omega}
         + \intomega{ \halfof{\rho(\phi)\abs{\sp}^2} }{x}{\Omega}.
    \end{align*}
  Special attention should be paid to the term $\sigma_1(\phi^{\epsilon})\bigbracket{\overline{\phi^{\epsilon}}-c}\mhi^{\epsilon}$. It follows from equation \eqref{eq:nsapdmch4}, the facts that $\phi^\epsilon\in (-1,1)$, $\psi^\epsilon\in (0,1)$ almost everywhere in $Q$, $(\mathbf{H3}*)$, H\"{o}lder's inequality and Young's inequality that
    \begin{align*}
         \Bigabs{\intomega{\sigma_1(\phi^{\epsilon})\mhi^{\epsilon}}{x}{\Omega}}
        & = \bigabs{ \inner{\sigma_1 (\phi^{\epsilon})}{ -\Delta \phi^{\epsilon} + \sigma_2 \mathcal{N} \bigbracket{\phi^{\epsilon} - \overline{\phi^{\epsilon}}} +  (\fhi^{\epsilon})'\bracket{\phi^{\epsilon}} + \pdif{G\bracket{\phi^{\epsilon},\psi^{\epsilon}}}{\phi} } }\\
        &  \le \frac{1}{4} \norm{\Delta \phi^{\epsilon}}^2 +C.
    \end{align*}
    In particular, the uniform (pointwise) boundedness of $\bigabs{ \sigma_1(\phi^{\epsilon}) (\fhi^{\epsilon})'(\phi^{\epsilon}) }$ due to \eqref{ass:sig1W} and \eqref{estimate_appontential_entropy} has been used. As a consequence, it holds
    \begin{align}
        & \eto^{\epsilon} \bigbracket{\sp^{\epsilon}(t),\phi^{\epsilon}(t),\psi^{\epsilon}(t)} + \intomega{\nu (\phi^{\epsilon}) |D \sp^{\epsilon}|^2+ m_{\phi}^{\epsilon}(\phi^{\epsilon})\abs{\gradmhi^{\epsilon}}^2 +m_{\psi}^{\epsilon}(\psi^{\epsilon}) \abs{\gradmsi^{\epsilon}}^2}{(x,\tau)}{Q_{(s,t)}}\nonumber\\
        & \quad \le \eto^{\epsilon} \bigbracket{\sp^{\epsilon}(s),\phi^{\epsilon}(s),\psi^{\epsilon}(s)}  + \frac{1}{4} \norm{\Delta \phi^{\epsilon}}^2_{L^2(Q_{(s,t)})} + C \norm{\sp^{\epsilon}}^2_{L^2(Q_{(s,t)})}
        + C(t-s).
        \label{eq:energyestimateapdm}
    \end{align}

Next, we derive entropy estimates for the approximate solution (cf. \cite{EG1997}). To this end, we choose $\zeta = (\whi^{\epsilon})'(\phi^{\epsilon})$ as the test function in \eqref{eq:nsapdmch3}. Then substituting equation \eqref{eq:nsapdmch4}, which holds almost everywhere in $Q$, into the resulting equality, we get
    \begin{align*}
        &\intomega{\pdif{\phi^{\epsilon}}{t} (\whi^{\epsilon})'(\phi^{\epsilon}) }{(x,\tau)}{Q_{(s,t)}} + \intomega{ \bigbracket{ \sp^{\epsilon} \cdot \nabla \phi^{\epsilon} + \sigma_1\bracket{\phi^{\epsilon}} \bigbracket{\overline{\phi^{\epsilon}}-c} } (\whi^{\epsilon})'(\phi^{\epsilon}) }{(x,\tau)}{Q_{(s,t)}} \\
        & \quad = \intomega{\bigbracket{-\Delta \phi^{\epsilon} + \sigma_2 \mathcal{N} \bigbracket{\phi^{\epsilon} - \overline{\phi^{\epsilon}}} + (\fhi^{\epsilon})'(\phi^{\epsilon}) + \pdif{G(\phi^{\epsilon},\psi^{\epsilon})}{\phi} } \Delta \phi^{\epsilon}}{(x,\tau)}{Q_{(s,t)}},
    \end{align*}
    where the fact $\div \bigbracket{\mbhi^{\epsilon} (\phi^{\epsilon}} \nabla \zeta) = \Delta \phi^{\epsilon}$ has been used. We note that
    \begin{align*}
        &\intomega{\pdif{\phi^{\epsilon}}{t} (\whi^{\epsilon})'(\phi^{\epsilon}) }{(x,\tau)}{Q_{(s,t)}} = \intomega{\whi^{\epsilon} \bigbracket{\phi^{\epsilon}(t)}}{x}{\Omega} - \intomega{\whi^{\epsilon}(\phi^{\epsilon}(s))}{x}{\Omega},\\
        & \intomega{\sp^{\epsilon} \cdot \nabla \phi^{\epsilon} (\whi^{\epsilon})'(\phi^{\epsilon}) }{(x,\tau)}{Q_{(s,t)}} = \intomega{\sp^{\epsilon} \cdot \nabla \bigbracket{\whi^{\epsilon}(\phi^{\epsilon})} }{(x,\tau)}{Q_{(s,t)}}=0.
    \end{align*}
    Besides, it holds
    \begin{align*}
    \intomega{\bigbracket{ (\fhi^{\epsilon})'(\phi^{\epsilon}) } \Delta \phi^{\epsilon}}{(x,\tau)}{Q_{(s,t)}}
    = -  \intomega{\bigbracket{ (\fhi^{\epsilon})''(\phi^{\epsilon}) } |\nabla \phi^{\epsilon}|^2}{(x,\tau)}{Q_{(s,t)}},
    \end{align*}
    which follows from the same argument as in \cite[Section 3]{ADG2013dm} by approximating $\phi^\epsilon$ with $\phi^\epsilon_\eta=\eta\phi^\epsilon$ for $\eta\in (0,1)$ and then letting $\eta\to 1$ with Lebesgue's dominated convergence theorem. Hence, from the above identities and
$$ \sigma_1(\phi^{\epsilon}) \bigabs{(\whi^{\epsilon})'(\phi^{\epsilon})} \le \sigma_1(\phi^{\epsilon}) \bigabs{\whi'(\phi^{\epsilon})} \le C,
$$
we can deduce that
    \begin{align*}
        &\intomega{\whi^{\epsilon} \bigbracket{\phi^{\epsilon}(t)}}{x}{\Omega} + \norm{\Delta \phi^{\epsilon}}^2_{L^2(Q_{(s,t)})} + \sigma_2 \bignorm{\phi^{\epsilon}-\overline{\phi^{\epsilon}}}^2_{L^2(Q_{(s,t)})} \\
        & \qquad + \intomega{(\fhi^{\epsilon})''(\phi^{\epsilon}) \abs{\nabla \phi^{\epsilon}}^2}{(x,\tau)}{Q_{(s,t)}}\\
        & \quad \le \intomega{\whi^{\epsilon}(\phi^{\epsilon}(s))}{x}{\Omega} + \intomega{ \pdif{G}{\phi} \bigbracket{\phi^{\epsilon},\psi^{\epsilon}} \Delta \phi^{\epsilon}}{(x,\tau)}{Q_{(s,t)}} + C(t-s)\\
        & \quad \le \intomega{\whi^{\epsilon}(\phi^{\epsilon}(s))}{x}{\Omega} + \frac14 \norm{\Delta \phi^{\epsilon}}^2_{L^2(Q_{(s,t)})} + C (t-s),
    \end{align*}
where $C>0$ is independent of $s,\ t$ and $\epsilon$.
Analogously, similar estimates hold for $\psi^\epsilon$.
Adding the entropy estimates of $\phi$, $\psi$ and the energy estimate \eqref{eq:energyestimateapdm} together, we obtain
    \begin{align}
        &\eto^{\epsilon} \bigbracket{\sp^{\epsilon}(t),\phi^{\epsilon}(t),\psi^{\epsilon}(t)} + \intomega{\whi^{\epsilon}(\phi^{\epsilon}(t)) + \wsi^{\epsilon}(\psi^{\epsilon}(t))}{x}{\Omega}
        \notag \\
        & \qquad + \intomega{\left(\nu (\phi^{\epsilon}) |D \sp^{\epsilon}|^2+ m_{\phi}^{\epsilon}(\phi^{\epsilon})\abs{\gradmhi^{\epsilon}}^2 +m_{\psi}^{\epsilon}(\psi^{\epsilon}) \abs{\gradmsi^{\epsilon}}^2\right)}{(x,\tau)}{Q_{(s,t)}}
        \notag\\
        & \qquad + \sigma_2 \bignorm{\phi^{\epsilon}-\overline{\phi^{\epsilon}}}^2_{L^2(Q_{(s,t)})} + \frac{1}{2} \norm{\Delta \phi^{\epsilon}}^2_{L^2(Q_{(s,t)})} + \frac{\beta}{2} \norm{\Delta \psi^{\epsilon}}^2_{L^2(Q_{(s,t)})}
        \notag \\
        & \qquad + \intomega{(\fhi^{\epsilon})''(\phi^{\epsilon}) \abs{\nabla \phi^{\epsilon}}^2}{(x,\tau)}{Q_{(s,t)}} + \intomega{(\fsi^{\epsilon})''(\psi^{\epsilon}) \abs{\nabla \psi^{\epsilon}}^2}{(x,\tau)}{Q_{(s,t)}}
        \notag\\
        & \quad \le \eto^{\epsilon} \bigbracket{\sp^{\epsilon}(s),\phi^{\epsilon}(s),\psi^{\epsilon}(s)} + \intomega{\whi^{\epsilon}(\phi^{\epsilon}(s)) + \wsi^{\epsilon}(\psi^{\epsilon}(s))}{x}{\Omega}
        \notag\\
        & \qquad + C \norm{\sp^{\epsilon}}^2_{L^2(Q_{(s,t)})} + C(t-s)
        \notag\\
        & \quad \le \eto^{\epsilon} \bigbracket{\sp^{\epsilon}(s),\phi^{\epsilon}(s),\psi^{\epsilon}(s)} + \intomega{\whi^{\epsilon}(\phi^{\epsilon}(s)) + \wsi^{\epsilon}(\psi^{\epsilon}(s))}{x}{\Omega}
        \notag\\
        & \qquad + C \intinterval{ \eto^{\epsilon} \bigbracket{\sp^{\epsilon}(\tau),\phi^{\epsilon}(\tau),\psi^{\epsilon}(\tau)} }{\tau}{s}{t} + C(t-s),
        \label{eq:energyestimateapdm-sum}
    \end{align}
    for almost all $s \in [0,\infty)$ (including $s=0$) and all $t \in [s,\infty)$.

    Since $c_\epsilon\in (0,1]$, it follows from \eqref{upper-W} and \eqref{upper-F} that the initial energy $\eto^{\epsilon} \bigbracket{\sp^{\epsilon}(0),\phi^{\epsilon}(0),\psi^{\epsilon}(0)}=\eto^{\epsilon} \bigbracket{\sp_0,\phi_0,\psi_0}$ is uniformly bounded with respect to $\epsilon$. Besides, we see that
    $$\intomega{\whi^{\epsilon}(\phi^{\epsilon}(0)) + \wsi^{\epsilon}(\psi^{\epsilon}(0))}{x}{\Omega}\leq \intomega{\whi(\phi_0) + \wsi(\psi_0)}{x}{\Omega}$$ is also uniformly bounded.
    Then an application of Gronwall's inequality to \eqref{eq:energyestimateapdm-sum} with $s=0$, $t=T>0$ yields the uniform estimate \eqref{unif-i}.

\medskip

(ii) Testing \eqref{eq:nsapdmch4} by $\phi^{\epsilon} - \overline{\phi^{\epsilon}}$, we get
\begin{align}
    \biginner{\mhi^{\epsilon}}{\phi^{\epsilon} - \overline{\phi^{\epsilon}}} & = - \biginner{\Delta \phi^{\epsilon}}{\phi^{\epsilon} - \overline{\phi^{\epsilon}}} + \intomega{(\fhi^{\epsilon})'(\phi^{\epsilon}) (\phi^{\epsilon} - \overline{\phi^{\epsilon}})}{x}{\Omega}
    \notag \\
    & \quad + \intomega{\pdif{G(\phi^{\epsilon},\psi^{\epsilon})}{\phi}(\phi^{\epsilon} - \overline{\phi^{\epsilon}}) } {x}{\Omega} + \sigma_2 \bignorm{\phi^{\epsilon} - \overline{\phi^{\epsilon}}}^2_*.
    \label{es-mu-e1}
\end{align}
Recalling the definition of $\tfhi^{\epsilon}$ and the assumption $(\mathbf{H1})$ for $\fhi$, we can apply the argument in \cite{MZ2004} and find a positive constant $C_\delta$ independent of $\epsilon$, such that
\begin{align}
    & \intomega{\bigabs{ \bigbracket{\tfhi^{\epsilon}}'(\phi^{\epsilon})}}{x}{\Omega} \le C_\delta \Bigbracket{\intomega{\bigbracket{\tfhi^{\epsilon}}'(\phi^{\epsilon}) (\phi^{\epsilon} - \overline{\phi^{\epsilon}})}{x}{\Omega} + 1}, \label{L1_tfhi}
\end{align}
and
\begin{align}
    &  c_\epsilon\intomega{\bigabs{(\flnhi)'(\phi^{\epsilon})}}{x}{\Omega} \le C_\delta \Bigbracket{c_\epsilon\intomega{(\flnhi)'(\phi^{\epsilon}) (\phi^{\epsilon} - \overline{\phi^{\epsilon}})}{x}{\Omega} + 1},
    \label{L1_lnfhi}
\end{align}
where $\delta >0$ is a sufficiently small constant such that $\overline{\phi^\epsilon}, \overline{\phi_0}, c \in (-1+\delta,1-\delta)$. The proof of \eqref{L1_lnfhi} can be found, for instance, in \cite{M2019} and we sketch the proof of \eqref{L1_tfhi} in the Appendix.
To obtain the desired estimate, we adapt the argument in \cite{ADG2013dm} with the necessary modifications. Firstly, it follows from \eqref{es-mu-e1}, the above inequalities and the estimate \eqref{unif-i} that
\begin{align*}
     c_{\epsilon} \norm{\flnhi' (\phi^{\epsilon})}_{\lp{1}}
    + \norm{(\tfhi^{\epsilon})'(\phi^{\epsilon})}_{\lp{1}}
    & \le C_\delta \bigbracket{\norm{( \mhi^{\epsilon} - \overline{\mhi^{\epsilon}}) (\phi^{\epsilon} - \overline{\phi^{\epsilon}})}_{\lp{1}} +\norm{\nabla \phi^{\epsilon}}^2 +1}\\
    & \le C \bigbracket{ \norm{\nabla \mhi^{\epsilon}} + 1 },
\end{align*}
where $C>0$ is independent of $\epsilon$. This together with \eqref{eq:nsapdmch4} yields
$$|\overline{\mu_\phi^\epsilon}|\leq C \bigbracket{ \norm{\nabla \mhi^{\epsilon}} + 1 }.$$
By the Poincar\'e--Wirtinger inequality, we further get
$$\|\mu_\phi^\epsilon\|\leq C \bigbracket{ \norm{\nabla \mhi^{\epsilon}} + 1 }.$$
Then it follows from \eqref{eq:nsapdmch4} and \eqref{unif-i} that
\begin{equation*}
\norm{(\fhi^{\epsilon})'(\phi^{\epsilon})}^2_{L^2(Q_T)} \le C \bracket{\norm{\nabla \mhi^{\epsilon}}_{L^2(Q_T)}^2 + 1}.
\end{equation*}
By construction, it holds
$$\bigbracket{\tfhi^{\epsilon}}'(s) \flnhi'(s) \ge 0,\quad \forall\,s \in (-1,1).$$
This useful fact implies that
\begin{align*}
     c_{\epsilon}^2 \norm{\flnhi'(\phi^{\epsilon})}_{L^2(Q_T)}^2
    &\leq  c_{\epsilon}^2 \norm{\flnhi'(\phi^{\epsilon})}_{L^2(Q_T)}^2
     + \norm{(\tfhi^{\epsilon})'(\phi^{\epsilon})}_{L^2(Q_T)}^2 \\
    & \quad + 2 c_{\epsilon} \inner{\flnhi'(\phi^{\epsilon})}{(\tfhi^{\epsilon})'(\phi^{\epsilon})}_{Q_T}
    \\
    &=\norm{(\fhi^{\epsilon})'(\phi^{\epsilon})}_{L^2(Q_T)}^2.
\end{align*}
Noticing that $\mbhi^{\epsilon} \ge \epsilon$, we find
\begin{align*}
    \epsilon c_{\epsilon}^2 \norm{\flnhi'(\phi^{\epsilon})}_{L^2(Q_T)}^2
     &\le C \epsilon \Bigbracket{ 1 + \intomega{ \abs{ \nabla \mhi^{\epsilon} }^2 }{x}{Q_T}} \\
    & \le C \Bigbracket{ 1 + \intomega{ \mbhi^{\epsilon}(\phi^{\epsilon}) \abs{ \nabla \mhi^{\epsilon} }^2 }{x}{Q_T}}.
\end{align*}
A similar argument gives
\begin{equation*}
    \epsilon c_{\epsilon}^2 \norm{\flnsi'(\psi^{\epsilon})}_{L^2(Q_T)}^2
    \le C \Bigbracket{1 + \intomega{ \mbsi^{\epsilon}(\psi^{\epsilon}) \abs{ \nabla \msi^{\epsilon} }^2 }{x}{Q_T}}.
\end{equation*}
Combining the above estimates with \eqref{unif-i}, we arrive at our conclusion.
\epf

\subsection{Passage to the limit as $\epsilon\to 0$.}\label{proof_existence_dm}
In this subsection, we pass to the limit as $\epsilon \to 0$ in \eqref{testapdmns}--\eqref{testapdmmsi} to recover the weak formulations in Definition \ref{wsnschdmd} for the weak solution to problem \eqref{eq:nsch1}--\eqref{eq:nschi}.

The goal is achieved by weak and strong convergence properties of the approximate solutions (always understood in the sense of a suitable subsequence).
\medskip

\bpf{\bf (Proof of Theorem \ref{wsnschdme})}
From Lemma \ref{energyestimateapdm}, we infer that for any given $T>0$, it holds\smallskip
\begin{itemize}
    \item $\sp^{\epsilon}$ is uniformly bounded in $L^{\infty} \bigbracket{0,T;\ls} \cap L^2 \bigbracket{0,T;\hssigma{1}}$,\smallskip
    \item $\jhi^{\epsilon},\jsi^{\epsilon},\widehat{\mathbf{J}}_{\phi}^{\epsilon},\widehat{\mathbf{J}}_{\psi}^{\epsilon}$ are uniformly bounded in $\bm{L}^2(Q_T)$,\smallskip
    \item $\phi^{\epsilon},\psi^{\epsilon}$ are uniformly bounded in $L^{\infty} \bigbracket{0,T;V}\cap L^2 \bigbracket{0,T;W}\cap H^1(0,T;V^*)$.
\end{itemize}
\smallskip
Then the Banach--Alaoglu theorem implies that there exist some functions $\big(\sp,\phi,\psi$, $\jhi,\jsi,\widehat{\mathbf{J}}_{\phi},\widehat{\mathbf{J}}_{\psi}\big)$ such that
\begin{align*}
    & \sp^{\epsilon} \mathop{\rightharpoonup}\limits^* \sp \quad \mathrm{in} \ L^{\infty}\bigbracket{0,T;\ls}
    \  \ \text{and}\ \
    \sp^{\epsilon} \rightharpoonup \sp \quad \mathrm{in} \ L^2\bigbracket{0,T;\hssigma{1}},
    \\
    & \phi^{\epsilon} \mathop{\rightharpoonup}\limits^* \phi, \ \ \psi^{\epsilon} \mathop{\rightharpoonup}\limits^* \psi \quad \mathrm{in} \ L^{\infty}\bigbracket{0,T;V}
    \ \ \text{and}\ \
    \phi^{\epsilon} \rightharpoonup \phi, \ \ \psi^{\epsilon} \rightharpoonup \psi \quad \mathrm{in} \ L^{2}\bigbracket{0,T;W},
    \\
    & \widehat{\mathbf{J}}_{\phi}^{\epsilon} \rightharpoonup \widehat{\mathbf{J}}_{\phi}, \ \ \widehat{\mathbf{J}}_{\psi}^{\epsilon}
    \rightharpoonup \widehat{\mathbf{J}}_{\psi},
    \ \
    \jhi^{\epsilon} \rightharpoonup \jhi, \ \ \jsi^{\epsilon} \rightharpoonup \jsi
    \quad \mathrm{in} \ \bm{L}^2(Q_T),
\end{align*}
as $\epsilon\to 0$ along a non-relabeled subsequence $\Big\{\sp^{\epsilon},\phi^{\epsilon},\psi^{\epsilon},\jhi^{\epsilon},\jsi^{\epsilon},\widehat{\mathbf{J}}_{\phi}^{\epsilon},\widehat{\mathbf{J}}_{\psi}^{\epsilon}\Big\}$.

 Hence, with a diagonal argument as in \cite{BF2013}, we can extract a further subsequence (non-relabeled for simplicity) such that there exists a limit $\bigbracket{\sp,\phi,\psi,\jhi,\jsi,\widehat{\mathbf{J}}_{\phi},\widehat{\mathbf{J}}_{\psi}}$ defined on $[0,\infty)$ in the sense that all the above weak or weak-$^*$ convergence results hold for any $T>0$.
In addition, exploiting the weak and weak-$^*$ lower semi-continuity of norms, we find
\begin{align*}
    &\sup\limits_{0\le t \le T} \Big( \intomega{ \underline{\rho} \frac{\abs{\sp(t) }^2}{2} + \half \abs{\nabla \phi (t)}^2 + \halfof{\beta} \abs{\nabla \psi (t)}^2 }{x}{\Omega} \Big) \nonumber\\
    & \qquad + \intomega{\underline{\nu} \abs{D \sp }^2 }{(x,t)}{Q_T} + \bignorm{\widehat{\bmj}_{\phi} }^2_{L^2(Q_T)} + \bignorm{\widehat{\bmj}_{\psi}}^2_{L^2(Q_T)} + \norm{\Delta \phi}^2_{L^2(Q_T)}  + \norm{\Delta \psi}^2_{L^2(Q_T)} \\
    & \quad \le \liminf\limits_{\epsilon \to 0^+} \sup\limits_{0\le t \le T} \Big( \intomega{ \underline{\rho} \frac{\abs{\sp^{\epsilon} }^2}{2} + \half \abs{\nabla \phi^{\epsilon} (t)}^2 + \halfof{\beta} \abs{\nabla \psi^{\epsilon} (t)}^2 }{x}{\Omega} \nonumber\\
    & \qquad + \halfof{\sigma_2} \bignorm{ \phi^{\epsilon}  - \overline{\phi^{\epsilon} }  }^2_* + \intomega{\fhi^{\epsilon} \bigbracket{\phi^{\epsilon} (t)} + \fsi^{\epsilon} \bigbracket{\psi^{\epsilon} (t)} }{x}{\Omega} \Big)\nonumber\\
    & \qquad + \liminf\limits_{\epsilon \to 0^+} \intomega{\underline{\nu} \abs{D \sp^{\epsilon} }^2 }{(x,t)}{Q_T} + \liminf\limits_{\epsilon \to 0^+} \bignorm{\widehat{\bmj}_{\phi}^{\epsilon} }^2_{L^2(Q_T)} + \liminf\limits_{\epsilon \to 0^+} \bignorm{\widehat{\bmj}_{\psi}^{\epsilon}}^2_{L^2(Q_T)} \\
    & \qquad + \liminf\limits_{\epsilon \to 0^+} \norm{\Delta \phi^{\epsilon}}^2_{L^2(Q_T)}  + \liminf\limits_{\epsilon \to 0^+} \norm{\Delta \psi^{\epsilon}}^2_{L^2(Q_T)} \\
    & \quad \le C,
\end{align*}
where $\underline{\rho} = \min\limits_{s \in [-1,1]} \rho(s)>0$.

For any $T>0$, by the Aubin--Lions--Simon lemma, we obtain the strong convergence
\begin{align*}
    & \phi^{\epsilon} \to \phi, \ \ \psi^{\epsilon} \to \psi \ \ \text{in} \ \ L^2 \bigbracket{0,T;V}
\end{align*}
and as a consequence,
\begin{align*}
\phi^{\epsilon} \to \phi, \ \ \psi^{\epsilon} \to \psi,\ \ \nabla \phi^{\epsilon} \to \nabla \phi, \ \ \nabla \psi^{\epsilon} \to \nabla \psi\quad \ae \ \mathrm{in} \ Q_T.
\end{align*}
The construction of approximate solutions entails that $\phi^{\epsilon} \in (-1,1)$ and $\psi^{\epsilon} \in (0,1)$ almost everywhere in $Q$, which is due to the singularity of $\fhi^{\epsilon}$ and $\fsi^{\epsilon}$.
Therefore, after passing to the limit as $\epsilon \to 0$, we can show that $\phi\in [-1,1]$, $\psi\in [0,1]$ almost everywhere in $Q$.
Hence, we have
$$\rho^\epsilon(\phi^\epsilon)\to \rho(\phi),\ \ m_\phi^\epsilon(\phi^\epsilon)\to m_\phi(\phi),\ \ m_\psi^\epsilon(\psi^\epsilon)\to m_\psi(\psi)\quad \ae \ \mathrm{in} \ Q_T.
$$
Furthermore, using $(\mathbf{H3}*)$, we can apply the method in \cite[Section 3]{B1999} to conclude
\begin{align*}
&(\widetilde{F}_{\phi}^{\epsilon})''\bracket{\phi^{\epsilon}}m_\phi^\epsilon(\phi^\epsilon) \to F_{\phi} '' \bracket{\phi} m_\phi(\phi),
\quad \ \ \, \ae \ \mathrm{in} \ Q_T,
\\
&(\widetilde{F}_{\psi}^{\epsilon})''\bracket{\psi^{\epsilon}}m_\psi^\epsilon(\psi^\epsilon) \to F_{\psi} '' \bracket{\psi} m_\psi(\psi),
\quad \ \ae \ \mathrm{in} \ Q_T.
\end{align*}
The convergence in \eqref{testapdmns}, \eqref{testapdmch1} and \eqref{testapdmch2} is similar to that in Section \ref{proofexistencendm}.
The slight difference here is that we only have the uniform boundedness of $\phi^{\epsilon}$ in $\lptbig{2}{0,T;\hs{2}}$, which leads to the uniform boundedness of $\psigma \bigbracket{ \rho^{\epsilon} \sp^{\epsilon}}$ in $\lptbig{2}{0,T;\bm{W}^{1,\frac{3}{2}}(\Omega)}$.
Thanks to the compact embedding $\wkp{1}{\frac{3}{2}} \hookrightarrow \ltwo$ in two and three dimensions, the Aubin--Lions--Simon lemma as well as the uniform boundedness of $\pdif{\psigma \bigbracket{ \rho^{\epsilon} \sp^{\epsilon}}}{t}$ in $\lptbig{\frac{8}{7}}{0,T;\wkpsigma{1}{4}^*}$, we can conclude the strong convergence
$$\psigma \bigbracket{ \rho^{\epsilon} \sp^{\epsilon}}\to \psigma \bigbracket{ \rho \sp}\quad \text{in}\ \lptbig{2}{0,T;\bm{L}^2(\Omega)}.$$
This further implies that $\sp^\epsilon\to \sp$ in $\lptbig{2}{0,T;\bm{L}^2(\Omega)}$. The continuity of $\sp$ follows from the same argument as in \cite[Section 5.2]{ADG2013ndm}.

Next, we investigate the convergence in \eqref{testapdmmhi} and \eqref{testapdmmsi}.
For any given $\bmeta \in C_c \bigbracket{[0,\infty);\bm{H}^1(\Omega)} \cap \bm{L}^{\infty}(Q)$, there exists $T\in (0,\infty)$ such that $\{ t\in \r_+: \bmeta(t) \neq 0 \} \subset [0,T]$.
Thus, for the first term and the last term on the left-hand side of \eqref{testapdmmhi}, it follows from the fact $0\le (\tfhi^{\epsilon})''\mbhi^{\epsilon} \le \alpha$ and Lebesgue's dominated convergence theorem that
\begin{align*}
    &\intomega{(\tfhi^{\epsilon})''(\phi^{\epsilon})\mbhi^{\epsilon}(\phi^{\epsilon}) \nabla \phi^{\epsilon} \cdot \bmeta}{(x,t)}{Q_T} \to \intomega{(\fhi)''(\phi)\mbhi(\phi) \nabla \phi \cdot \bmeta}{(x,t)}{Q_T},
    \\
    &\intomega{\frac{\p^2 G}{\p \phi^2}(\phi^{\epsilon},\psi^{\epsilon})\mbhi^{\epsilon}(\phi^{\epsilon}) \nabla \phi^{\epsilon} \cdot \bmeta}{(x,t)}{Q_T} \to \intomega{\frac{\p^2 G}{\p \phi^2}(\phi,\psi)\mbhi(\phi) \nabla \phi \cdot \bmeta}{(x,t)}{Q_T},
    \\
    &\intomega{\frac{\p^2 G}{\p \phi \p \psi}(\phi^{\epsilon},\psi^{\epsilon})\mbhi^{\epsilon}(\phi^{\epsilon}) \nabla \psi^{\epsilon} \cdot \bmeta}{(x,t)}{Q_T} \to \intomega{\frac{\p^2 G}{\p \phi \p \psi}(\phi,\psi) \mbhi(\phi) \nabla \psi \cdot \bmeta}{(x,t)}{Q_T},
    \\
    & \inner{ \nabla \mathcal{N} \bracket{ \phi^{\epsilon} - \overline{\phi^{\epsilon}} }}{\mbhi^{\epsilon}(\phi^{\epsilon})\bmeta}_{Q_T} \to \inner{ \nabla \mathcal{N} \bigbracket{ \phi - \overline{\phi} }}{\mbhi(\phi)\bmeta}_{Q_T}.
\end{align*}
The convergence of $\inner{\Delta \phi^{\epsilon}}{\div (\mbhi^{\epsilon}(\phi^{\epsilon}) \bmeta)}_{Q_T}$ can be treated essentially in the same way as in \cite{ADG2013dm,B1999,EG1996}, so we omit the details here. On the other hand, although the proof of the convergence
$$ \intomega{ c_{\epsilon} \flnhi'(\phi^{\epsilon}) \div (\mbhi^{\epsilon}(\phi^{\epsilon})\bmeta)}{(x,t)}{Q_T} \to 0 $$
follows the methodology in \cite{ADG2013dm}, we provide the detailed argument below since our choice of $\mbhi$ is more general. To this end, we split the integral into three parts as follows
\begin{align*}
    & \intomega{
    c_{\epsilon} \flnhi'(\phi^{\epsilon}) \div (\mbhi^{\epsilon}(\phi^{\epsilon})\bmeta)}{(x,t)}{Q_T}
    \\
    &\quad =  \intomega{ c_{\epsilon} \flnhi'(\phi^{\epsilon})  \div (\mbhi^{\epsilon}(\phi^{\epsilon})\bmeta)}{(x,t)}{\{ (x,t) \in Q_T \ : \ -1+\delta_{\phi,1}^{\epsilon} \le \phi^{\epsilon}(x,t) \le 1-\delta_{\phi,2}^{\epsilon} \}}
    \\
    &\qquad  +\intomega{ c_{\epsilon} \flnhi'(\phi^{\epsilon})  \div (\mbhi^{\epsilon}(\phi^{\epsilon})\bmeta)}{(x,t)}{\{ (x,t) \in Q_T \ : \ \phi^{\epsilon}(x,t) > 1-\delta_{\phi,2}^{\epsilon} \}}
    \\
    &\qquad  +\intomega{ c_{\epsilon} \flnhi'(\phi^{\epsilon}) \div (\mbhi^{\epsilon}(\phi^{\epsilon})\bmeta)}{(x,t)}{\{ (x,t) \in Q_T \ : \ \phi^{\epsilon}(x,t) < -1+\delta_{\phi,1}^{\epsilon} \}}.
\end{align*}
It follows directly from the definition of $c_{\epsilon}$ that the first part converges to $0$ as $\epsilon \to 0$. Concerning the second part, because of the definition of $\mbhi^{\epsilon}$, we find
\begin{align*}
    & \intomega{ c_{\epsilon} \flnhi'(\phi^{\epsilon}) \div (\mbhi^{\epsilon}(\phi^{\epsilon})\bmeta)}{(x,t)}{\{ (x,t) \in Q_T \ : \ \phi^{\epsilon}(x,t) > 1-\delta_{\phi,2}^{\epsilon} \}}
    \\
    &\quad =  \intomega{ 2 \epsilon c_{\epsilon} \flnhi'(\phi^{\epsilon}) \div \bmeta}{(x,t)}{\{ (x,t) \in Q_T \ : \ \phi^{\epsilon}(x,t) > 1-\delta_{\phi,2}^{\epsilon} \}}
    \\
    & \quad \le C c_{\epsilon} \epsilon \norm{\flnhi'(\phi^{\epsilon})}_{L^2(Q_T)}
    \\
    & \quad \le C \sqrt{ \epsilon } \bracket{ \sqrt{ \epsilon } c_{\epsilon} \norm{\flnhi'(\phi^{\epsilon})}_{L^2(Q_T)}}.
\end{align*}
Thus, the estimate obtained in Lemma \ref{energyestimateapdm}-(ii) easily entails the convergence to $0$ as $\epsilon \to 0$.
A similar argument applies to the third part as well. This completes the convergence in \eqref{testapdmmhi}. Using the same reasoning, we can establish the convergence in \eqref{testapdmmsi}. Finally, the following identities
$$\jhi = \sqrt{\mbhi(\phi)} \widehat{\mathbf{J}}_{\phi}\quad\text{and}\quad \jsi = \sqrt{\mbsi(\psi)} \widehat{\mathbf{J}}_{\psi}$$
can be verified as in \cite{ADG2013dm}.

We now turn to the special case with $\sigma_1 \equiv 0$ or $\overline{\phi_0} = c$ and prove the energy inequality \eqref{energyinequalitydm} using the argument at the end of Section \ref{proofexistencendm}.

In this case, the energy inequality for the approximate solutions reduces to the following form:
    \begin{align*}
        & \eto^{\epsilon} \bigbracket{ \sp^{\epsilon}(t), \phi^{\epsilon}(t), \psi^{\epsilon}(t) } \\
        &\qquad + \intomega{\nu \bracket{\phi^{\epsilon}} \abs{D \sp^{\epsilon}}^2 + \mbhi^{\epsilon}\bracket{ \phi^{\epsilon} } \bigabs{\nabla \mhi^{\epsilon}}^2 + \mbsi^{\epsilon}\bracket{ \psi^{\epsilon} } \bigabs{\nabla \msi^{\epsilon}}^2  }{(x,\tau)}{Q_{(s,t)}} \\
        & \quad \le \eto^{\epsilon} \bigbracket{ \sp^{\epsilon}(s), \phi^{\epsilon}(s), \psi^{\epsilon}(s) },
    \end{align*}
for almost all $s \in [0,\infty)$ (including $s=0$) and all $t \in [s,\infty)$.
The only difference is that we need the weak convergence $\mathbf{\widehat{J}}_{\phi}^{\epsilon} \rightharpoonup \mathbf{\widehat{J}}_{\phi}$ and $\mathbf{\widehat{J}}_{\psi}^{\epsilon} \rightharpoonup \mathbf{\widehat{J}}_{\psi}$ to conclude
\begin{align*}
    & \intomega{ \bigabs{\mathbf{\widehat{J}}_{\phi}}^2 }{(x,\tau)}{Q_{(s,t)}} \le \liminf\limits_{\epsilon \to 0} \intomega{ \mbhi^{\epsilon} \bracket{\phi^{\epsilon}} \abs{ \nabla \mhi^{\epsilon} }^2 }{(x,\tau)}{Q_{(s,t)}}, \\
    & \intomega{ \bigabs{\mathbf{\widehat{J}}_{\psi}}^2 }{(x,\tau)}{Q_{(s,t)}} \le \liminf\limits_{\epsilon \to 0} \intomega{ \mbsi^{\epsilon} \bracket{\psi^{\epsilon}} \abs{ \nabla \msi^{\epsilon} }^2 }{(x,\tau)}{Q_{(s,t)}}.
\end{align*}
The proof of Theorem \ref{wsnschdme} is complete.
\epf
\medskip

\appendix
\section{Existence of a Pressure Field}\label{existence_pressure}
Regarding the regularity of $\rho \bracket{\phi} \sp$ in time as well as the existence of a pressure field, we have the following result:
\bpp\label{recover_pressure}
    Let $\bracket{\sp,\phi,\psi,\mhi,\msi}$ (resp. $\bracket{\sp,\phi,\psi,\jhi,\jsi}$) be a weak solution to the problem \eqref{eq:nsch1}--\eqref{eq:nschi} in the case of nondegenerate (resp. degenerate) mobilities in the sense of Definition \ref{wsnschndmd} (resp. Definition \ref{wsnschdmd}). We have
    $$ \pdif{\psigma \bracket{\rho \sp}}{t} \in \lptbigloc{\frac{8}{7}}{[0,\infty);\wkpsigma{1}{4}^*}. $$
    Moreover, there exists a unique pressure $\pi \in \wkptbigloc{-1}{\infty}{[0,\infty);\lp{\frac{4}{3}}}$ with $\overline{\pi} = 0$ (corresponding to the given weak solution) such that it holds
    $$ \pdif{ \bracket{\rho \sp} }{t} + \div \bracket{ \sp \otimes \bracket{\rho \sp - \gamma \jhi} } - \div \bracket{\nu\bracket{\phi} D \sp} + \nabla \pi = \bm{f}, $$
    in the distribution sense, where
    \begin{equation} \notag
        \bm{f}=\left\{
        \begin{aligned}
            & \mhi \nabla \phi + \msi \nabla \psi, \quad && \mathrm{if} \ \mbhi \ \mathrm{and} \ \mbsi \ \mathrm{nondegenerate}, \\
            & -\Delta \phi \nabla \phi - \beta \Delta \psi \nabla \psi + \sigma_2 \mathcal{N} \bracket{\phi-\overline{\phi}}, \quad && \mathrm{if} \ \mbhi \ \mathrm{and} \ \mbsi \ \mathrm{degenerate}.
        \end{aligned}
        \right.
    \end{equation}
\epp

\bpf Firstly, it follows from the Definition \ref{wsnschndmd} that
$$ \bm{f} \in \lptbigloc{2}{[0,\infty);\bm{L}^{\frac{3}{2}}(\Omega)} \subset \lptbigloc{\frac{8}{7}}{[0,\infty);\bm{L}^{\frac{4}{3}}(\Omega)}, $$
in the case of nondegenerate mobilities.
Similarly, in the case of degenerate mobilities, it follows from the Definition \ref{wsnschdmd} that
$$ \bm{f} \in \lptbigloc{1}{[0,\infty);\bm{L}^{\frac{3}{2}}(\Omega)} \cap \lptbigloc{2}{[0,\infty);\bm{L}^{1}(\Omega)}, $$
which implies
$$ \bm{f} \in \lptbigloc{\frac{8}{7}}{[0,\infty);\bm{L}^{\frac{4}{3}}(\Omega)}. $$
Therefore, in both cases, it holds
\begin{align*}
    \inner{\pdif{\psigma \bracket{\rho \sp}}{t}}{\bmtheta}_{Q_T} & = - \inner{\rho \sp}{\pdif{\bmtheta}{t}}_{Q_T} \\
    & = \inner{\rho \sp \otimes \sp}{\nabla \bmtheta}_{Q_T} - \inner{\nu(\phi) D \sp}{D \bmtheta }_{Q_T} \\
    & \quad - \gamma \biginner{\bracket{\sp \otimes \jhi}}{\nabla \bmtheta }_{Q_T} + \inner{f}{\bmtheta}_{Q_T} \\
    & \le C \norm{\bmtheta}_{\lpt{8}{0,T;\wkp{1}{4}}},
\end{align*}
for all $T>0$ and $\bmtheta \in C_0^{\infty}\bracket{Q_T}$ with $\div \bmtheta = 0$,
where we have used the following estimates:
\begin{align*}
    \rho \sp \otimes \sp \ &\mathrm{in} \ L^2 \bigbracket{0,T;\bm{L}^{\frac{3}{2}}(\Omega)}, \\
    D \sp \ &\mathrm{in} \ \lptbig{2}{0,T;\bm{L}^2(\Omega)},\\
    \sp \otimes \jhi \ & \mathrm{in} \ \lptbig{\frac{8}{7}}{0,T;\bm{L}^{\frac{4}{3}}(\Omega)}.
\end{align*}
It then follows from the density of $\bmtheta \in C_0^{\infty}\bracket{Q_T}^d$ with $\div \bmtheta = 0$ in $\lptbig{8}{0,T;\wkpsigma{1}{4}}$ that
$$ \pdif{\psigma \bracket{\rho \sp}}{t} \in \lptbigloc{\frac{8}{7}}{[0,\infty);\wkpsigma{1}{4}^*}. $$
On the other hand, since it holds $\sp \in C_w \bigbracket{[0,\infty);\ls}$ and $\rho \bracket{\phi} \in C\bigbracket{[0,\infty);H} \cap L^{\infty} \bracket{Q}$, we have
$$ \rho \sp \in C_w \bigbracket{[0,\infty);\bm{L}^2(\Omega)}. $$
Now let us define
$$ \bm{G}(t) = \rho\bracket{\phi(t)} \sp(t) + \intinterval{ \bigbracket{ \div \bracket{ \sp \otimes \bracket{\rho \sp - \gamma \jhi } } - \div \bracket{\nu \bracket{\phi} D \sp} - \bm{f} } }{\tau}{0}{t}. $$
It follows from the definition of weak solutions to \eqref{eq:nsch1}-\eqref{eq:nschi} that
\begin{equation*}
\pairing{\bm{G}(t)}{\bmtheta}_{\hs{-1},\hszero{1}} = 0,\quad \forall\, t \ge 0\ \ \text{and}\ \ \bmtheta \in C^{\infty}_{0,\sigma} \bracket{\Omega}.
\end{equation*}
Thanks to the foregoing estimates, it holds $\bm{G} \in C_w \bigbracket{[0,\infty);\bm{W}^{-1,\frac{4}{3}}(\Omega)}.$ Then, the de Rham theorem (see \cite[Theorem IV.2.3]{BF2013}) guarantees that, for all $t \ge 0$, there exists a unique $p (t) \in \lp{\frac{4}{3}}$ with $\overline{p} = 0$ such that
\begin{equation} \label{grad_p_2}
    \bm{G}(t) = - \nabla p(t).
\end{equation}
In addition, as in \cite{BF2013}, it can be proven that $p \in C_w \bigbracket{[0,\infty);\lp{\frac{4}{3}}}$. Thus, it holds
$$p \in \lptbig{\infty}{0,T;\lp{\frac{4}{3}}},\quad \forall\, T >0.$$
We now introduce the distribution $\pi = \pdif{p}{t}$, so it holds $\pi \in \wkptbig{-1}{\infty}{0,T;\lp{\frac{4}{3}}}$ for all $T >0$.
By taking test functions of the form $\pdif{\bmtheta}{t}$ with $\bmtheta \in C_0^{\infty}\bigbracket{\Omega \times (0,\infty)}^d$ in \eqref{grad_p_2} and recalling the definition of $\bm{G}$, we can show that it holds
\begin{equation*}
    \pdif{\bracket{\rho \sp}}{t} + \div \bracket{ \sp \otimes \bracket{\rho \sp - \gamma \jhi } } - \div \bracket{\nu \bracket{\phi} D \sp} + \nabla \pi = \bm{f},
\end{equation*}
in the sense of distributions in $Q$.
\epf
\medskip

\section{Estimate for the Approximate Potential Function}\label{estimate_approximate_potential}
Below we prove the estimate \eqref{L1_tfhi} using the argument in \cite[Appendix]{MZ2004} (see also \cite{M2019}).
\medskip
\bl\label{uniform_bound_approximate_potential}
Suppose that the assumptions (\textbf{H1}) and (\textbf{H3*}) hold except the requirements on the continuity of $\fhi$ and $\fsi$ at the corresponding pure phases, that is, here it is allowed that
\begin{align*}
& \fhi(s_1) \to + \infty, \ \ \ \text{when} \ s_1 \to -1^+ \ \text{or} \ s_1 \to 1^-,
\\
& \fsi(s_2) \to + \infty, \ \ \ \text{when} \ s_2 \to 0^+ \ \text{or} \ s_2 \to 1^-.
\end{align*}
Let $\tfhi^\epsilon$ and $\tfsi^\epsilon$ be defined as in Section \ref{approximate_degenerate}.
For any given $m_1 \in (-1,1)$ (resp. $m_2 \in (0,1)$), there exist two positive constants $c_{m_1}$ and $c_{m_1}'$ (resp. $c_{m_2}$ and $c_{m_2}'$) independent of $\epsilon$ such that
$$ c_{m_1} \bigabs{\bracket{\tfhi^{\epsilon}}'(s_1)} - c_{m_1}' \le \bracket{\tfhi^{\epsilon}}'(s_1) (s_1-m_1), \quad \forall\, s_1 \in (-1,1)$$
and
$$ c_{m_2} \bigabs{\bracket{\tfsi^{\epsilon}}'(s_2)} - c_{m_2}' \le \bracket{\tfsi^{\epsilon}}'(s_2) (s_2-m_2),\quad \forall\, s_2 \in (0,1).
$$
Moreover, $c_{m_1}$ and $c_{m_1}'$ (resp. $c_{m_2}$ and $c_{m_2}'$) continuously depend on $m_1$ (resp. $m_2$).
\el
\medskip

\bpf Let us first introduce some notation:
\begin{align*}
& \Delta_{\phi,1}^{\epsilon} \stackrel{\rm{def}}{=} \sup \set{s \in(0,1):\mbhi(-1+s) \le \epsilon}, \\
& \Delta_{\phi,2}^{\epsilon} \stackrel{\rm{def}}{=} \sup \set{s \in(0,1):\mbhi(1-s) \le \epsilon}.
\end{align*}
Obviously, as $\epsilon \to 0$, $\Delta_{\phi,1}^{\epsilon}$ and $\Delta_{\phi,2}^{\epsilon}$ strictly decrease and tend to $0$.
In addition, for all
$$s \in I_1^{\epsilon} \stackrel{\rm{def}}{=} [-1 + \Delta_{\phi,1}^{\epsilon},1-\Delta_{\phi,2}^{\epsilon}],$$
it holds $\mbhi(s) \ge \epsilon$, thus we have
$$ \half \fhi''(s) \le \bigbracket{\tfhi^{\epsilon}}''(s) \le \fhi''(s). $$
Recall that for $\fhi$ and any given $m\in(-1,1)$, if we take $c_m = \min \big\{\halfof{1+m},\halfof{1-m} \big\}$, then there exists a constant $c_m'$ continuously depending on $m$ (see \cite[Proposition 4.3]{M2019}) such that
$$ c_m \bigabs{\fhi'(s)} - c_m' \le \fhi'(s) (s-m), \quad \forall\, s \in (-1,1).
$$
Therefore, for all $s \in I_1^{\epsilon}$, we have
\begin{align*}
    c_m \bigabs{\bigbracket{\tfhi^{\epsilon}}'(s)} - c_m' & \le c_m \bigabs{\fhi'(s)} - c_m' \le \fhi'(s) (s-m) \\
    & \le \max \left\{ \bigbracket{\tfhi^{\epsilon}}'(s) (s-m),\ \half \bigbracket{\tfhi^{\epsilon}}'(s) (s-m) \right\},
\end{align*}
which implies
$$ c_m \bigabs{\bigbracket{\tfhi^{\epsilon}}'(s)} - 2 c_m' \le \bigbracket{\tfhi^{\epsilon}}'(s) (s-m).$$
Note that there exists $\epsilon_1 \in (0,1]$ such that for all $\epsilon \in (0,\epsilon_1]$, it holds $m + c_m, m - c_m \in I_1^{\epsilon}$.
Hence, for all $\epsilon \in (0,\epsilon_1]$ and $s \in [-1,1] \bigcap \bracket{I_1^{\epsilon}}^c$, we have
\begin{align*}
    & \bigbracket{\tfhi^{\epsilon}}'(s) (s - m - c_m) \ge 0,\quad  \bigbracket{\tfhi^{\epsilon}}'(s) (s - m + c_m) \ge 0,
\end{align*}
which yields
$$
c_m \bigabs{\bigbracket{\tfhi^{\epsilon}}'(s)} \le \bigbracket{\tfhi^{\epsilon}}'(s) (s-m).
$$
On the other hand, for all $\epsilon \ge \epsilon_1$ and $s \in (-1,1)$, we have
$$ \bigbracket{\tfhi''}^{\epsilon}(s) \le \frac{1}{\epsilon_1}\fhi''(s)\mbhi(s) \le \frac{\alpha}{\epsilon_1},$$
which implies
$$ c_m \bigabs{\bigbracket{\tfhi^{\epsilon}}'(s)} - \bigbracket{\tfhi^{\epsilon}}'(s) (s-m) \le (c_m + 2) \frac{\alpha}{\epsilon_1}. $$
Combining the above estimates, we obtain
$$ c_m \bigabs{\bigbracket{\tfhi^{\epsilon}}'(s)} - \left( 2 c_m' + (c_m + 2) \frac{\alpha}{\epsilon_1} \right)  \le \bigbracket{\tfhi^{\epsilon}}'(s) (s-m), $$
for all $\epsilon \in \big(0,\half \mbhi(0)\big] \bigcup \big(0,\half \mbsi(\half)\big]$ and all $s \in [-1,1]$.
The estimate for $\bigbracket{\tfsi^{\epsilon}}'$ is almost the same and we just omit it here.
\epf
\bigskip

\noindent
\textbf{Declarations}
\medskip
\\
\noindent
\textbf{Conflict of interest.} The authors have no competing interests to declare that are relevant to the content of this article.
\smallskip
\\
\textbf{Fundings.} The research of H. Wu was partially supported by Natural Science Foundation of Shanghai (No. 25ZR1401023).
\smallskip
\\
\noindent
\textbf{Data availability.} Data sharing not applicable to this article as no datasets were generated or analysed during the current study.
\smallskip
\\
\noindent
\textbf{Acknowledgments.} The authors are grateful to the anonymous referee for several valuable comments and suggestions.
M. Grasselli is a member of Gruppo Nazionale per l'Ana\-li\-si Matematica, la Probabilit\`{a} e le loro Applicazioni (GNAMPA), Istituto Nazionale di Alta Matematica (INdAM). His research is part of the activities of ``Dipartimento di Eccellenza 2023--2027'' of Politecnico di Milano.
H. Wu is a member of Key Laboratory of Mathematics for Nonlinear Sciences (Fudan University), Ministry of Education of China.


\end{document}